%
%
\magnification\magstep1 \input amstex \input amsppt.sty


\font\bbf=cmbx10 scaled \magstep1   
\def\cal#1{{\Cal#1}}
\def\O{{}^{}\Cal O}
\def\Univ{\hbox{\font\=cmssbx10\U}} 
\def\lei{\hbox{\kern.45mm$_{^\downarrow}\kern-1.280mm\cap\kern.85mm$}}
\def\Ze{Z\!\!\!Z} 
\def\Zepp{{{{Z\!\!\!Z^{\phantom{l}}}^{{}_{{}^{\!}\!+}}}}}  
\def\Zep{{{{{{{{Z\!\!\!Z_{}}_{}}_{}}_{}}}_{{}^{\!+}}}}}  
\def\supp{\roman{supp}{}^{}\,}
\def\Card{{}^{}\roman{Card}{}^{}\,}
\def\inve{\lower.85mm\hbox{$^{^-}$}\kern-.5mm{}^\iota}
\def\Eps{\hbox{\font\=cmmi10 scaled\magstep1\\char'017}\kern0.15mm}
\def\Nu{\hbox{\font\=cmmi10 scaled\magstep1\\char'027}\kern0.15mm}
\def\rajou{{}^{}\Cal B_{s\,}} 
\def\value{\hbox{\kern.2mm\font\=cmr10\\char'022\kern-.2mm}} 
\def\image{\hbox{\font\=cmr10\\char'022\kern-1mm\char'022}} 
\def\rist{\,{\times}_{{}_{{}_{\!\!\!\!\bold{t}\ }}}}  
\def\risti#1{{}^{}\,{\times}_{{\!}_{{}_{\!\!\!{#1}{}^{\!}\ }}}}
\def\tanB{{\tau_{_B}}} 
\def\Circ{\,\text{\font\=cmbsy10\\char'016}\sp\,}
\def\Examplee{{\font\=cmssi10\E\kern.15mmx\kern.15mma\kern.15mmm\kern.14mmp\kern.17mml\kern.15mme}\kern.3mm. }
\def\Examples{{\font\=cmssi10\E\kern.15mmx\kern.15mma\kern.15mmm\kern.14mmp\kern.17mml\kern.15mme\kern.15mms}\kern.3mm. }
\def\Remarkk{{\font\=cmssi10\R\kern.15mme\kern.15mmm\kern.15mma\kern.15mmr\kern.15mmk\kern.15mm. }}
\def\Remarkss{{\font\=cmssi10\R\kern.15mme\kern.15mmm\kern.15mma\kern.15mmr\kern.15mmk\kern.15mms}\kern.3mm. }
\def\N{{I\!\!N}} 
\def\No{{I\!\!N\kern-.54mm\lower.15mm\hbox{$_{\roman o}$}}} 
\def\Nopot#1{I\!\!N\kern-.54mm\lower.15mm\hbox{$_{\roman o}$}\kern-.7mm{}^{#1}} 
\def\potNo{^{\kern.37mm I\!\!{N_{}}_{\kern-.22mm\roman o}}} 
\def\plusinfty{\lower1.05mm\hbox{$^+$}\infty}
\def\Re{{I\!\!R}}
\def\Rep{{{I\!\!R^{\phantom{l}}}^{{}_{{}^{\!}\!+}}}}
\def\Repp{{{{{{{I\!\!R_{}}_{}}_{}}_{}}}_{{}^{\!+}}}}
\def\Ce{{\hbox{$C\kern-2.5mm\raise.2mm\hbox{\font\=cmssqi8\I}\kern1.48mm$}}}
\def\imag{\kern.15mm\lower.6mm\hbox{$^{^*}$}\kern-1.8mm\imath\kern.1mm} 
\def\Ke{{{}^{}I\!\!K}}
\def\bit#1{{\font\†=cmmib10\text{\†\char'#1}}} 
\def\biit#1{\hbox{\font\=cmmib10\#1}} 
\def\bmi#1#2{\hbox{\font\=cmmib#1\\char'#2}}
\def\eeta{021}\def\xxi{030}\def\varPii{\char'005}
\def\fssi#1{\hbox{\font\=cmssi10\#1}\kern0.15mm} 
\def\smb#1{\hbox{\font\†=cmmi8\†#1\kern.3mm}} 
\def\ssmb#1{\hbox{\font\=cmmi6\#1}} 
\def\ecal#1{\kern.1mm\hbox{\font\†=cmsy8\†#1\kern.3mm}} 
\def\conc{\!\bold{\hat{\phantom w}}\!}
\def\concc{{}^{}{}^{}\!\bold{\hat{\phantom w}}\!{}^{}}
\def\idv{\hbox{\font\=cmr10\id}\kern.25mm\lower.8mm\hbox{\font\=cmr7\v}\kern.3mm} 
\def\seq#1{\langle#1\rangle}
\def\ymp{{}^{}\Cal N_o\,}
\def\vecs{\upsilon\kern-0.3mm\lower.15mm\hbox{$_s$}\kern0.2mm} 
\def\vecss{\hbox{\font\=cmitt10\v}\kern-0.1mm\lower.15mm\hbox{$_s$}\kern0.2mm} 
\def\bnull#1{\hbox{\font\=cmssbx10\0}{}_{\font\=cmmi6\lower.15mm\hbox{\kern-.1mm\#1\kern.15mm}}} 
\def\dom{{{}^{}\roman{dom}\,{}_{{}^{}}}}

\def\rng{{}^{}\roman{rng}\,{}_{{}^{}}}
\def\CPi#1{C\kern-.2mm\lower.05mm\hbox{$_{_\Pi}$}\kern-1.52mm{}^{#1}}
\def\Lip#1{{}^{}{{{{{\Cal L{}^{}ip}^{}}^{}}^{}}^{}}^{#1{}^{}}         }
\def\Cinfty{C\kern.4mm\raise.3mm\hbox{$^\infty$}\kern.15mm}
\def\Cinftyzero{\hbox{$C\kern.4mm\raise.3mm\hbox{$^\infty$}\kern.15mm\kern-3.5mm_{\font\=cmr6\lower.15mm\hbox{\kern.1mm\0}}\kern1.9mm$}}
\def\RHB#1#2{\raise#1mm\hbox{$#2$}} 
\def\LHB#1#2{\lower#1mm\hbox{$#2$}} 
\def\sixroman#1{\hbox{\font\=cmr6\#1\kern.1mm}}
\def\sNor#1{\kern.25mm\lower.38mm\hbox{$_{#1}$}}
\def\sNorr#1{\kern-.2mm\lower.38mm\hbox{$_{#1}$}}
\def\ai#1{{}_{\font\=cmmi6\lower.15mm\hbox{\kern-.1mm\#1\kern.15mm}}} 
\def\yi#1{^{\font\=cmmi6\raise.0mm\hbox{\kern-.1mm\#1\kern.15mm}}} 
\def\ar#1{{}_{\font\=cmr6\lower.15mm\hbox{\kern.1mm\#1}}} 
\def\yr#1{^{\font\=cmr6\raise.0mm\hbox{\kern.3mm\#1}}} 
\def\yrai^#1_#2{^{\kern.4mm\hbox{\font\=cmr6\{#1}}}_{\kern.2mm{#2}}}
\def\lupar{\kern.2mm\lower1mm\hbox{$^{^(}$}} 
\def\rupar{\lower1mm\hbox{$^{^)}$}\kern-.15mm} 
\def\yyi#1{^{\font\=cmmi6\lower.6mm\hbox{\kern-.25mm\#1\kern-.05mm}}} 
\def\yyr#1{^{\font\=cmr6\lower.45mm\hbox{\kern-.25mm\#1\kern-.15mm}}} 
\def\yplus{\lower1mm\hbox{$^{^+}$}} 
\def\yminus{\lower1mm\hbox{$^{^-}$}} 
\def\ClT{\roman{Cl}\kern.25mm\lower.4mm\hbox{$_{\Cal T}$}\kern0.2mm} 
\def\IntT{\sp\roman{Int}\kern.2mm\lower.4mm\hbox{$_{\Cal T}$}\kern0.2mm} 
\def\inc{\subseteq}
\def\iinc{\supseteq}
\def\exi#1{\exists\,#1\kern.2mm\,;}
\def\all#1{\forall\,#1\kern.2mm\,;}
\def\imply{\Rightarrow}
\def\equivv{\Leftrightarrow}
\def\embed{\hookrightarrow}
\def\sp{\kern0.15mm} 
\def\ssp{\kern0.37mm} 
\def\snn{\kern-0.2mm} 
\def\sn{\kern-0.3mm} 
\def\ssn{\kern-0.63mm} 
\def\KP#1{\kern#1mm} 
\def\KN#1{\kern-#1mm} 
\def\hyppy#1{$\phantom{}$\hskip#1}
\def\mhyppy#1{\text{\hskip#1mm}}
\def\NS{\vskip1.7mm}
\def\NSN{\vskip1.7mm\noindent}
\def\VBOX/#1/#2/HEREend{\vbox{#2\vskip-#1mm}\vfill\null\eject}
\def\œ$#1${\hbox{$#1$}} 
\def\"{\"a} \def\"{\"o}
\def\newProCla#1\par#2\par{\vskip1.7mm\noindent\bf#1\it#2\vskip1.7mm}
\def\Prooff{{\font\=cmssi10\P\kern.37mmr\kern.37mmo\kern.37mmo\kern.37mmf\kern.37mm. }\rm}

\def\QED{\hfill\hbox{$\ \sqcap\kern-2.45mm\sqcup$}}
\def\newQED{\hfill\hbox{$\ \sqcap$\hskip-2.45mm$\sqcup$}\vskip1.7mm}
\def\abstract#1{{\eightpoint\parindent5mm\narrower\baselineskip3.36mm#1\par}}
\def\noin{\noindent}
\def\Newline{\kern-10mm\newline}
\font\rp=cmr8
\def\eps{\varepsilon}
\def\RunMyHead#1#2#3#4{%
 \headline{\ifnum\pageno=\firstpage\hfil%
           \else{\ifodd\pageno{\rp#3\phantom\folio\hfil#4\hfil\phantom{#3}\folio}%
                 \else{\rp\folio\phantom{#2}\hfil#1\hfil\phantom\folio#2}%
                 \fi}%
           \fi}%
 \footline{\ifnum\pageno=\firstpage\hfil{\rp[\,\folio\,]}\hfil%
           \else\hfil%
           \fi}%
}%
\def\subhead#1\par#2\par{\vskip4mm\smallbreak\null\smallskip\vbox{\noindent\bbf#1\hfill\kern1.5mm#2\hfill\phantom{#1}\vskip2.5mm\nopagebreak}\nopagebreak\noindent}
\def\subheadd#1\par#2\par#3\par{\vskip4mm\smallbreak\null\smallskip\vbox{\noindent\bbf#1\hfill#2\hfill\phantom{#1}\vskip1.5mm\centerline{#3}\vskip2.5mm\nopagebreak}\nopagebreak\noindent}
\def\CinftyPi{\Cinfty\kern-3.5mm_{_{\bold\Pi}}\kern1.45mm}
\def\CinftyS{\Cinfty\kern-3.9mm_{_{\Cal S}}\kern1.45mm}
\def\wave{\hbox{\font\†=cmsy10\†\hbox{\char'164}\kern-2.35mm\hbox{\char'165}\kern.4mm}}
\def\wavee{\hbox{\font\†=cmsy8\†\hbox{\char'164}\kern-2.0mm\hbox{\char'165}\kern.4mm}} 
\def\barmj{\kern.25mm\bar{\hbox{\font\=cmr10\\char'021}}\kern.4mm}

\def\fssit{\font\=cmssi10\}
\def\sigrd{\sigma\kern-.3mm_{_{rd}}\kern.15mm} 
\def\taurd{\tau_{_{rd}}\kern.15mm}
\def\tsigrd{\tau\sigma\kern-.3mm_{_{rd}}\kern.15mm} 
\def\tauR#1{\tau_{_{I\!\!R}}\kern-1.5mm^{#1}}
\def\RN{I\!\!R\kern.3mm^{\hbox{\font\=cmmi6\N}}} 
\def\QTN{Q\kern.1mm_{\lower.2mm\hbox{\font\=cmmi6\T}}^{\kern.2mm\hbox{\font\=cmmi6\N}}} 
\def\leLCSr{{\le}{}_{_{\roman{LCS}}}}
\def\leLCS-{{\le}{}_{_{\roman{LCS}}}\text{\sp-\sp}}

\def\Centerline#1\par#2\par#3{\noindent#1\phantom{#3}\hfill#2\hfill\phantom{#1}#3}
\def\varPii{\char'005}


\hsize125mm \vsize240truemm \parskip.5mm \overfullrule=0mm \hoffset30truemm
\vsize243truemm \voffset23truemm 
\vsize243truemm \voffset8truemm
\hcorrection{-1truein}\vcorrection{-1truein}

\document \baselineskip4.5mm \newcount\firstpage\firstpage=1\pageno=\firstpage

\RunMyHead{S.\ Hiltunen}{Differentiation,}{implicit functions, and applications $\ldots$}{}

{\font\=cmss12 scaled\magstep1\
\centerline{%
              Differentiation, implicit functions, and}\vskip1mm\centerline{%
             applications to generalized well-posedness}}\vskip3mm\centerline{\font\=cmr12\%
                               by}\vskip2mm\centerline{\font\=cmmi12\%
                      Seppo\kern1mm Hiltunen}\vskip5mm\noin
\abstract{{\bf Abstract.} This article is centered around generalizing a
  previous implicit function theorem of the author to be applicable for maps $
f:E\sqcap F\to F$ which can be lifted to Keller $C^k_\pi\,$--\,maps $f_i:E
\sqcap F_i\to F_i$ with $F_i$ Banach and $F=\projlim F_i\ssp$. \hfill We prove
theorems about existence and differentiability of functions $g$ satisfying $
f\sp(\sp x\ssp,g\ssp(x))=b$ constant. In addition to these abstract theorems,
we give several examples of applications to proving smooth dependence of the
solution on initial/boundary values and the nonlinearity in nonlinear
(partial) differential equations.\vskip1mm
\noin{\bf%
Keywords:} implicit functions, smooth solution maps, well-posedness,
  projective limits, locally convex spaces, Banach spaces, differentiability,
smoothness, holomorphy.\vskip1mm
\noin{\bf%
Subject classification:} 58C15, 47J07, 46T20, 46T25, 35B30. (AMS 2000)}


\subhead 0

                                Introduction

On \hfil $Q\ai I=I\sn\times\sn\Re\ssp\yi N$ \hfil with \hfil $I=\Repp$ \hfil
or \hfil $I=[\,0\,,\smb T\sp\,]$ \hfil for \hfil $\smb T\in\Rep,$ \hfil
consider a nonlinear\linebreak partial differential equation \hfil $
\partial\ar 0\sp y=P\ssp y+\varphi\circ[\,\roman{id\,},R\ssp y\,]$ \hfil with
initial condition\linebreak $y\ssp(\sp 0\,,\cdot)=y\ar 0$ for functions $y:
Q\ai I\to Y\sp$. Here $Y$ is a finite dimensional real vector space, and $
P\sp,R$ are linear partial differential operators not containing "time\ssp"
derivatives, and we define $[\,\roman{id\,},z\,]\ssp(\xi)=(\sp\xi\ssp,
z\ssp(\xi))$ for $\xi=(\sp t\ssp,\eta\sp)\in Q\ai I\ssp$. For example, the
nonlinear scalar wave equation $\,\wave u={u_{}}_{tt}-\Delta\ssp u=\psi\ssp(\sp
t\ssp,\eta\,;u\ssp,{u_{}}_t\sp,\roman{grad\,}u\sp)$ with $u\ssp(\sp 0\,,\cdot)
=u\ar 0$ and ${u_{}}_t\sp(\sp 0\,,\cdot)=u\ar 1$ can be put in the previous
form by defining $P\ssp y=$ $[\,\bnull{}\ssp,\Delta\ssp u\,]$ and $R\ssp y=
[\,u\ssp,v\ssp,\roman{grad\,}u\,]\ssp,$ when we have $y=[\,u\ssp,v\,]\ssp,$
and $\varphi\ssp(\sp t\ssp,\eta\,;\xi\ar 1\sp,\xi\ar 2\ssp,\zeta\sp)$ $=(\sp
\xi\ar 2\ssp,\psi\ssp(\sp t\ssp,\eta\,;\xi\ar 1\sp,\xi\ar 2\ssp,\zeta\sp))\ssp$.

{\it Well-posedness\ssp} of the above initial value problem traditionally
refers to a result of the type that for each initial value $y\ar 0$ in a set $
O$ in some function space $E\ar 0$ there is a unique solution $y$ in a space $
F=F\ssn\ai I$ of functions $Q\ai I\to Y\sp$. One might require $O$ to be open
in $E\ar 0\ssp$, and moreover, that the map $g:E\ar 0\to F$ given by $O\owns y
\ar 0\mapsto y$ is continuous. Instead of continuity of this {\it solution map\ssp}
$g\ssp,$ one might be interested to know whether it is of class $C^{\ssp 1}$
or smooth or analytic.

The choice of the space $F$ depends on $E\ar 0$ and the conditions imposed on
$\varphi\ssp$. If $\varphi$ is smooth, one might be able to prove {\it
regularity\ssp} results saying that for example the solutions $y$ are smooth,
and one then wishes to choose as $F$ some space of smooth functions. In the
articles [\,Po1\,] and [\,Po2\,]\ssp, Poppenberg considers nonlinear
Schr\"dinger equations $\imag\,\partial\ar 0\ssp y=\varphi\circ
[\,\roman{id\,};y\ssp,\roman{grad\,}y\ssp,\Delta\ssp y\,]$ with $y\ssp(\sp
0\,,\cdot)=y\ar 0\ssp$. Using inverse and implicit function theorems of
Nash\,--\,Moser type, he obtains existence of {\it continuous\ssp} and $
C^{\ssp 1}$ solution maps $y\ar 0\mapsto y\ssp$. Poppenberg's most strict
space $F$ is isomorphic to $C^{\ssp 1\sp}(\sp I\sp,E\ar 0)$ which is obtained
in [\,Po1\,] with $E\ar 0=H^\infty(\Re\sp)\ssp$.

In [\,PS\ssp; Thm.\ 1, p.\ 171\,]\ssp, the authors sketch a "local\ssp"
well-posedness assertion about a $\Lip 0$ solution map \hfil $g:E\ar 0\iinc O
\to F$ \hfil of Banach spaces for a scalar wave\linebreak equation \hfil $
\wave u=u^{\,k}\cdot\prod{}\yi N\KN{2.5}\LHB{.5}{_{i=0\,}}(\sp\partial_iu\sp
)^{\,\alpha_i}$ \hfil obviously required to satisfy by $u$ only in some\linebreak
weak sense. In [\,KR\,]\ssp, no continuity of the solution map $g$ is even
asserted. Further parallel theorems about at most continuous normable space
solution maps for nonlinear wave equations can be found in [\,Ge\,]\ssp,
[\,KS\,]\ssp, [\,Z\,]\ssp, to give some examples for background. These
well-posedness assertions with solution maps $g$ are roughly of the kind where
$E\ar 0\cong H^{\sp\sigma}(\Re\ssp\yi N)\sqcap\cdots\,H^{\,\sigma\sp-\sp i\sp}
(\Re\ssp\yi N)$ with $\sigma<\infty$ and $i\in\No\ssp,$ and $F$ isomorphic to
a subspace of $C\ssp(\sp I\sp,E\ar 0)\ssp$. A tendency in research has been to
find minimal $\sigma\ssp$, that is, as loose initial value space $E\ar 0$ as possible.

We want to go in the other direction and see what can be said about a solution
map if $E\ar 0$ is taken a non-normable space of smooth functions with
correspondingly more strict $F\sp$. Further, we also wish to vary the
"nonlinearity\ssp" $\varphi$ in a suitably wide locally convex space $E\ar 1$
of smooth functions, and still obtain a {\it smooth\ssp} solution map $
\varSigma:E=E\ar 0\sqcap E\ar 1\iinc U\to F$ defined by $x=(\sp y\ar 0\ssp,
\varphi\sp)\mapsto y\ssp$. In the still unfinish- ed paper [\,Hi4\,]\ssp, we
study this kind of {\it generalized\ssp} well-posedness for the previous wave
equation $\,\wave u=\psi\circ[\,\roman{id\,};u\ssp,\partial\ar 0\sp u\ssp,
\roman{grad\,}u\,]\ssp$.

As a tool in this work we use a corollary (= Theorem 4.3 below) of the
implicit function theorem [\,Hi1\ssp; p.\ 235\,] of the present author\ssp.
Since a detailed proof of existence and smoothness of $\varSigma$ requires a
fair amount of space it would have been impractical to include it together
with the required background material in one article. Hence, we have chosen to
develop here the required differential calculus and present the proofs of the
implicit function theorems. Then (in Section 5 below) we give various less
laborious examples which illustrate the technique of application of our
implicit function theorems 4.1\ssp, 4.3\ssp, and 4.5 to (generalized)
well-posedness questions of the type exemplified by the above wave equation.

The physical motivation for this kind of study may be seen as follows.
Physicists tend to propose various particular nonlinear functions as some kind
of "corrections\ssp" to traditional linear partial differential equations used
to model natural phenomena. There may not be any firm ground for the
assumption that a particular proposition would represent the "right\ssp"
natural law. Instead, one can think that it is only an approximation, and that
the "right\ssp" law is in some neighborhood of the proposed one. One then
would wish that a small smooth change in the parameter law causes a
corresponding smooth change in the observable phenomenon represented by the
solution of the equation. Thanks to our implicit function theorems of Section
4 below, we can handle parameter changes in open sets of general locally
convex spaces. A restriction of the Nash\,--\,Moser theorems is that there the
parameter must vary in an open set of a graded Fr\'echet space.

However, our theorems of Section 4 have a flaw associated with the so-called
"loss of derivatives\ssp". A partial remedy (at the cost of increase of
complexity of the required proofs in applications) for this flaw is given by
our newest implicit function theorem in [\,Hi3\,]\ssp. Besides proving the new
implicit function theorem in that paper, we develop a general scheme, the
concept of {\it Banach space situation\ssp}, in order to be able to formulate
various implicit and inverse function theorems within a unified framework so
that their differences can be seen more clearly.
                                        
To illustrate the motivation for the approach used in [\,Hi3\,]\ssp, we take
the previous initial value problem $\partial\ar 0\sp y=P\ssp y+\varphi
\circ[\,\roman{id\,},R\ssp y\,]$ with $y\ssp(\sp 0\,,\cdot)=y\ar 0
$ which also can\Newline be written as an integral equation \hfil $
y\ssp(\sp t\ssp,\cdot)=y\ar 0+\int_{\,0\,}^{\,t}(P\ssp y+\varphi\circ
[\,\roman{id\,},R\ssp y\,]\sp)\ssp(\sp t\ssp,\cdot)\,d\sp t\ssp$.\linebreak
This reformulation, however, generally leaves us a "loss of derivatives\ssp"
--- a phrase without precise meaning unless we specify a sufficiently
structured situation. Alternatively, in favourable cases we can invert the
linear map $y\mapsto(\sp y\ar 0\ssp,\partial\ar 0\ssp y-P\ssp y\sp)$ by some $
\varLambda\ssp,$ in which case we can write the problem as $y=\varLambda\ssp(\sp
y\ar 0\ssp,\varphi\circ[\,\roman{id\,},R\ssp y\,]\sp)\ssp$.

For our theorems of Section 4 to be applicable here, it is required that the
partial differential operator $P$ suitably "dominates\ssp" $R$ in the sense
that the space $F$ can be expressed as a projective limit of a system $\cal F=
\seq{\,F_i\ssp}$ of Banach spaces in such a way that we then obtain $\CPi k$
maps $f_i:E\sqcap F_i\to F_i$ defined by the prescription $(\sp x\ssp,y\sp)=
(\sp y\ar 0\ssp,\varphi\ssp,y\sp)\mapsto y-\varLambda\ssp(\sp y\ar 0\ssp,
\varphi\circ[\,\roman{id\,},R\ssp y\,]\sp)\ssp$.

Nash\,--\,Moser theorems as well as our newest implicit function theorem
in [\,Hi3\,] are applicable in situations with "loss of derivatives\ssp" ---
this loosely referring to the fact that we only have maps $f_{ij}:E\sqcap F_j
\to F_i$ with $i$ sufficiently preceding $j$ in the partial order associated
with the projective system $\cal F$. This kind of theo- rems are sometimes
called "hard\ssp", obviously not because they would be particu- larly hard to
prove but because their application requires lots of work.

Nash\,--\,Moser implicit function theorems for a map $f:E\sqcap F\to F$
require both spaces $E\ssp,F$ to be (graded) Fr\'echet, whereas for our
theorems of Section 4 below and in [\,Hi3\,] neither space needs to be
Fr\'echet. This is utilized in [\,Hi2\,] where we construct the Lie group Diff$\ssp
(M\sp)$ of diffeomorphisms of a smooth finite dimensional paracompact (but not
necessarily second countable) manifold $M$ using Theorem 4.3 as a tool. There
we have $E=F=\cal D\ssp(\sp\tanB\sp M\sp)\ssp,$ the inductive limit space of
com- pactly supported smooth sections of the tangent bundle $\tanB\sp M$ of $M
$ which is ana- logous to the test function spaces $\cal D\ssp(\Omega)\ssp$.

We remark that in [\,Hi4\,] we also consider solution maps for the more
general quasilinear \hfill ("perturbations of\," the) \hfill wave \hfill
equations \hfill $L\ssp(\sp y)\ssp y=\varphi\circ(\sp\barmj y\sp)\ssp,$ \hfill
where\linebreak $\barmj y=[\,\sp\roman{id}\,;y\ssp,\partial\ar 0\sp y\ssp,
\partial\ar 1\sp y\ssp,\ldots\,\partial\ai N\sp y\,]$ and $L\ssp(\sp y)\ssp v=
\sum_{\iota_1=0}^{\,\ssmb N}\sum_{\iota_2=0}^{\,\ssmb N}
\varphi^{\ssp\iota_{\sn1}\iota_2}\circ(\sp\barmj y\sp)\,.\,(\sp
\partial_{\iota_1}\partial_{\iota_2}v\sp)\ssp$, and we also allow $y$ to have
values in a possibly infinite dimensional Hilbert space. For these, however,
we have to apply the methods of [\,Hi3\,]\ssp, and the differentiability
properties of the resulting solution maps are generally worse than in the case
of the "pure\ssp" wave equations $\ssp\wave y=\varphi\circ(\sp\barmj y\sp)\ssp$.

The organization of the present article is as follows.

In Section 1\sp, we introduce our formalism for treating locally convex spaces
which we always require to be Hausdorff. It is a common habit to denote the
algebraic\,--\,topological structure and the underlying set of a topological
vector space by the same symbol, say '\sn$E\ssp$'. This convention has the
defect that in contexts where an explicit distiction has to be made, one has
to insert additional verbal explanations and use various ad hoc notations. To
avoid the need to resort to such unaesthetic means, and to have a unified
formalism, we have abandoned the usual ambiguous convention. Instead, we
understand that a {\it locally convex space is\ssp} the algebraic\,--\,
topological structure $E=(X\sp,\cal T\sp)\ssp,$ where further the pair $X=(\sp
a\ssp,c\sp)$ is the underlying {\it vector structure\ssp} with $a$ the vector
addition and $c$ the scalar multiplication. The underlying set is denoted by $
\vecs E\ssp,$ and is determined by any of the sets $a\ssp,c\,,\cal T,$ since
for example it equals the range of both functions $a$ and $c\,,$ and is the
largest member of the topology $\cal T$.

In Section 2\ssp, we introduce the locally convex spaces $F$ of $C^{\ssp i\,}
$--\,functions $y:\varPi\iinc Q\to\varPi\ar{\sn 1}$ needed to obtain the
(generalized) well-posedness results in Section 5 be- low. Here $\varPi\sp,
\varPi\ar{\sn 1}$ are locally convex with $\varPi$ always finite dimensional. We
mainly need the case where also $\varPi\ar{\sn 1}$ is such, and then $F$ is
Fr\'echet or Banach. Of central importance in obtaining the results of Section
5 is the knowledge of smoothness of \hfil maps \hfil like \hfil $f:(\sp x\ssp,
y\sp)\mapsto x\circ[\,\roman{id\,},y\,]=z\ssp,$ \hfil where \hfil $z\ssp(\eta)
=x\ssp(\sp\eta\ssp,y\ssp(\eta))$ \hfil for \hfil $\eta\in Q\ssp$.\linebreak
Lemma 2.7 prepares us for Theorem 3.6 which gives smoothness of $f\sp$.

In [\,Hi2\,] for the construction of the differentiable structure for the
group of diffeomorphisms, we need some spaces of compactly supported smooth
functions. Section 2 here being the natural place for the introduction of
these spaces, we have included their definitions there although we do not need
them in the present article. An analogous remark applies to the spaces of
functions with bounded partial derivatives which are needed in the treatment
of the wave equation.

Section 3 contains the background about differentiability concepts and results
for maps between general locally convex spaces needed to state and prove the
implicit function theorems of Section 4 and their applications later.

In Section 5\ssp, we give several examples to illustrate the technique of
application of the theorems of Section 4\ssp. In these examples, we prove
smooth dependence of the solution on initial\sp/\sp boundary data and the
nonlinearity for nonlinear (\sp partial\sp) differential equations. We have
chosen the examples so that fairly complete proofs could be presented in
reasonable space, and hence the results here are not supposed to be an
addition of utmost importance to the mainstream of research.

The purpose of Section 6 is to present a formalization of some well-posedness
concepts which in the literature are generally used only in a loose manner. We
just propose some definitions, and then give several examples where previous
known results are formulated with the aid of our precise concepts. The purpose
of the introduction of these concepts is to make possible the exact formulation
of different "ituitive\ssp" well-posedness results so that one could more
clearly see the differences between the "nature\ssp" of these results despite
the fact that they concern (\sp partial\sp) differential equations of various kind.

To be able to express our assertions concisely and precisely at the same time,
and to avoid the use of various ad hoc notations, we follow a rather strict
set theoretic formalism in the spirit of [\,Ky\ssp; Appendix, pp.\
250\,--\,281\,]\ssp, with only few exceptions. To spare the reader some of the
annoying trouble of searching the definition of each particular notation from
the body of the text, we have here collected some of our general set theoretic
notations under the title

\NS{\fssit
Some conventions.} The symbol $\bold K$ will generally denote either the field
$\bold R$ of {\it real numbers\ssp} or $\bold C$ of {\it complex numbers\ssp}.
The underlying sets of these fields we write $\Ke\ssp,\Re\,,\Ce\ssp,$
respectively. We put $\Rep=\{\,s\in\Re:s>0\,\}$ and $\Repp=\Rep\cup
\{\sp 0\sp\}\ssp$. The set of {\it integers\ssp} is $\Ze=\bigcap\ssp\{\,Z:0\in
Z\inc\Re$ and $\all{n\in Z}\,n+1\ssp,n-1\in Z\,\}\ssp,$ and we put $\Zepp=\Ze
\cap\Rep$ and $\Zep=\Zepp\sn\cup\{\sp 0\sp\}\ssp$. The set of {\it natural
numbers\ssp} is $\No=\bigcap\ssp\{\,N:\emptyset\in N$ and $\all{k\in N}\,k^+
\in N\ssp\}\ssp,$ where $k^+=k\cup\{\sp k\sp\}\ssp$. Also let $k^{++}=(\sp
k^+)^+$. We put $\N=\No\sn\setminus\sn\{\sp\emptyset\sp\}\ssp$. There is a
natural bijection $\beta:\No\to\Zep\ssp,$ and we shall generally "identify\ssp"
$k\in\No$ and $\beta\ssp(k\sp)\in\Ze$ by using the same symbol for them. For $
k\in\No$ we then may write $k^+=k+1=\{\,0\,,1\ssp,\ldots\,k\,\}\ssp$. If we
had to make an explicit distinction, we would write ${\iota_{}}_{\roman n}\sp
k=\beta\ssp(k\sp)$ and $\nu_{_{\sp\roman I}}\sp n=n.=\beta\sp\inve(n)\ssp$.
For a finite function $\sigma:D\to\Zep\ssp,$ we write\par\centerline{$
\sum\sigma=\beta\ssp(\Card(\ssp\bigcup\ssp\{\sp\{\sp i\sp\}\sn\times\sn\beta
\sp\inve(n):(\sp i\ssp,n\sp)\in\sigma\,\}))\ssp,$}
\noin which more descriptively could also be written $\sum_{\,i\ssp\in\sp D}
\sigma\ssp(i\sp)\ssp$.

We use the word "family\ssp" as a synonym for "function\ssp". We put $\seq{\,
{t_{}}_z:z\in S\,}=\mathbreak\{\,w:\exi{z\in S}\,w=(\sp z\ssp,{t_{}}_z)\ssp\}
\ssp$. A {\it finite sequence\ssp} is any function $x$ with $k=\dom x\in\mathbreak
\No\ssp$. Then we may write $x=\seq{\,x_1,\ldots\,x_k\ssp}=\seq{\,
{x_{}}_{i\ssp+\sp1}:i\in k\,}\ssp$. If also $y=\seq{\,y_1,\ldots\,y_l\ssp}\ssp
$,\Newline then $x\sp\conc y=\seq{\,x_1\sp,\ldots\,x_k\sp,y_1\sp,\ldots\,
y\ssp\ai l\sp}\ssp$. Finite sequences are things different from "tuples\ssp",
which are defined by the recursion schema $(x)=x$ and $(x_1,\ldots\,x_l,
x_{l+1})=$ $((x_1,\ldots$ $x_l),x_{l+1})\ssp$. We shall apply the definition
schema\par$\mhyppy{3}
(\sp x^1_1\sp,\ldots\,x^1_{i_1}\ssp;x^2_1\sp,\ldots\,x^2_{i_2}\ssp;\ldots\,
x^k_1\sp,\ldots\,x^k_{i_k}\ssp)$\newline$\mhyppy{46}
=((\sp x^1_1\sp,\ldots\,x^1_{i_1})\ssp,(\sp x^2_1\sp,\ldots\,x^2_{i_2})\ssp,
\ldots\,(\sp x^k_1\sp,\ldots\,x^k_{i_k}))\ssp,$

\noin from which we in particular get $(\sp x_1,\ldots\,x_k\ssp;y_1,\ldots\,
y_l)=(\sp x_1,\ldots\,x_k\sp,(\sp y_1,\ldots\,y_l))\ssp$.

The {\it value\,}, cf.\ [\,Ky\ssp; p.\ 261, Def.\ 68\,]\ssp, of a function $f$
at $x$ is $f\sp\value x=f\sp\value(x)\ssp,$ which we generally (following the
usual customs) write simply $f\sp(x)\ssp,$ and sometimes\,\footnote{Note that
  $x\ssp(i)=x_{i\sp+1}$ when $x=\seq{\,x_1\sp,\ldots\,x_k\sp}\ssp$.}
even\Newline $f\sp x$ or $f{_{\sn}}_x\,$. Let $f\ssp[\sp A\ssp]=
f\ssp\image\sn\!A=\{\,y:\exi{x\in A}\,(\sp x\ssp,y\sp)\in f\,\}\ssp$. Adapting
[\,Ky\ssp; p.\ 259, Defs.\ 51, 52\,]\ssp, for any ordered pair $P=(A\,,B\sp)
\ssp,$ we write $A=\sigrd P$ and $B=\taurd P\sp$. Recursively, we define $
\sigrd\KN2\yi i\KP1Z$ by $\sigrd\KN{2.3}\yr 0\KP{.7}Z=Z$ and $
\sigrd\KN{2}\yi i\yplus\yyr 1\sp Z=\sigrd(\sp\sigrd\KN2\yi i\KP1\sp Z\sp)\ssp.
$ The {\it inverse\ssp} of\Newline a class $f$ is $f\sp\inve=\{\ssp(\sp y\ssp,x\sp):
(\sp x\ssp,y\sp)\in f\,\}\ssp$. We write $\cal T\lei A=\{\,U\cap A:U\in\cal T
\,\}\ssp,$ for example trace of a topology or filter.

For functions $f\sp,\sp g$ the symbol $[\,f\sp,\sp g\,]$ denotes the function
$(\dom f\sp)\cap(\dom g\sp)\owns x\mapsto(\sp f\sp(x)\ssp,\sp g\ssp(x))\ssp,$
and $f\risti 2 g$ is defined by $(\sp x\ssp,y\sp)\mapsto(\sp f\sp(x)\ssp,\sp
g\ssp(\sp y))\ssp$. We write $\roman{id}=\mathbreak\{\ssp(\sp x\ssp,x\sp):x=x
\,\}\ssp,$ the identity function on the universe $\Univ=\{\,x:x=x\,\}\ssp$.

Recalling (cf.\ [\,Ky\ssp; Def.\ 1, p.\ 252\,]\ssp, [\,Du\ssp; Notes A2, p.\
408\,]\sp) that a class $x$ is called a {\it set\ssp} if{}f $x\in y$ for some
$y\ssp,$ and otherwise a {\it proper\ssp} class, a {\it large family\ssp} is a
function which is not a set, otherwise a {\it small\ssp} family.

In general, we follow the principle that if we have a term '\sp t\ssp$x\ssp$'
with one free variable '$x\sp$', we apply the scheme for recursive definition
that t$^0x=x$ and t$^{k+1}x=\roman t\ssp(\sp\roman t^kx\sp)$ for $k\in\No\ssp.
$ This gives meaning for example to $\roman{dom}\yr 2\sn f$ once we have
defined $\dom f=\{\,x:\exi y\,(\sp x\ssp,y\sp)\in f\,\}\ssp$. Also note $\rng
f=\dom f\sp\inve$ for any class $f\sp$.

With only some conventional exceptions, for example $x+y\,z=x+(\sp y\,z\sp)
\ssp,$ if we have "binary symbols\ssp" $\roman b\sp\ai k\sp,\ldots\,
\roman b\ar 1\ssp,$ we follow the principle of {\it reduction of parentheses\ssp} 
given by the scheme $x\ai k\ssp\roman b\ai k\ssp\ldots\,x\ar 1\ssp
\roman b\ar 1\ssp x\ar 0=(\sp x\ai k\ssp\roman b\ai k\ssp\ldots\,x\ar 1)\,
\roman b\ar 1\ssp x\ar 0\ssp$. In this manner, for example \hbox{$A\cup B\cup
C$} and $[\,f\sp,g\ssp,h\,]$ are defined. Usually $x\value y\value z\not=
x\value(\sp y\value z\sp)\ssp$.

We further list some definitions:

$\cal T${\it--\,closure\ssp} of $A$ = $\ClT A=(\sp\bigcup\sp\cal T\sp)\sn
\setminus\sn\bigcup\ssp\{\,U\sn\in\cal T:A\cap U=\emptyset\,\}\ssp,$

$\cal T${\it--\,interior\ssp} of $A$ = $\IntT A
=\bigcup\ssp\{\,U\sn\in\cal T:U\inc A\,\}\ssp,$

$\cal T\rist\cal U=\{\,\bigcup\cal A:\cal A\inc\{\,U\sn\times\sn V:U\sn\in
\cal T$ and $V\ssn\in\cal U\,\}\ssp\}$\newline\hyppy{14.2mm}
= the Tihonov/\ssp Tychonoff product topology of $\sp\cal T\sn$ and $\sp\cal U\sp{},$

$f\sp.\,g=\{\ssp(\sp x\ssp,y\sp):\exi{v\ssp,z}\,(\sp x\ssp,v\sp)\in f$ and $
 (\sp x\ssp,z\sp)\in g$ and $(\sp z\ssp,y\sp)\in v\,\}\ssp,$\par\hyppy{10.5mm}
if $f\sp,\sp g$ are functions with $f\sp\value x$ a function having\newline\hyppy{14mm}
$g\value x\in\dom(\sp f\sp\value x\sp)$ for $x\in D=(\dom\sn f\sp)\cap(\dom g
 \sp)\ssp,$ then\newline\hyppy{14mm}
$f\sp.\,g$ is the function $D\owns x\mapsto f\sp\value x\value(\sp g\value x
\sp)=(\sp f\sp(x))\ssp(\sp g\ssp(x))\ssp$.

Hoping that these examples have sufficiently illustrated our principles of
defining notations, sometimes (hopefully) self-explanatory notations and
abbreviations will be used, as for example $\seq{\,{t_{}}_i\ssp}=\seq{\,
{t_{}}_i\sn:i\in I\,}\ssp,$ when $I$ is clear from the context. The internal
references we make so that for example Theorem 3 means that theorem in the
same section, Theorem 4.3 refers to theorem 3 of section 4.


\subhead 1

                  Conventions for locally convex spaces

We consider real or complex Hausdorff locally convex topological vector
spaces. Precisely, this means the following. Letting $\bold K$ be either $
\bold R$ or $\bold C\ssp$, we understand that \hbox{$\bold K=(\sp
a_{_{I\!\!K}}\sp,c_{_{I\!\!K}})\ssp,$} where \hbox{$a_{_{I\!\!K}}\sp,
c_{_{I\!\!K\!}}:\Ke^{\times2}=\Ke\times\Ke\to\Ke$} are the field addition and
multiplication, respectively. Then a vector space over $\bold K\,$, or a $
\bold K\,$--\,vector space is a pair \hbox{$X=(\sp a\ssp,c\sp)$} where for
some nonempty set $S$ we have the functions \hbox{$a:S^{\times2}\to S$} and \hbox{$
c:\Ke\times S\to S\ssp$,} the vector addition and scalar multi\-plication,
respectively. Writing \hbox{$\biit K=(\sp\bold K\ssp,\tau_{_{I\!\!K}})\ssp,$}
where $\tau_{_{I\!\!K}}$ is the natural\sp/\sp standard topology of $\Ke\sp$,
we have $\biit K$ a topological field. A topogical $\biit K\,$--\,vector space
now is any pair\sp/\sp triplet $E=(X\sp,\cal T\sp)=(\sp a\ssp,c\,,\cal T\sp)$
where $X=\sigrd E$ is a $\bold K\,$--\,vector space and $\cal T$ is a topology
(= set of open sets) on the {\it underlying set\ssp} $S=\vecs E=\vecss\sp X=
\rng a$ $=\bigcup\cal T$ such that $a$ and $c$ are continuous $\cal T\rist
\cal T\to\cal T$ and $\tau_{_{I\!\!K}}\sn\rist\cal T\to\cal T_{\sp}$, resp.

We write \hbox{TVS = TVS$\ssp(\biit K\ssp)$} for the class of all topological
$\biit K\,$--\,vector spaces $E$ having \hbox{$\cal T=\taurd E$} a Hausdorff
topology, the notation indicating that '$(\biit K\ssp)\sp$' may be dropped if
the scalar field is understood from the context, or unessential. We denote by
LCS and BaS the subclasses of TVS having as members the locally convex ones or
Banach\sp((iz)\sp able) spaces, respectively.

The zero vector of $E\in{}$TVS is $\bnull E=\bigcup\ssp\{\,x:(\sp x\ssp,x\ssp,
x\sp)\in\sigrd\KN{2.3}\yr 2\KP{.7}E\,\}\ssp$. The filter of zero neighborhoods
is written $\ymp E\ssp,$ and the von Neumann bornology (= set of bounded sets)
is $\rajou E\ssp$. In a connection where we are dealing with a Banachable
space $E$ with a specified norm $x\mapsto\|\ssp x\ssp\|\ssp,$ we may use the
notations $B\ai E\sp(\smb R)=\{\,x:\|\ssp x\ssp\|<\smb R\ssp\}$ and $
\bar B\ai E\sp(\smb R)=\{\,x:\|\ssp x\ssp\|\le\smb R\ssp\}\ssp,$ although this
is not quite logical since the topology does not determine the norm, cf.\
[\,Jr\ssp; p.\ 115\,]\ssp.

For \hbox{$X=(\sp a\ssp,c\sp)\in\roman{VS}\ssp(\bold K)$} and any $S$ writing
\hbox{$X_{\sp|\ssp S}=(\sp a\,|\,S^{\times2\!},c\,|\,(\sp\Univ\sn\times\sn
S\sp))\ssp$}, we further put \hbox{$E_{\sp/\sp S}=(X_{\sp|\ssp S}\ssp,\cal T
\lei S\sp)$} when \hbox{$E=(X\sp,\cal T\sp)$} is any topological vector space.
If we have a topological vector space $F$ of functions $Q\to\Ke$ with $\sigrd
F$ a vector substructure of $\bold K^{\,Q}$ for some set $Q\ssp,$ then we
write $0\sp\ai Q=Q\sn\times\sn\{\sp 0\sp\}=\bnull F\ssp$.

By a {\it vector map\ssp} (of l.c.s.\ = locally convex spaces) we here mean
any triple $\tilde f=(E\ssp,F\sp,f\sp)\ssp,$ where $E\ssp,F\in$ LCS and $f$ is
a function with $f\inc(\sp\vecs E\sp)\times(\sp\vecs F\sp)\ssp,$ cf.\
[\,Hi1\ssp; p.\ 238\,]\ssp. The set of continuous linear maps $E\to F$ we
denote by $\cal L\ssp(E\ssp,F\sp)\ssp,$ and then the {\it class\ssp} \hbox{of
continuous linear maps (of l.c.s.) is $\cal L=\{\ssp(E\ssp,F\sp,\sp\ell\,):
E\ssp,F\in\roman{LCS}$ and} $\ssp\ell\in\cal L\ssp(E\ssp,F\sp)\ssp\}\ssp$. The
product of $E\ssp,F\in$ LCS we write $E\sqcap F\sp,$ and of a nonempty small
family $\biit E\in\roman{LCS}^{\,I}$ correspondingly $\prod\biit E\ssp,$ \hfill
and further $E^{\,I}=\prod\,(\sp I\times\{\sp E\sp\})\ssp$. Then $\taurd
\biit K^{\,k}=\tau_{_{I\!\!K}}\KN{1.8}^k\,$. We put $|\ssp x\ssp|=\big(
\sum_{i=1}^{\,k}|\ssp x_i\sp|^{\,2}\big){}^{\sp\frac12}$ for $x=\seq{\,
x\ar 1\sp,\ldots\,x\sp\ai k\sp}\in\Ke^{\,k}$.

A {\it direction\ssp} on a set $I\ar 1$ is any reflexive and transitive relation
$\Delta$ with $\dom\Delta=I\ar 1\ssp,$ and such that any $i\ssp,\sp j\in
I\ar 1$ have some $k$ with $(\sp i\ssp,k\sp)\ssp,(\sp j\sp,k\sp)\in\Delta\,$.
A family $\cal P$ $=\{\,(i,j;F_j,F_i,\rho_{ij}):(i,j)\in\Delta\,\}\in\cal L^{\,
\Delta}$ we call a {\it projective system\ssp} in LCS if{}f $\Delta$ is a
direction on some $I\ar 1\sp,$ and $\rho_{ii}\inc\roman{id}$ and $\rho_{ij}
\circ\rho_{jk}=\rho_{ik}$ hold for $(i,j),(j,k)\in\Delta\,$. Then a {\it cone\ssp}
over $\cal P$ is any pair $(E\ssp,\bit{032}\sp)\ssp,$ where $\bit{032}=\seq{\,
(E\ssp,{F_{\sn}}_i\ssp,{\rho_{}}_i\sp):i\in I\ar 1\ssp}\in\cal L^{\,I_1}$ with
$\rho_i=\rho_{ij}\circ\rho_j$ for all $(i,j)\in\Delta\ssp$. A {\it projective
limit\ssp} of $\cal P$ is any cone $(F\sp,\bit{031}\sp)$ with $\bit{031}=\seq{\,
(F\sp,{F_{\sn}}_i\ssp,{p_{}}_i\sp):i\in I\ar 1\ssp}$ such that for every cone
$(E\ssp,\bit{032}\sp)\ssp,$ there is a unique $\ell\in\cal L\ssp(E\ssp,F\sp)$
with ${\rho_{}}_i={p_{}}_i\circ\ell$ for all $i\in I\ar 1\sp$.

One observes that for nonempty $I\ar 1\ssp,$ in a cone $(E\ssp,\bit{032}\sp)$
over some $\cal P$ the space $E$ is determined by $\bit{032}$ since $E=\bigcup
\,(\dom\sn(\dom\sn(\rng\bit{032}\sp)))\ssp$. Thus if one wants to restrict for
nonempty $\bit{032}$ only, a cone could equally be defined as the family $
\bit{032}\,$. A projective limit $(F\sp,\emptyset\ssp)$ necessarily has $F$
trivial, i.e., $\vecs F=\{\ssp\bnull F\}\ssp$. In the sequel, we shall need
only nontrivial projective limits. Thus to avoid the need to write down
unnecessary information, we now adopt the convention that a {\it projective
limit is\ssp} the family $\bit{031}$ above which we always assume to be nonempty.

For further facts about locally convex projective limits the reader may consult 
[\,Ho\ssp; Exe.\ 2.11.3, p.\ 155\,] and [\,Jr\ssp; pp.\ 37--38\,]\ssp. In
particular $\taurd F$ is the smallest (= weakest) locally convex topology for
$\sigrd F$ such that $p_i\in\cal L(F,F_i)$ for all $i\ssp$. We also have \hfil
$\bigcap\ssp\{\,p_i^{-\iota}[\sp\{\ssp\bnull F\LHB{.3}{_{_i}}\}\sp]:i\in I_1
\ssp\}=\{\ssp\bnull F\}\ssp,$ \hfil and \hfil $[\ {p_{}}_i:i\in I\ar 1\ssp]:x
\mapsto\seq{\,{p_{}}_i\sp x:i\in I\ar 1\ssp}$\linebreak defines a topological
linear isomorphism of $F$ onto the topological linear subspace $G$ of $\prod\,
\seq{\,F_i:i\in I\ar 1\ssp}$ with $\vecs G$ consisting of $\seq{\,x_i:i\in
I\ar 1\ssp}$ having $x_i=\rho_{ij}x_j$ for all $(i,j)\in\Delta\,$. The limit
space $F$ is already determined by any "tail\ssp" because of

\newProCla 1 Lemma.

If $\seq{\,(F\sp,{F_{\sn}}_i\sp,{p_{}}_i):i\in I\ar 1\ssp}$ is a projective
limit of $\{\,(i,j;F_j,F_i,\rho_{ij}):(i,j)\in\Delta\,\}$ and $\emptyset\not=J
\inc I_1\inc\Delta\inve\sp[\ssp J\,]\ssp,$ then also $\seq{\,(F,F_i,p_i):i\in
J\,}$ is a pro- jective limit of $\{\,(i,j;F_j,F_i,\rho_{ij}):(i,j)\in\Delta
\cap J^{\times2\,}\}\ssp$.

\Prooff Given a locally convex space $E$ and ${\rho_{}}_i\in\cal L\ssp(E\ssp,
{F_{\sn}}_i)$ with ${\rho_{}}_i={\rho_{}}_{ij}\circ{\rho_{}}_j$ for $(i,j)\in
\Delta\cap J^{\times2},$ it suffices to show that there is a unique $\ell\in
\cal L\ssp(E\ssp,F\sp)$ with ${\rho_{}}_i={p_{}}_i\circ\sp\ell$ for all $i\in
J\sp$. To see existence of $\ell\,,$ we first observe that for every $i\in
I\ar 1\ssn\setminus\sn J\sp,$ there is a unique ${\rho_{}}_i$ with $
{\rho_{}}_i={\rho_{}}_{ij}\circ{\rho_{}}_j$ for all $j\in J$ with $(i,j)\in
\Delta\,$. Then we obtain existence of a unique $\ell$ with ${\rho_{}}_i=
{p_{}}_i\circ\sp\ell$ for $i\in I_1\ssp$. Uniqueness of $\ell$ under the
restriction $i\in J$ now easily follows.                              \newQED

Note that in the above proof we needed the assumption that $\Delta$ is a
direction. In [\,Ho\ssp; Exe.\ 2.11.3, p.\ 155\,]\ssp, this is not (at least
explicitly) assumed, cf.\ [\,Ho\ssp; Exe.\ 2.12.4, p.\ 174\,]\ssp, where it is
explicit. Observe that Horv\'ath's order (loc.\ cit., p.\ 2) includes
antisymmetry, which we do not require for a direction.

\NS For locally convex spaces $E=(X\sp,\cal T\sp)$ and $F=(\sp Y,\sp\cal U\sp)
\ssp,$ we write $E\le F$ iff $\idv F$ $=\roman{id}\,|\,(\sp\vecs F\sp)\in
\cal L\ssp(F,E\sp)\ssp,$ i.e., the identity function is linear $Y\to X\sp,$
which means that $Y$ is a vector substructure of $X\sp,$ and continuous $
\cal U\to\cal T$. The proper class $\leLCSr=\{\,(E\ssp,F\sp):E\le F\,\}$ is a
reflexive, antisymmetric, and transitive relation. If\Newline $\cal F\not=
\emptyset$ is lower $\leLCSr\ssp$--\,bounded and downwards\sp/\ssp left $
\leLCSr\ssp$--\,directed, we then have $\leLCS-\inf\cal F=\bigcap\ssp\{\ssp E:
E$ greatest lower $\leLCSr\ssp$--\,bound of $\cal F\,\}\in{}$LCS\ssp. \hfill
Left directed- ness here means that $\leLCSr\KN4\inve\KP1\cap\sp
\cal F^{\times2}$ is a direction.

For a $\bold K\,$--\,vector space $X\sp,$ we write ${\le}\ai X=\leLCSr\cap
\O\ai X^{\,\times2}\ssp,$ when $\O\ai X$ denotes the set of all $F\in{}$LCS
with $\sigrd F$ a vector substructure of $X\sp$. For $E\ssp,F\in\roman{LCS}
\ssp,$ we call $E$ a {\it topological linear subspace\ssp} of $F$ in case we
have $E=\leLCS-\inf\ssp\{\,G:F\le G$ and $\vecs G\inc\vecs E\,
\}\ssp$. Here observe that $\leLCS-\inf\cal F=\Univ\not\in{}$LCS
in case the greatest lower bound does not exist. We now have

\newProCla 2 Lemma.

For $E\ssp,F\in\roman{LCS}$ holds $E=F_{/\upsilon_s E}$ if and only if $E$ is
a topological lin- ear subspace of $F\sp$.

\Prooff Leaving the simpler "$\Rightarrow$\ssp" case as an exercise to the
reader, we show "$\Leftarrow$\ssp". Thus writing $\cal F=\{\,G:F\le G$ and $
\vecs G\inc\vecs E\,\}\ssp,$ and letting $S$ be the lin-\Newline ear $\sigrd F
\,$--\,span of $(\sp\vecs E\sp)\cap(\sp\vecs F\sp)\ssp,$ we have $F_{/S}$ a
lower $\leLCSr\ssp$--\,bound of $\cal F,$ whence $F_{/S}\le E\ssp$.
Consequently, now $S=\vecs E$ is a vector subspace of $\sigrd F\sp$. Then we
have\Newline $F_{/S}\in\cal F,$ whence $E\le F_{/S}\le E\ssp,$ and hence $E=
F_{/S}=F_{/\upsilon_s E}\ssp$.                                           \QED


\subhead 2

                       Particular function spaces

For our applications of theorems of Section 4\ssp, we shall need some locally
convex spaces of differentiable functions \hbox{$Q\to\vecs F\sp,$} which we
now introduce. We have here \hbox{$Q\inc\vecs\varPi$} with \hbox{$\varPi\ssp,F
\in\roman{LCS}\ssp(\biit K\ssp)$} and $\varPi$ finite dimensional. For the
function spaces we shall generally use a symbol of the form $
S\ssp(\sp Q\ai\varPii\sp,F\sp)\ssp$.

If we have a finite sequence $\seq{\,{F_{\sn}}_0\ssp,\ldots\,{F_{\sn}}_k\ssp}$
with $F_0=\biit K$ and $F_k=\varPi\ssp,$ and each $l\in k^{\ssp+}\sn\setminus
\sn\{\sp 0\sp\}$ having some $i\ssp,\sp j\in l$ and $\smb N\in\No$ with $F_l=
F_i\sqcap F_j$ or $F_l=F_i^{\,N}\ssp,$ then any nonempty $Q\inc\vecs\varPi$
uniquely\footnote{assuming $\Ke$ is constructed so that $1.=\{\emptyset\}
\not\in\Ke\sp$.} determines $\varPi\sp$. Since $S\ssp(\ssp\emptyset_\varPi
\ssp,F\sp)$ does not depend on $\varPi\ssp,$ we then may define $S\ssp(\sp
Q\ssp,F\sp)=S\ssp(\sp Q_\varPi\ssp,F\sp)\ssp$. We put $S\ssp(Q)=S\ssp(\sp
Q\ssp,\biit K\,)\ssp$. \hfill In\!\newline case $\biit K=\biit R\,,$ we put $
S\sp(Q)_{\sp\bold c}=S\ssp(\sp Q\ssp,\biit C_{_{I\!\!R}})\ssp,$ where $
\biit C_{_{I\!\!R}}=(\sp a_{_C}\ssp,c_{_C}\ssp|\,(\Re\sn\times\sn\Ce\sp)\ssp,
\tau_{_C})$ is the topological field $\biit C$ considered as a real
topological vector space.

\NSN{\bf
1} {\fssit Continuous functions.} We put $C\ssp(\sp Q_\varPi\sp,F\sp)=(X\sp,
  \cal T\sp)\ssp,$ where $X$ is the vector sub- structure of $(\sp\sigrd F\sp)
^{\,Q}$ with underlying set formed by the functions $x$ continuous $\taurd
\varPi\to\taurd F\sp,$ and $\cal T$ is the compact-open topology. The space $E
=C\ai{bd}\sp(\sp Q_\varPi\sp,F\sp)$ is defined analogously with the following
deviations: In addition to be continuous, for $x$ to belong to $\vecs E\ssp,$
we require $\rng x\in\rajou F\sp$. As a zero basis for $\taurd E$ we take the
filter basis $\{\,N\sp(\sp V):V\in\ymp F\,\}\ssp,$ where $N\sp(\sp V)=\{\,x\in
\vecs E:\rng x\inc V\,\}\ssp$. We write $C.(\sp Q_\varPi\sp,F\sp)$ for the
closed topological linear subspace of $C\ai{bd}\sp(\sp Q_\varPi\sp,F\sp)$
formed by the functions $x$ which "vanish far away\ssp", this meaning that
every $V\in\ymp F$\newline has some $(\sp\taurd\varPi\sp)\lei Q\,$--\,compact
$K$ with $x\ssp[\,Q\sn\setminus\sn K\,]\inc V$. We have the implication\vskip.5mm\centerline{$
E\in\{\,C\ai{bd}\sp(\sp Q_\varPi\sp,F\sp)\ssp,C.(\sp Q_\varPi\sp,F\sp)\,\}$
and $F\in\roman{BaS}\imply E\in\roman{BaS}\,$.}

\noin In case $Q$ is $\taurd\varPi\,$--\,compact, we trivially have $
C.(\sp Q_\varPi\sp,F\sp)=C\ai{bd}\sp(\sp Q_\varPi\sp,F\sp)\ssp$.

\NSN{\bf
2} {\fssit Variations.} For $k\in\No$ and a vector map $\tilde f=(E\ssp,F\sp,
  f\sp)$ adapting [\,AS\ssp; pp.\ 206, 229\,]\ssp, we have the first and $
k^{\,\roman{th}}$ order {\it variation\ssp} maps $\delta\ssp\tilde f=(E\sqcap
E\ssp,F\sp,\delta\ssp f\sp)$ and $\delta^{\,k\sn}\tilde f=(E^{\,k\ssp+\sp 1\!}
,F\sp,\delta^{\,k\sn}f\sp)\ssp$, where the variation functions $\delta\sp f=
\delta\ai{EF}\sp f$ and $\delta^{\,k\sn}f=\delta^{\,k}\KN{2.1}\LHB{.3}{\ai{EF}}\sp
f$ are defined as follows. For $u\in\vecs E\ssp,$ we define the {\it directional
derivative\ssp} ${d_{}}_u\sp f$ by\vskip.5mm\centerline{$
{d_{}}_u\sp f\sp(x)=\lim_{\,t\,\to\,0}t^{-1}(\sp f\sp(\sp x+t\,u\sp)-f\sp(x))$}

\noin with the understanding that we have $x\in\dom(\sp{d_{}}_u\sp f\sp)$
exactly in case $x\in\dom f$ and 0 is a limit point of $\{\,t\in\Ke\sn
\setminus\sn\{\sp 0\sp\}:x+t\,u\in\dom f\,\}$ and the previous limit exists
w.r.t.\ the topologies $\tau_{_{I\!\!K}}$ and $\taurd F\sp$. For $\biit u\in(\sp\vecs E
\sp)^{\,k}$ with $k\in\No\ssp,$ we then define $d\ssp(\biit u)\ssp f$ by the
recursion $d\ssp(\emptyset)\ssp f=f$ and $d\ssp(\sp\biit u\sp\conc\seq{\sp
u\sp})\ssp f={d_{}}_u\sp (\sp d\ssp(\biit u)\ssp f\sp)\ssp$. Finally we put $
\delta^{\,k\sn}f\sp(\seq{\sp x\sp}\conc\biit u\sp)=(\sp d\ssp(\biit u)\ssp
f\sp)\ssp(x)$ and $\delta f\sp\value(\sp x\ssp,u\sp)=\delta^{\,1\sn}f\sp\value
\seq{\ssp x\ssp,u\ssp}$ with the understanding that we have $\seq{\sp x\sp}
\conc\biit u\in\dom(\sp\delta^{\,k\sn}f\sp)$ if{}f $\biit u\in(\sp\vecs
E\sp)^{\,k}$ and $x\in\dom(\sp d\ssp(\biit u)\ssp f\sp)\ssp$. One sees that $
\seq{\sp x\sp}\conc(\sp k\sn\times\sn\{\ssp\bnull E\})\in\dom(\sp
\delta^{\,k\sn}f\sp)$ for any $x\in\dom f\sp$.

In case $\varPi=\biit K^{\,k}$ with standard basis $\{\ssp e\sp\ar 1\sp,\ldots\,
e\sp\ai k\sp\}\ssp,$ where we hence have $e\sNor{l\ssp+\sp 1}=$ \hbox{$
\{\ssp(\sp l\ssp,1)\ssp\}\cup(
(\sp k\sn\setminus\sn\{\sp l\sp\})\sn\times\sn\{0\})\ssp$}, we define the {\it
partial derivatives\ssp} $\partial_j\sp f$ for $j=1\ssp,\ldots\,k\ssp$,
sometimes taking $j=0\,,\ldots\,k-1\ssp$, equivalently written $j\in k\ssp$,
and $\partial^{\,\alpha}f$ for every $\alpha=\seq{\,\alpha\sp\ar 1\sp,\ldots\,
\alpha\sp\ai k\sp}\in\Nopot k$ by $\partial_j\sp f=
{d_{}}_{\hbox{\font\=cmmi8\e}_j}\sp f$ and $\partial^{\,\alpha}f=
d\ssp(\biit u)\ssp f\sp$, \hfill where we have taken\vskip-1.2mm\noin $\biit u
=(\sp\alpha\sp\ar 1\sn\times\sn\{\ssp e\sp\ar 1\})\concc\sn\ldots\,
(\sp\alpha\sp\ai k\sn\times\sn\{\ssp e\sp\ai k\})\ssp$. We write $|\ssp\alpha\sp|
=\sum_{\nu=1}^{\,k}\alpha_\nu\ssp$, the {\it order\ssp} of $\alpha\ssp$.

We shall need these variations (or directional derivatives) in the following
section for general \hbox{$E\in\roman{LCS}\ssp$}, but here mainly for finite
dimensional ones when constructing the following spaces of

\NSN{\bf
3} {\fssit Differentiable functions.} Putting the restriction $Q\inc
  \roman{Cl_{\ssp}}_{{\tau_{}}_{rd}\sp\varPi}\ssp O\ssp,$ where we write $O=
\roman{Int_{\ssp}}_{{\tau_{}}_{rd}\sp\varPi}\ssp Q\ssp$, \hfil for \hfil a \hfil
function \hfil $f\in(\sp\vecs F\sp)^{\,Q}$ \hfil and \hfil for \hfil any \hfil
$\biit u\in(\sp\vecs\varPi\sp)^{\,k}\ssp$, \hfil we \hfil define\linebreak the \hfil
extended \hfil iterated \hfil directional \hfil derivative \hfil $
\bar d\ssp(\biit u)\ssp f$ \hfil by \hfil $(\sp\bar d\ssp(\biit u)\ssp f\sp)
\ssp(x)=(\sp d\ssp(\biit u)\ssp f\sp)\ssp(x)$\linebreak for \hfil $x\in P=O
\cap(\dom(\sp d\ssp(\biit u)\ssp f\sp))\ssp,$ \hfil and \hfil by \hfil $(\sp
\bar d\ssp(\biit u)\ssp f\sp)\ssp(x)=\lim_{\,v\to x\,}(\sp d\ssp(\biit u)\ssp
f\,|\,O\sp)\ssp(v)$\linebreak for $x\in Q\sn\setminus\sn P$ in case a (unique)
limit exists. \hfil Agreeing that \hfil $\infty+1=\infty\,$, \hfil for any\linebreak
$i\in\No\cup\{\infty\}\ssp,$ we then define

$\mhyppy{2}C^{\,i}(\sp Q_\varPi\sp,F\sp)=\leLCS-\inf\,\{\,G:C\ssp(\sp
Q_\varPi\sp,F\sp)\le G$ and $\all{k\ssp,\biit u}\,k\in\No$ and

$\mhyppy{24}\biit u\in(\sp\vecs\varPi\sp)^{\,k}$ and $k<i+1\imply
\bar d\ssp\ai G(\biit u)\in\cal L\ssp(\sp G\ssp,C\ssp(\sp Q_\varPi\sp,F\sp))\ssp\}\,,$

\noin where $\bar d\ssp\ai G(\biit u)=\{\ssp(\sp f\sp,g\sp):f\in\vecs G$ and $
g=\bar d\ssp(\biit u)\ssp f\,\}\ssp$.

We let $\cal D_K(\sp Q_\varPi\sp,F\sp)$ be the topological linear subspace of $
\Cinfty(\sp Q_\varPi\sp,F\sp)$ formed by\Newline the $x$ having $\supp x=
\roman{Cl_{\ssp}}_{{\tau_{}}_{rd}\sp\varPi}\sp(\sp x\sp\inve\sp[\,(\vecs F\sp)
\sn\setminus\sn\{\ssp\bnull F\}\,]\sp)\inc K\sp$. Then we put $\,
\cal D\ssp(\sp Q_\varPi\sp,F\sp)$ $=\Cinftyzero(\sp Q_\varPi\sp,F\sp)=
\leLCS-\inf\,\{\,\cal D_K(\sp Q_\varPi\sp,F\sp):K$ is $(\sp\taurd\varPi\sp)
\lei Q\,$--\,compact $\}\,$. Still let

$\mhyppy{2}C^{\,i}_{bd}(\sp Q_\varPi\sp,F\sp)=\leLCS-\inf\,\{\,G:C\ai{bd}\sp(\sp
Q_\varPi\sp,F\sp)\le G$ and $\all{k\ssp,\biit u}\,k\in\No$ and

$\mhyppy{24}\biit u\in(\sp\vecs\varPi\sp)^{\,k}$ and $k<i+1\imply
\bar d\ssp\ai G(\biit u)\in\cal L\ssp(\sp G\ssp,C\ai{bd}\sp(\sp Q_\varPi\sp,F\sp))\ssp\}\,$.

In case $\biit K=\biit C$ and $\varOmega\in\taurd\varPi\ssp,$ we write $
C^{\ssp 1\sp}(\sp\varOmega_\varPi\sp,F\sp)=H\sp(\sp\varOmega_\varPi\sp,F\sp)
\ssp,$ and further let $H\ai b\sp(\sp Q_\varPi\sp,F\sp)=G$ be the topological
linear subspace of $C^{\,\sp 0\sp}(\sp Q_\varPi\sp,F\sp)=E$ with $\vecs G=\{\,
x\in\vecs E:x\,|\,\varOmega\in\vecs H\sp(\sp\varOmega_\varPi\sp,F\sp)\ssp\}
\ssp,$ where $\varOmega=\roman{Int_{\ssp}}_{{\tau_{}}_{rd}\sp\varPi}\ssp Q\ssp
$. If $Q$ is $\taurd\varPi$ --\,compact, we have the implication $F\in
\roman{BaS}\ssp(\biit C\ssp)\imply H\ai b\sp(\sp Q_\varPi\sp,F\sp)\in
\roman{BaS}\ssp(\biit C\ssp)\ssp$. The vectors of $H\sp(\sp\varOmega_\varPi\sp,
F\sp)$ are called {\it holomorphic\ssp} functions $\varPi\iinc\varOmega\to F\sp$.

\NSN{\bf
4} \Remarkss Our method above of defining the spaces of differentiable
  functions as certain locally convex greatest lower bounds is only a means of
stating the definitions concisely, still precisely. Apriori, it is not at all
clear that for example \hbox{$E=C^{\ssp i\sp}(\sp Q\ai\varPii\sp,F\sp)$} even
is a topological vector space, or that $E\in\roman{LCS}\ssp(\biit K\ssp)\ssp$.

According to our definitions, it is perfectly possible for \hbox{$E=\Univ$} to
hold, and to prevent this, a proof is required. There one proceeds in the
manner which one usually puts as part of the definition (thus mixing
definitions, propositions, and their proofs)\ssp. Namely, one first specifies
a suitable set \hbox{$S\inc(\sp\vecs F\sp)^{\,Q}$} and a prospective filter
basis $\cal V$ for the neighborhoods at $Q\sn\times\sn\{\sp\bnull F\}\ssp$,
then constructs the topology $\cal T$ for $S$ by translations of $\cal V_{\sp}
$, and finally puts $E=(((\sp\sigrd F\sp)^{\,Q\sp})_{\ssp|\ssp S}\ssp,
\cal T\sp)\ssp$. It then remains to prove that indeed $E$ is the required $
\leLCSr\,$--\,infimum. The standard details we leave as exercises for the
mathematically matured reader.

Under our basic assumption \hbox{$Q\inc
\roman{Cl_{\ssp}}_{{\tau_{}}_{rd}\sp\varPi}\sp(\sp
\roman{Int_{\ssp}}_{{\tau_{}}_{rd}\sp\varPi}\ssp Q\sp)$} which suffices for
the unique definition of the extended directional derivatives at boundary
points, not even the space \hbox{$E=C\ssp(Q)$} generally needs to be Fr\'echet,
i.e., the neighborhood filter $\ymp E$ does not need to have a countable base.
To give an example, with $\xi=(\sp 0\,,0\sp)$ we take $Q=\{\sp\xi\sp\}\cup
\ssp{]}\,0\,,1\ssp{[}^{\,\times2}$, hence $\varPi=\biit R\sp\sqcap\biit R\,$
and $F=\biit R\,$. Supposing $\ymp E$ has a countable base, there exists a
sequence $\seq{\,{K_{}}_n:n\in\No\ssp}$ of compact subsets of $Q$ such that
every compact $K\inc Q$ has some $n$ with $K\inc{K_{}}_n\ssp$. Then we
construct a sequence $\biit s=\seq{\,(\sp{s_{}}_n\sp,(\sp n+1)^{\sp-1\sp}):n
\in\No\ssp}\to\xi$ in top $\taurd\sp\varPi$ with $\biit s\value n\in
Q\sn\setminus\sn{K_{}}_n$ for all $n\in\No\ssp$. Taking $K=\{\sp\xi\sp\}\cup
(\rng\biit s\sp)\ssp$, a contradiction follows.

A more specific result is the following: Let the real Banach space $F$ be at
least one dimensional, \hfil and let \hfil $Q=O\cup B$ \hfil where \hfil
$O\in\taurd\varPi$ and $B\inc C={\partial_{\ssp}}_{{\tau_{}}_{rd}\sp\varPi}\ssp
O\ssp$.\linebreak Then \hfil $C^{\,i}(\sp Q_\varPi\sp,F\sp)$ is Fr\'echet $
\equivv B\in(\sp\taurd\varPi\sp)\lei C\ssp,$ \hfil i.e., \hfil the \hfil set \hfil
$B$ \hfil is \hfil open \hfil in \hfil the\linebreak relative topology of the
boundary $C$ of the open set $O\ssp$.

We also mention the following subtlety in the previous definitions. If we had
chosen \hfil $P=\dom(\sp d\ssp(\biit u)\ssp f\sp)$ \hfil when defining \hfil $
\bar d\ssp(\biit u)\ssp f\sp$, \hfil the resulting spaces \hfil $C^{\ssp i}(\sp
Q_\varPi\sp,F\sp)$\linebreak would have had smaller underlying set, and
generally not even the normable spaces $C^{\sp 1}(Q)$ for compact $Q$ would
have been complete, hence Banach. A counterexample can be constructed for $Q=
\{0\}\sn\cup\bigcup\ssp\{\ssp{I_{}}_n:n\in\N\ssp\}\inc\Re\,$, where with $
{x_{}}_n=2^{-n}$ we\Newline take ${I_{}}_n=[\,{x_{\sn}}_n\sp,{x_{\sn}}_n+x_n^{\,2}\ssp]\ssp$.
\hfil One \hfil considers \hfil the \hfil Cauchy \hfil sequence $\biit f=\seq{\,
{f_{\sn}}_i:i\in\No\ssp}$\linebreak with \hfil ${f_{\sn}}_i=\{\ssp(\sp 0\ssp,
0\sp)\ssp\}\sn\cup\bigcup\ssp\{\ssp\seq{\,x-{x_{}}_n+{\eps_{}}_{in}:x\in
{I_{}}_n\ssp}:n\in\N\ssp\}\,,$ \hfil where\linebreak ${\eps_{}}_{in}=
({-}1)^{\ssp n}\sp(\sp i+1\sp)^{-1}{x_{}}_n\ssp$.

For \hfil $f\in\vecs C^{\ssp i}(\sp Q_\varPi\sp,F\sp)$ \hfil and \hfil a \hfil
boundary \hfil point \hfil $x$ \hfil of \hfil $Q\ssp,$ \hfil we \hfil need \hfil
not \hfil have \hfil the\linebreak equality $d\ssp(\biit u)\ssp f\sp(x)=
\bar d\ssp(\biit u)\ssp f\sp(x)$ even if both sides are defined with $\dom
\biit u\inc i\ssp$. An example for the previous $Q$ is $f=\taurd C^{\sp 1}(Q)
\text{\,-\sp}\lim\sp\biit f\,$ having\vskip.5mm\centerline{$
d\ssp(\sn\seq{\sp 1\sp}\sn)\ssp f\sp(0)={d_{}}_1\sp f\sp(0)=f\sp'(0)=0\not=1=
\underset{x\to 0}\to\lim\,f\sp'(x)=\bar d\ssp(\sn\seq{\sp 1\sp}\sn)\ssp f\sp(0)\ssp$.}

Although we have above assumed $\roman{dim_{_H}}\sp\varPi\in\No\ssp$, i.e.,
that $\varPi$ has finite Hamel dimension, our definitions for the spaces of
differentiable functions are meaningfull also for $\roman{dim_{_H}}\sp\varPi
\not\in\No\ssp$. However, if we have $Q\inc
\roman{Cl_{\ssp}}_{{\tau_{}}_{rd}\sp\varPi}\sp(\sp
\roman{Int_{\ssp}}_{{\tau_{}}_{rd}\sp\varPi}\ssp Q\sp)$ and $\roman{dim_{_H}}\sp
\varPi\not\in\No\ssp$, then $f=Q\sn\times\sn\{\sp\bnull F\}$ is the only
continuous $f:Q\to\vecs F$ with compact support,\Newline whence the spaces $
\cal D\ai K\sp(\sp Q\ai\varPii\sp,F\sp)$ and $
\cal D\ssp(\sp Q\ai\varPii\sp,F\sp)$ are trivial.

For \hbox{$E=C^{\ssp i\sp}(\varOmega\ai{\varPii}\sp,F\sp)$} in the case where
we have \hbox{$\roman{dim_{_H}}\sp\varPi\not\in\No$} and \hbox{$\varOmega\in
\taurd\sp\varPi\sp$}, that is \hbox{$Q=\varOmega$} open in $\varPi\sp$, there
arises the question about the relation of the set $\vecs E$ to e.g.\ the
differentiability classes occurring in [\,Ke\,]\ssp. In the second simplest
case $i=1.=\{\emptyset\}\ssp$, one might first ask the following

\noin{\bf
?} {\fssit Question\ssp}. Suppose we only know that $f$ and $d{_{}}_u\sp f$
  are defined on $\varOmega$ and continuous $\taurd\sp\varPi\to\taurd F$ for
every fixed $u\in\vecs\varPi\sp$. \hfill {\it Does it then follow\ssp} that
the function $f\sp'(x)=\seq{\,d{_{}}_u\sp f\sp(x):u\in\vecs\varPi\,}$ is
continuous $\ssp\taurd\sp\varPi\to\taurd F\,$ for each fixed $x\in\varOmega$\ ?\KN{10}

A proof similar to the one for $w=\bnull F$ in the proof of Theorem 3.5 below
would show that $f\sp'(x)$ is linear $\sigrd\sp\varPi\to\sigrd F\sp$, whence
in the case of a positive answer to Question ?\ssp, we would have $\rng f\sp'
\inc\cal L\ssp(\varPi_{\sp},F\sp)\ssp$. Allowing general $i$ and assuming $f
\in\vecs E\ssp$, the proof of [\,Ke\ssp; Thm.\ 2.4.0, p.\ 90\,] shows that $
d^{\,l\sn}f\sp(x)=\seq{\,(\sp d\ssp(\biit u)\ssp f\sp)\value x:\biit u\in(\sp
\vecs\varPi\sp)^{\,l}\ssp}$ is symmetric $(\sp\vecs\varPi\sp)^{\,l}\to\vecs F$
for all fixed $l\in i+1$ and $x\in\varOmega\ssp$.

If we put the assumption that $\varPi$ is {\it metrizable and barrelled\,},
e.g., metrizable (locally convex) Baire, Fr\'echet, or Banach, the Banach\,--\,
Steinhaus theorem [\,Jr\ssp; Thm.\ 11.1.3, p.\ 220\,] then inductively would
show (\sp cf.\ [\,Jr\ssp; Thm.\ 5.1.4, p.\ 89\,]\sp) that $d^{\,l\sn}f\sp(x)$
is a continuous multilinear map $\varPi^{\,l}\to F$ for all $l\in i+1\,$ and $
x\in\varOmega\ssp$. Even more, \hfill the variation \hfill $\delta^{\,l\sn}f$ \hfill
is continuous \hfill $\taurd\sp\varPi^{\,l\ssp+\sp 1}\to\taurd F\sp$, \hfill
whence writing $C_c^{\,i\sp}(\varOmega\,;\varPi_{\sp},F\sp)=(\sp
C_c^{\,i\sp}\image\sn\{\ssp(\varPi_{\sp},F\sp)\ssp\})\cap(\sp\vecs
F\sp)^{\,\varOmega}=\{\,f:(\varPi_{\sp},F_{\sp},f\sp)\in C_c^{\,i}$ and $\dom
f=\varOmega\,\}\ssp$,\Newline see the three lines paragraph just before
Proposition 3.4 below, {\it we would get\ssp} $\ssp\vecs E$ $=
C_c^{\,i\sp}(\varOmega\,;\varPi_{\sp},F\sp)$ \hfill {\it in the case where a
positive answer to\ssp} Question ? {\it exists\,}, {\it and\,} we have $\varPi
$ {\it metrizable and barrelled\ssp}.

Observe that for \hbox{$\biit K=\biit C$} a positive answer to Question ? is
obtained from the discussion just after formula (d) in the proof of Theorem
3.8 below. For real scalars, we trivially have a positive answer if $
\taurd\sp\varPi$ is the strongest locally convex topology for $\sigrd\sp\varPi
\sp$. Then also $\varPi$ is barrelled by [\,Jr\ssp; Prop.\ 11.3.1\sp(b)\ssp,
p.\ 223\,] since it is linearly homeomorphic to the locally convex direct sum,
or coproduct $\coprod{_{_{\sn\roman{LCS}}}}(\sp\alpha\sn\times\sn
\{\biit R\ssp\})$ where $\alpha=\roman{dim_{_H}}\sp\varPi\sp$.
However, it is metrizable only if $\sp\alpha\in\No\ssp$.

An example with $\varPi$ real and normable, but not even barrelled, where we
have a negative answer to \hbox{Question ?} can be constructed as follows.
Assuming $1<p\le\plusinfty\,$, \hfill we take $\,F=\biit R\,$ and $\,\varPi=
(X\sp,\cal T\sp)\,$ where $\,X=\coprod{_{\!_{\roman{VS}}}}(\No\ssn\times\sn
\{\sp\bold R\sp\})\,$ and $\,\cal T\,$ is the topology induced from $
\ell\RHB{.3}{^{\,p\sp}}(\No)\ssp$. \hfill With $\,\varphi=
\seq{\,s\,(1+s^{\,2\sp})^{\sp-1}\sn:s\in\Re\,}\ssp$, \,we define $f\sp(x)=
\sum_{\sp n\sp\in\sp I\!\!N_{\roman o}}
2^{\sp-n\sp}\varphi\ssp(\sp 2^{\,n\sp}x{_{}}_n)\,$ for $\,x=\seq{\,x{_{}}_n:n
\in\No\ssp}\in\vecs\varPi\sp$. \hfill One easily verifies that we have $f$ and
$d{_{}}_u\sp f:x\mapsto\sum_{\sp n\sp\in\sp I\!\!N_{\roman o}}
\varphi\ssp'(\sp 2^{\,n\sp}x{_{}}_n)\,u{_{}}_n$ continuous $\ssp\cal T\to
\tauR{\ }$ for every\Newline $u=\seq{\,u{_{}}_n\ssp}\in\vecs\varPi\sp$. \hfill
We see that $f\sp'(\sp\bnull\varPii)\not\in\cal L\ssp(\varPi_{\sp},F\sp)$
from $f\sp'(\sp\bnull\varPii)\ssp u=\sum_{\sp n\sp\in\sp I\!\!N_{\roman o}}
u{_{}}_n$ by considering $u=\seq{\,(\sp n+1)^{\sp-1}\sn:n\in\smb N\,}\cup(
(\No\ssn\setminus\smb N\sp)\sn\times\sn\{0\})$ with $\smb N\to\infty\,$.

A positive answer to \hbox{Question ?} for real scalars is obtained if we
assume that \œ$\varPi\in\roman{LCS}$ is Baire, e.g., a Fr\'echet space.
Indeed, fixing \œ$x\in\varOmega\ssp$, and for an arbitrary continuous seminorm
$\Nu$ for $F$ considering the function

$\varphi=\seq{\,\sup\ssp\{\,\Nu\,(\sp n\,(\sp f\sp\value z-f\sp\value x\sp)):z
=x+n^{\sp-1\sp}u\in\varOmega$ and $n\in\Zepp\ssp\}:u\in\vecs\varPi\,}\ssp$,

\noin one first proves $(\sp\varphi\inve\sp)\image{]}\,n\ssp,\plusinfty\,{[}
\in\taurd\sp\varPi$ for every $n\in\Zepp$. The Baire assumption then gives $\,
\roman{Int_{\ssp}}_{{\tau_{}}_{rd}\sp\varPi}\ssp((\sp\varphi\inve\sp)\image
[\,0\,,n\,]\sp)\not=\emptyset\,$ for some $\ssp n\ssp$. \hfill It follows that
$\,f\sp'(x)\ssp[\,U\,]\inc\mathbreak(\sp\Nu\sp\inve\sp)\image[\,0\,,1\,{[}\sp$
holds for some $U\in\ymp\varPi\sp$, whence the conclusion. From our discussion
above we further get \hfill $\vecs C^{\ssp i\sp}(\varOmega\ai{\varPii}\sp,
F\sp)=C_c^{\,i\sp}(\varOmega\,;\varPi_{\sp},F\sp)$ {\it for any $\varPi\sp,F
\in\roman{LCS}\ssp(\biit R\sp)$ and $\varOmega\in\taurd\sp\varPi$ with $\varPi
$ metrizable Baire\,}.

\NSN{\bf
5} \Examples (a) \ Let $\varPi\sp,E\in\roman{BaS}\ssp(\biit R\sp)$ be such
  that $\varPi$ is finite dimensional, and let $Q=\bigcup\,(\rng\bit{024}\sp)
\ssp$, where $\bit{024}=\seq{\,{K_{}}_i\sn:i\in\No\ssp}$ is increasing, and $
K{_{}}_i$ is $\taurd\varPi\,$--\,compact with $K{_{}}_i=
\roman{Cl_{\ssp}}_{{\tau_{}}_{rd}\sp\varPi}\sp(\sp
\roman{Int_{\ssp}}_{{\tau_{}}_{rd}\sp\varPi}\ssp{K_{}}_i\sp)$ for $i\in
I\sn\ar 1=\No\ssp$. Assume that every $(\sp\taurd\varPi\sp)\lei Q\,$--\,
compact $K$ has some $i$ with $K\inc{K_{}}_i\ssp$. \hfil Then we write $F=
\Cinfty(\sp Q\ai\varPii\sp,E\sp)$ and $F{_{\!}}_i=\mathbreak C^{\ssp i}((
{K_{}}_i)\sp\ai\varPii\sp,E\sp)$ \hfil for \hfil $i\in\No\ssp$. \hfil With \hfil
$\Delta=\{\ssp(\sp i\ssp,\sp j\sp):i\ssp,\sp j\in\No$ and $i\le j\,\}\ssp,$ \hfil
we \hfil also \hfil let\linebreak ${\rho_{}}_{ij}=\seq{\,y\,|\,{K_{}}_i\sn:y
\in\vecs F_j\ssp}$ and ${\rho_{}}_i=\seq{\,y\,|\,{K_{}}_i\sn:y\in\vecs F\,}$
for $(\sp i\ssp,\sp j\sp)\in\Delta\,$. Then the pro- jective system \hbox{$
\{\ssp(\sp i\ssp,\sp j\,;F_j\ssp,F{_{\!}}_i\ssp,{\rho_{}}_{ij}):
(\sp i\ssp,\sp j\sp)\in\Delta\,\}$} has limit \hbox{$\seq{\,(F\sp,
F{_{\!}}_i\ssp,{\rho_{}}_i):i\in\No\ssp}\ssp$}. The space $F$ is Fr\'echet,
and each $F{_{\!}}_i$ is Banach. In case $\bit{024}$ is constant, we have $
{\rho_{}}_{ij}\cup{\rho_{}}_i\inc\roman{id}\,$ for $\,(\sp i\ssp,\sp j\sp)\in
\Delta\,,$ and then $F=\leLCS-\sup\,\{\ssp F{_{\!}}_i\sn:i\in\No\ssp\}\ssp$.

(b) \ Let $\emptyset\not=\cal F\inc\dom\sn({}{\le}\sp\ai X)$ for some $\bold K
\,$--\,vector space $X\sp,$ and let $S$ be a vector subspace of $X\sp$. If now \hbox{$
G=\leLCS-\sup\ssp\cal F_{\sp}$}, then \hbox{$G_{/S}=\leLCS-\sup\,\{\ssp F_{/S}
:F\in\cal F\,\}\ssp$}. If also $\cal F$ is right $\leLCSr\sp$--\,directed, i.e.,
we have $\Delta=\leLCSr\sn\cap\cal F^{\times2}$ a direction, then \hbox{$
\{\ssp(F\sn\ar 1\sp,F\ssn\ar 2\,;F\ssn\ar 2\ssp,F\sn\ar 1\sp,
\idv F\ssn\ar 2):(F\sn\ar 1\sp,F\ssn\ar 2)\in\Delta\,\}$} is an LCS
projective system with limit $\seq{\,(\sp G\sp,F\sp,\idv\sp G\sp):F\in\cal F\,
}\ssp$. For example, in (a) we might have $Q={K_{}}_i=[\,0\,,1\,]$ and $\varPi
=E=\biit R\,,$ and $X=\bold R^{\sp Q}$ with $\,S=\{\,x\in{I\!\!R}^{\sp\,Q\!}:
x\value 0=0\,\}$\newline\hyppy{46.2mm}or $\,S=\{\,x\in{I\!\!R}^{\sp\,Q\!}:
x\value 0=x\value 1=0\,\}\ssp$.

(c) \ We have $\cal D\ssp(\Re\sp)=\leLCS-\sup\,\{\ssp F_m:m\in\cal M\,\}\ssp,$
when the Banach spaces $F_m$\newline are defined as follows. We write $
{I_{}}_n=\Re\setminus[-n\ssp,n\,]$ for $n\in\N\sp,$ and put ${I_{}}_0=\Re\,$.
Letting $\cal M$ be the set of increasing sequences $m=\seq{\,{i_{}}_n\sp}:\No
\to\N\sp,$ for $m\in\cal M\ssp,$ we let $B_m$ be the set of all functions $x:
\Re\to\Re$ such that for all $n\in\No$ we have $x\,|\,{I_{}}_n\in\vecs
C^{\ssp i_n}(\sp{I_{}}_n)$ and ${i_{}}_n\ssp|\,x^{(i)}(s)\ssp|\le 1$ for all $
i\in i_n^{\ssp+}$ and all $s\in I_n\ssp$. We now let $F_m$ be the unique
topological vector space such that $\sigrd F_m$ is a vector substructure\linebreak
of $\bold R^{\sp I\!\!R}$ and $B_m\in(\ymp F_m)\cap(\rajou F_m)\ssp$. \hfill
Then $F_m\in\roman{BaS}\ssp(\biit R\ssp)$ with $B_m$ the closed unit ball for the
norm $\vecs F_m\owns x\mapsto\|\ssp x\ssp\|\LHB{.4}{\sp_m}=\inf\ssp\{\,t\in
\Rep\ssn:t^{-1}x\in{B_{}}_m\ssp\}\ssp$.

The description of a filter basis for $\ymp\cal D\ssp(\Re\sp)$ given in
[\,Ho\ssp; Exa.\ 2.12.7, pp.\ 170\,--\,171\ssp] should serve as a hint for the
reader to finding a proof of the assertion that $\cal D\ssp(\Re\sp)=
\leLCS-\sup\,\{\ssp F_m:m\in\cal M\,\}$ holds.

Alternatively, we could have defined the spaces $F_m$ by first constructing
the sets $B_m$ as follows. We let $\cal M$ be the set of all functions $m=
\seq{\,{i_{}}_n\sp}:\Ze\to\N\sp,$ and we write ${I_{}}_n={]}\,n\ssp,n+2\,{[}$
for $n\in\Ze\ssp$. For $m\in\cal M\ssp,$ we then let $B_m$ be the set of all $
x\in{I\!\!R}^{\,I\!\!R}$ such that $x\,|\,{I_{}}_n\in\vecs C^{\ssp i_n}(\sp
{I_{}}_n)$ and ${i_{}}_n\ssp|\,x\sp\lupar\yyi i\rupar(s)\ssp|\le 1$ hold for
all $n\in\Ze$ and every $i\in i_n^{\ssp+}$ and $s\in{I_{}}_n\ssp$.

\NS Most of our later examples concern maps which are particular cases of

\NSN{\bf
6} \Examplee Let $i\in\No\cup\{\infty\}\ssp,$ and let $\varPi\ssp,
  \varPi\ar{\sn 1}\sp,\varPi\ar{\sn 2}\in\roman{LCS}\ssp(\biit R\sp)$ be
finite dimensio- nal. Let $Q\inc\vecs\varPi$ be as specified in paragraph 3
above. Writing $Q\yr 0=Q\times(\sp\vecs\varPi\ar{\sn 1})$ and $\varPi\ar{\sn 0}
=\varPi\sqcap\varPi\ar{\sn 1}\ssp,$ let $O\in(\sp\taurd\varPi\ar{\sn 0})\lei
Q\yr 0$ be such that $\dom(\sp Q\yr 0\ssn\setminus\sn O\sp)$ is relatively $
(\sp\taurd\varPi\sp)\lei Q\,$--\,compact (intuitively: the "hole\ssp" in $O$
projected on $Q$ should be rela-\Newline tively compact)\ssp. Then defining \hbox{$
E=\Cinfty(\sp O\ai\varPii\LHB{.3}{_{\sn_0}}\sp,\varPi\ar{\sn 2})$} and \hbox{$
F\sn\ar 1=C^{\ssp i\sp}(\sp Q\ai\varPii\ssp,\varPi\ar{\sn 1})$} and $
F\sn\ar 2 =C^{\ssp i\sp}(\sp Q\ai\varPii\ssp,\varPi\ar{\sn 2})\ssp,$ we are
interested in the map $\tilde f=(E\sqcap F\sn\ar 1\ssp,F\sn\ar 2\ssp,f\sp)
\ssp$, where\vskip.5mm\centerline{$
f=\{\ssp(\sp x\ssp,y\ssp,z\sp):x\in\vecs E$ and $y\in\vecs F\sn\ar 1$ and $
                 z=x\circ[\,\roman{id\,},y\,]\in\vecs F\sn\ar 2\ssp\}\ssp$.}

The "hard part\ssp" of the proof of Theorem 3.6 below is contained in

\newProCla 7 Lemma.

In {\ssp\rm Example 6\ssp,} whenever $(\sp x\ssp,y\sp)\in\dom f\sp,h\in\vecs
F\sn\ar 1\ssp,W\in\ymp F\sn\ar 2\ssp,$ there are $U_{},V$ with $U\sn\times\sn
V\in\ymp(E\sqcap F\sn\ar 1)$ and $\delta>0$ such that for all $(\sp u\ssp,v\sp
)\in U\sn\times\sn V$ and\Newline $0<|\ssp t\ssp|<\delta$ hold the following\vskip.5mm\noin{\rm
(1)} \ $u\circ[\,\roman{id\,},y+v\,]\in W\sp,$\newline{\rm
(2)} \ $x\circ[\,\roman{id\,},y+v\,]-x\circ[\,\roman{id\,},y\,]\in W\sp,$\newline{\rm
(3)} \ $t^{-1}\sp(\sp x\circ[\,\roman{id\,},y+t\,h\,]-x\circ[\,\roman{id\,},y\,]
      \sp)-\partial\ar 2\ssp x\circ[\,\roman{id\,},y\,]\,.\,h\in W\sp$, \hfill {\,\rm
      where $\,\partial\ar 2\ssp x=$\linebreak\null \hfill $
\{\ssp(\sp\eta\ssp,\xi\ssp,p\sp):(\sp\eta\ssp,\xi\sp)\in O$ and $p=\seq{\,
\delta\ai\varPii\sn\LHB{.3}{_{_1}}\ai\varPii\sn\LHB{.3}{_{_2}}(\sp
x\ssp(\sp\eta\ssp,\cdot))\ssp(\sp\xi\ssp,\zeta\sp):\zeta\in\vecs
                                                 \varPi\ar{\sn 1}\ssp}\ssp\}\,$.}

\Prooff Fixing a norm in each space $\varPi\ssp,\varPi\ar{\sn 1}\sp,
\varPi\sn\ar 2\ssp,$ a basis of the filter $\ymp F\sn\ar 1$ is $\{\,
N\sn\ar 1\sp(K\sp,\eps\sp):K$ is $(\sp\taurd\varPi\sp)\lei Q\,$--\,compact and
$\eps>0\,\}\ssp,$ when we let\par$\mhyppy{3}
N\sn\ar 1\sp(K\sp,\eps\sp)=\{\,y\in\vecs F\sn\ar 1\ssn:\all{k\ssp,\bit\eeta}\,
k\in\No$ and $k<i+1$\newline\hyppy{26mm} and $k\,\eps<1$ and $\bit\eeta\in
(\bar B\ai\varPii\sp(1))^{\,k}\imply(\sp\bar d\ssp(\bit\eeta)\ssp
y\sp)\sp\image\sn K\inc B\ai\varPii\sn\LHB{.3}{_{_{1\!}}}(\eps)\ssp\}\,,$

\noin and similar characterizations hold for $E$ and $F\sn\ar 2\ssp$.

If we have $(\sp x\ssp,y\ssp,z\sp)\in f\sp,$ then $z=x\circ[\,\roman{id\,},
y\,]\in\vecs C^{\ssp i\sp}(\sp Q_{\varPi\,},\varPi_2)\ssp,$ whence $\dom z=Q=
\dom y\ssp,$ and so $y=\rng[\,\roman{id\,},y\,]\inc\dom x=O\ssp$. Conversely,
we have $w\in\dom f$ for any $w=(\sp x\ar 1\sp,y\ar 1)\in\vecs(E\sqcap
F\sn\ar 1)$ with $y\ar 1\inc O\ssp$. Now $\dom(\sp Q\yr 0\ssn\setminus\sn
O\sp)\inc K\ar 0$ for some\Newline compact $K\ar 0\inc Q\ssp$. Then so is $
y\,|\,K\ar 0\inc O\ssp$. A standard argument gives some $\eps\ar 0\in$ ${]}\,
0\,,1\ssp{[}$ such that $(\sp y+v\sp)\,|\,K\sn\ar 0\inc O$ for all $v\in
N\sn\ar 1(K\ar 0\ssp,2\,\eps\ar 0)\ssp$. But then we have also $y+v\inc O\ssp
$, whence $(\sp x\ar 1\sp,y+v\sp)\in\dom f$ for every $x\ar 1\in\vecs E\ssp$. \hfill
Keeping this result in mind, below we tacitly assume $K\ar 0\inc K$ and $
\eps\le\eps\ar 0\ssp$.

Recall the rule for change of order of differentiation when the (partial)
derivatives are continuous, theorem of H.\ A.\ Schwarz, which already was
implicitly used in the symmetry proof of [\,Ke\ssp; Thm.\ 2.4.0, p.\ 90\,]\ssp.
Using this, for each fixed $\bit\eeta\in(\sp\vecs\varPi\sp)^{\,k}$ with $k\in
i+1\ssp$, a recursion on $l\in k=\dom\bit\eeta$ in conjunction with an
associated inductive proof shows existence of a finite function $S=
S\sn\ar 1(\bit\eeta)$ with the following properties\par{\leftskip6mm\hskip-10.5mm\noin%
(a) \ For every \hbox{$\pi\in S$} there are a bijection \hbox{$\sigma:k\to k$}
   and \hbox{$n\in\Zepp$} and finite sequen- ces $\bit\eeta\yr 0$ and \hbox{$
\bit\eeta\yr 1=\seq{\,\bit\eeta\yrai^1_1\sp,\ldots\,\bit\eeta\yrai^1_{i_1}\ssp
}$} and \hbox{$\bit\eeta\yr 2=\seq{\,\bit\eeta\yrai^2_1\sp,\ldots\,
\bit\eeta\yrai^2_{i_2}\ssp}\ssp$}, where further each $\bit\eeta^{\,\iota}_j$
is a finite sequence, and such that we have \hbox{$\pi=(\sp\tau_{\sp},n\sp)=
(\sp\bit\eeta\yr 0\!,\bit\eeta\yr 1\!,\bit\eeta\yr 2\!,n\sp)$} and\linebreak \hbox{$
\bit\eeta^{\,0}\sp\concc\bit\eeta\yrai^1_1\sp\conc\!\ldots\,
\bit\eeta\yrai^1_{i_1}\ssn\concc\bit\eeta\yrai^2_1\sp\concc\!\ldots\,
\bit\eeta\yrai^2_{i_2}=\bit\eeta\ssp\circ\sp\sigma\,$}. Hence we then here
have $\tau=(\sp\bit\eeta\yr 0\!,\bit\eeta\yr 1\!,\bit\eeta\yr 2)\in\dom S=
D\sn\ar 1(\bit\eeta)$ with $n=S\sp\value\tau=\smb N{_{\sn}}_\tau\in\Zepp$. \hfill
We also have the estimate\linebreak $\sum_{\,\tau\sp\in\sp D_1(\bmi7\eeta)}
\smb N{_{\sn}}_\tau=\sum S\sn\ar 1(\bit\eeta)\le(\sp k+1)^{\,k}$.\vskip0mm\hskip-10mm\noin%
(b) \ If for $\tau$ as above and $i\ar 0=i\ar 1\sn+i\ar 2$ and\newline$\mhyppy{10}
  P_{\tau\sp}(\sp y\ssp,v\sp)=[\,\bar d\ssp(\sp\bit\eeta\yrai^1_1)\ssp y\ssp,
\ldots\,\bar d\ssp(\sp\bit\eeta\yrai^1_{i_1})\ssp y\ssp,\bar d\ssp(\sp
\bit\eeta\yrai^2_1)\ssp v\ssp,\ldots\,
\bar d\ssp(\sp\bit\eeta\yrai^2_{i_2})\ssp v\,]\,,$ \hfill we write\linebreak$
\roman{t_{}}_{\tau\sp}(\sp u\,;y\ssp,v\sp)=
(\sp\partial\ar 2^{\,i_0}\bar d\ar 1^{}(\sp\bit\eeta\yr 0)\ssp u\sp)\circ
[\,\roman{id\,},y+v\,]\,.\,P_{\tau\sp}(\sp y\ssp,v\sp)\,,$ \hfill
where further\linebreak$\mhyppy{6.4}
\partial\ar 2^{\,l}\ssp\bar d\ar 1^{}(\bit\eeta\yr 0)\ssp u=
\partial\ar 2^{\,l}\sp(\sp\bar d\ar 1(\bit\eeta\yr 0)\ssp u\sp)$ \hfill with\linebreak$\mhyppy{10.4}
\bar d\ar 1(\bit\eeta\yr 0)\ssp u=\{\ssp(\sp\eta\ssp,\xi\ssp,
(\sp\bar d\ssp(\bit\eeta\yr 0)\ssp(\sp u\ssp(\,\cdot\,,\xi\sp)))\value\eta\sp)
:(\sp\eta\ssp,\xi\sp)\in O\,\}\,$ \hfill and\linebreak$\mhyppy{16.8}
\partial\ar 2^{\,l}\ssp\bar u=\{\ssp(\sp\eta\ssp,\xi\ssp,\sp\ell\ssp):
(\sp\eta\ssp,\xi\sp)\in O$ and\newline$\mhyppy{38.7}\ell\sp=\seq{\,
(\sp d\ssp(\bit\xxi)\ssp(\sp\bar u\ssp(\sp\eta\ssp,\sp\cdot\sp)))\value\xi:
\bit\xxi\in(\sp\vecs\varPi\ar{\sn 1})^{\,l}\ssp}\ssp\}\,$,\newline then $\mhyppy{10}
\bar d\ssp(\bit\eeta)\ssp(\sp u\circ[\,\roman{id\,},y+v\,]\sp)=
\sum_{\,\tau\sp\in\sp D_1(\bmi7\eeta)}\smb N{_{\sn}}_\tau\ssp
\roman{t_{}}_{\tau\sp}(\sp u\,;y\ssp,v\sp)$\vskip0mm
\noin
holds for any $\,u\in\vecs E$ and $y\ssp,v\in\vecs F\sn\ar 1$ with $y+v\inc O\ssp$.\par}

\noin Observe that for example in the case where \hbox{$\tau=
(\sp\bit\eeta\,,\emptyset\,,\emptyset\sp)$} we above have \hbox{$
\smb N{_{\sn}}_\tau=1$} and $\roman{t_{}}_{\tau\sp}(\sp u\,;y\ssp,v\sp)=
(\sp\partial\ar 2^{\,0}\bar d\ar 1^{}(\bit\eeta)\ssp u\sp)\circ
[\,\roman{id\,},y+v\,]\,.\,\emptyset\ai Q=(\sp\bar d\ar 1(\bit\eeta)\ssp u\sp)
\circ[\,\roman{id\,},y+v\,]\ssp$.

Now to prove (1)\ssp, fixing $y$ and $W_{}$, we have some compact $K$ and $
\eps>0$ with $N\sn\ar 2\sp(K\sp,\eps\sp)\inc W_{}$. \hfil Putting $V=
N\sn\ar 1\sp(K\sp,\eps\sp)\ssp,$ \hfil there is some $(\sp\taurd\varPi\sn\ar 0)
\lei O\,$--\,compact $K\sn\ar 1$\linebreak with $(\sp y+v\sp)\,|\,K\inc
K\sn\ar 1$ for all $v\in V_{}$. \hfil Taking \hfil $U=N\sp(K\sn\ar 1\sp,
\smb M^{\ssp-\eps^{-1\!}}\smb N^{\ssp-1}\eps\sp)\ssp,$ \hfil where\par\noin
with $\mhyppy{10}N=\bigcup\ssp\{\ssp(B\ai\varPii\sp(1))^{\,k\!}:k\in\No$
and $k<i+1$ and $k\,\eps<1\,\}$\newline
we have $\mhyppy{5}\smb M=\sup\ssp\{\,1+|\,(\sp\bar d\ssp(\bit\eeta)\ssp y\sp)
\value\eta\,|\ar 1\ssn:\bit\eeta\in N$ and $\eta\in K\,\}$\newline
and $\mhyppy{12}\smb N=\sup\big\{\sum S\sn\ar 1(\bit\eeta):\bit\eeta\in N\sp
\big\}\ssp,$ \hfill by a computation left to the reader,

\noin we obtain $\{\,u\circ[\,\roman{id\,},y+v\,]:(\sp u\ssp,v\sp)\in U\sn
\times\sn V\,\}\inc N\sn\ar 2\sp(K\sp,\eps\sp)\inc W_{}$, hence (1)\ssp.

For (2) fixing $x\ssp,y\ssp,W$ with $N\sn\ar 2\sp(K\sp,\eps\sp)\inc W_{}$, we
consider Riemann integration of continuous curves from $I=[\,0\,,1\,]$ to the
(generally) nonnormable locally convex space $\varPi\ar{\sn 2}\KN{.5}^Q\le
F\ssn\ar 2\ssp$. For $v\in N\sn\ar 1(K\ar{\sn 0}\ssp,\eps\sp)\ssp$, taking $
\bar v=\seq{\,\partial\ar 2\ssp x\circ[\,\roman{id\,},y+s\,v\,]\,.\,v:s\in I\,
}$ we may write $x\circ[\,\roman{id\,},y+v\,]-x\circ[\,\roman{id\,},y\,]=
\int_{\sp I}\bar v\,$. For any $\bit\eeta\in(\vecs\varPi\sp)^{\,k}$ with $k\in
i+1$ and \hbox{$\eta\in K\inc Q$} writing \hbox{$\check v=
\check v\,(\sp\bit\eeta\,,\eta\sp)=\seq{\,\bar d\ssp(\bit\eeta)\ssp(\sp
\bar v\ssp(s))\value\eta:s\in I\,}\ssp$}, by classical results we  have $
\bar d\ssp(\bit\eeta)\sp\big(\sn\int_{\ssp I}\bar v\ssp\big)\sn\value\eta=
\int_{\ssp I}\check v\,$, where the latter integral is for a curve in $
\varPi\ar{\sn 2}\ssp$.

Since $\check v$ is continuous, in order to have \hbox{$\int_{\ssp I}\check v
\in B\ai\varPii\sn\LHB{.3}{_{_2}}(\eps)\ssp$}, it suffices to establish \hbox{$
\rng\check v\inc B\ai\varPii\sn\LHB{.3}{_{_2}}(\eps)\ssp$}. Altogether, we see
that for the proof of (2) it suffices to find $V$ with $V\ssn\ar 2=\{\,\bar v
\ssp(s):v\in V$ and $0\le s\le 1\,\}$ $\inc N\sn\ar 2\sp(K\sp,\eps\sp)\ssp$.

Defining $x\ar 2\ssp(\sp\eta\ssp,\xi\ar 1\sp,\xi\ar 2)=\partial\ar 2\ssp x\ssp
(\sp\eta\ssp,y\ssp(\eta)+\xi\ar 1)\,\xi\ar 2\ssp,$ we have $\bar v\ssp(s)=
x\ar 2\sn\circ[\,\roman{id\,},s\,v\ssp,v\,]$ and observe that $x\ar 2\sp(\sp
\eta\ssp,\xi\ar 1\sp,\,\cdot\,)$ is linear. A deduction (left to the reader)
similar to the one used to establish (a) and (b) above shows that we may write
\ $\bar d\ssp(\bit\eeta)\ssp(\sp\bar v\ssp(s))=$\vskip.5mm\centerline{$
\sum_{\,\tau\sp\in\sp D_2(\snn\bmi7\eeta)}s^{\,l\sp}\smb N\yrai^2_{\sn\tau}\ssp
(\sp\partial\ar 2^{\sp\,l}\bar d\ar 1^{}(\sp\bit\eeta\sp\ar 0)\ssp x\ar 2)
\circ[\,\roman{id\,},s\,v\ssp,\bar d\ssp(\sp\bit\eeta\sp\yr 0)\ssp v\,]\,.\,[\,
\bar d\ssp(\sp\bit\eeta\sp\yr 1)\ssp v\ssp,\ldots\,
\bar d\ssp(\sp\bit\eeta\,\yi l)\ssp v\,]\,$.}

\noin Consequently $V\ssn\ar 2\inc N\sn\ar 2\sp(K\sp,\eps\sp)$ follows for $V=
N\sn\ar 1\sp(K\sp,\smb M^{\ssp-1}\smb N^{\ssp-1}\eps\sp)\ssp$, when taking $N$
as\Newline before, and $\,P=\bigcup\ssp\{\ssp(
{B_{}}_{\varPi\ssp\sqcap\ssp\varPi_1\sp\sqcap\ssp\varPi_1}(1))^{\,k\!}:k\in\No
$ and $k<i+1$ and $k\,\eps<1\,\}\ssp$, \hfil we\linebreak put $\
\smb N=\sup\big\{\sum S\sn\ar 2(\bit\eeta):\bit\eeta\in N\sp\big\}$ \,and\newline$\mhyppy{7}
\smb M=\sup\ssp\{\,1+|\ssp(\sp\bar d\ssp(\bit{020}\sp)\ssp x\ar 2)\ssp(\sp
\eta\ssp,\xi\ar 1\sp,\xi\ar 2)\ssp|\ar 2\ssn:\bit{020}\in P$ and\newline$\mhyppy{46}
\eta\in K$ and $|\,\xi\ar 1|\ar 1\le\eps\ar 0$ and $|\,\xi\ar 2\sp|\ar 1\le 1
\,\}\ssp$. \hfil The norm\linebreak we are using for $\varPi\sqcap
\varPi\ar{\sn 1}\sn\sqcap\varPi\ar{\sn 1}$ is $(\sp\eta\ssp,\xi\ar 1\sp,\xi\ar 2)
\mapsto\sup\ssp\{\ssp|\ssp\eta\ssp|\ssp,|\,\xi\ar 1|\ar 1\sp,|\,\xi\ar 2\sp|\ar 1\}\ssp$.

To prove (3)\ssp, we first observe that the expression which should be got in
$W,$ can be written $x\ar 1\sn\circ[\,\roman{id\,},y\ar 1\sn+t\,h\,]-x\ar 1\sn
\circ[\,\roman{id\,},y\ar 1\ssp]\ssp,$ when we take $y\ar 1=
\bnull F\sn\LHB{.3}{_{_1}}$ and define\vskip.5mm\centerline{$
x\ar 1\sp(\sp\eta\ssp,\xi\sp)=\int_{\,0}^{\,1}(\sp\partial\ar 2\ssp x\ssp(\sp
\eta\ssp,y\ssp(\eta)+s\,\xi\sp)-\partial\ar 2\ssp x\ssp(\sp\eta\ssp,y\ssp(\eta
)))\ssp(h\ssp(\eta))\,d\ssp s\,$.}

\noin Noticing that (2) requires only derivatives of $x\ssp(\,\cdot\,,\xi\sp)$
of order less than $i+1\ssp,$ we can apply (2) with $x=x\ar 1$ and $y=y\ar 1$
and $v=t\,h\ssp$. This gives (3)\ssp.                                    \QED


\subhead3

                   Differentiability in general dimensions

We shall use the classes $C^{\ssp k}$ of order $k$ continuous
differentiabilities of [\,Hi1\,] "pulled back\ssp" to LCS\ssp, which are
exactly the Keller classes $\CPi k$ of [\,Ke\,]\ssp. Put precisely, we have $
\CPi k=\{\,\tilde f:\tilde f$ vector map (of l.c.s.) and $
\bit{025}\sp_{_{\bold{tv}}}\sp\tilde f\in C^{\ssp k\,}\}\ssp,$ when $
\bit{025}\sp_{_{\bold{tv}}}\sp(E\ssp,F\sp,f\sp)=(\sp
\lambda\sp_{_{\bold{tv}}}\sp E\ssp,\lambda\sp_{_{\bold{tv}}}\sp F\sp,f\sp)$
with $\lambda\sp_{_{\bold{tv}}}\sp(X\sp,\cal T\sp)=(X\sp,\varLambda\sp)\ssp,$
where in turn $\varLambda$ is the convergence (= $\lambda$imit structure)
corresponding to the topology $\cal T$. Hence we have $(\sp x\ssp,\ecal F\sp)
\in\varLambda$ if{}f $x\in T=\bigcup\cal T$ and $\ecal F$ is a filter on $T$
with $\{\,V:\exi{U\in\cal T}\,x\in U\inc V\inc T\,\}\inc\ecal F\sp$.

For a map $\tilde f=(E,F,f)$ we may also write $f:E\iinc U\to F$ when $U=\dom
f\sp$. We say $\tilde f$ {\it smooth\ssp} if{}f $\tilde f\in\CinftyPi$ holds.
We remark that although in [\,Ke\,] only real scalars have been considered,
Keller's definitions and results relevant for us are equally valid for $\biit
K=\biit C\ssp$. We shall write $\CPi k=\CPi k\ssp(\biit K\ssp)$ in case we
want to stress the coefficient field.

We say that $\tilde f$ is {\it almost Gateaux differentiable\ssp} (or {\it
directionally\ssp} differentiable, or has directional derivatives) if{}f $\dom
\sn(\sp\delta\ai{\sp E\sp F}\sp f\sp)=(\dom\sn f\sp)\sn\times\sn(\sp\vecs E\sp
)\ssp$. With complex scalars assuming that $F$ is Mackey complete (see [\,KM\ssp;
p.\ 15, Lemma 2.2\,]\sp) we then say that $\tilde f$ is {\it holomorphic\ssp},
and write $\tilde f\in\cal H_{_T}$ in case it is directionally diffe- rentiable
and $\tilde f\in\CPi{\ssp 0}\sp(\biit C\ssp)$ holds. By [\,KM\ssp; p.\ 81, Thm.\
7.4, $(2)\imply(1)$\,]\ssp, direction- al differentiability follows if the
function $t\mapsto\ell\,(\sp f\sp(\sp x+t\,u\sp))$ is an ordinary holo- morphic
function $\{\,t:x+t\,u\in\dom\ssn f\,\}\to\Ce$ for any $x\in\dom f$ and $u\in
\vecs E$ and $\ell\in\cal L\ssp(F\sp,\biit C\,)\ssp$. In case also $E$ is
Fr\'echet, our concept of holomorphy equals [\,KM\ssp; p.\ 83, Def.\ 7.8\,]\ssp,
see [\,KM\ssp; p.\ 88, Thm.\ 7.19, $(1)\equivv(3)$\,]\ssp, cf.\ also [\,Pi\ssp;
p.\ 183 ({\font\=cmr8\T})\,]\ssp. 

The higher order {\it differential\sp}\,functions $
d^{\,k}\KN{2.1}\LHB{.3}{\ai{EF}}\sp f=d^{\,k\sn}f$ of map $\tilde f=(E\ssp,
F\sp,f\sp)$ we define\Newline as follows. Letting $\varOmega$ be the set all $
x$ with $\{\ssp\seq{\sp x\sp}\conc\biit u:\biit u\in(\sp\vecs E\sp)^{\,k\,}\}
\inc\dom(\sp\delta^{\,k\sn}f\sp)\ssp$, and such that the function $\biit u
\mapsto\delta^{\ssp k\sn}f\sp(\sn\seq{\sp x\sp}\conc\biit u\sp)$ is
multilinear $(\sp\sigrd E\sp)^{\,k}\to\sigrd F\sp$, we put $d^{\ssp k\sn}f=
\seq{\ssp\seq{\,\delta^{\ssp k\sn}f\sp(\sn\seq{\sp x\sp}\conc\biit u\sp):
\biit u\in(\sp\vecs E\sp)^{\,k\,}}:x\in\varOmega\,}\ssp$.

For $k\in\No\cup\{\infty\}$ and a vector map $\tilde f=(E,F,f)$ with $E,F\in$
LCS and $\varOmega=\dom f\in\taurd E$, we now have $\tilde f\in\CPi k$ if{}f $
\dom(d^if)=\varOmega$ and $\rng(d^if)\inc\cal L^i(E,F)$, set of continuous
multilinear maps, hold for all $i<k+1$, and further we have\vskip.5mm
$\KP{10}\all{l<k+1\ssp,x\in\varOmega\ssp,W\in\ymp F}\,\exi{U\in\ymp E}$

$\KP{10}\all{\eps>0}\,\exi{V\in\ymp E}\,\all{v\in V}\,(\sp
 d^{\sp\,l}\sn f\sp(\sp x+v\sp)-d^{\sp\,l}\sn f\sp(x))\ssp[\,
 U^{\ssp l\ssp}]\inc\eps\ssp W\sp$.\hfill($\pi$)\vskip.5mm
\noin
The condition ($\pi$) is an explicit expression for $d^{\,l\sn} f$ to be
lim\,-\,continuous $\taurd(\lambda\sp_{_{\bold{tv}}}\sp E\sp)$ $\to\taurd
\Cal L^{\,l\sp}_\pi(E\ssp,F\sp)\ssp,$ see [\,Ke\ssp; p.\ 40, Prop.\ 0.5.1 and
p.\ 45, Thm.\ 0.5.12\,]\ssp.

\newProCla 1 Proposition.

If $(E\ssp,F\sp,f\sp)\in\CPi k$ and $F\sn\ar 1$ is a sequentially closed
topological lin- ear subspace of $F\sp,$ then the implication $\rng f\inc\vecs
F\sn\ar 1\imply(E\ssp,F\sn\ar 1\sp,f\sp)\in\CPi k$ holds.

\Prooff First using sequential closedness and induction on $l<k+1\ssp,$ we
prove $\rng(\sp\delta^{\,l\sn}f\sp)\inc\vecs F\sn\ar 1\ssp,$ and then apply
($\pi$)\ssp.                                                          \newQED

\noin{\bf
2} {\font\=cmssi10\D\sp e\sp f\sp i\sp n\sp i\sp t\sp i\sp o\sp n\ssp}.
  For a vector map $\tilde f=(E\ssp,F\sp,f\sp)\ssp,$ we say that $\cal F$ is a {\it
projective ext- ension\ssp} of $\tilde f$ via $\cal R$ if{}f we here have
functions $\cal R=\{\,(\sp i\ssp,\sp j\,;\sp\rho^{\,1}_{\sp i\sp j}\ssp,
\rho^{\,2}_{\sp i\sp j}\sp):(\sp i\ssp,\sp j\sp)\in\Delta\ssp\}$ and\Newline
$\cal F=\seq{\,({E_{}}_i\ssp,{F_{}}_i\ssp,{f_{}}_i\,;\sp\rho^{\,1}_{\sp i}\ssp,
\rho^{\,2}_{\sp i}\sp):i\in I\sn\ar 1\ssp}\ssp,$ where $I\sn\ar 1=\dom\Delta\,,
$ and $\seq{\,(F\sp,{F_{}}_i\ssp,\sp\rho^{\,2}_{\sp i}\sp):i\in I\sn\ar 1\ssp}
$ is an LCS projective limit of $\{\ssp(\sp i\ssp,\sp j\,;F_j\ssp,{F_{}}_i\ssp,
\sp\rho^{\,2}_{\sp i\sp j}\sp):(\sp i\ssp,\sp j\sp)\in\Delta\ssp\}\ssp,$ and\vskip.5mm\noin%
(1) \ $\rho^{\,1}_{\sp i}\in\cal L\ssp(E\ssp,{E_{}}_i\sp)$ and $
    \rho^{\,1}_{\sp i\sp j}$ linear $\sigrd E_j\to\sigrd E_i$ and $
    (E_i,F_i,f_i)$ a vector map\par\noin
(2) \ $\rho^{\,1}_{\sp i}=\rho^{\,1}_{\sp i\sp j}\circ\rho^{\,1}_{\sp j}$ and
    $\rho^2_i\circ f=f_i\circ\rho^1_i$ and $\rho^2_{ij}\circ f_j=f_i\circ\rho^1_{ij}$

\noin hold for all $(\sp i\ssp,\sp j\sp)\in\Delta\,$. We call $\cal F$ a
projective extension of $\tilde f$ if{}f it is such via some $\cal R\,,$ and
we speak of a projective $\CPi k\ssp$--\,extension in case $\roman{\,dom\,}(
\rng\cal F\sp)\inc\CPi k$ holds.

Observe that the projective limit assumption implicitly contains the
requirements $\rho^{\,2}_{\sp i\sp k}=\rho^{\,2}_{\sp i\sp j}\circ
\rho^{\,2}_{\sp j\sp k}$ and $\rho^{\,2}_{\sp i}=\rho^{\,2}_{\sp i\sp j}\circ
\rho^{\,2}_{\sp j}$ for $(i,j),(j,k)\in\Delta\,$. One also observes that if
in (2) we drop $\rho^2_i\circ f=f_i\circ\rho^1_i\ssp,$ the remaining two
equalities already guaran-\Newline tee existence of a unique $f$ with (2)\ssp.
Conversely, if instead of $\rho^2_{ij}\circ f_j=f_i\circ\rho^1_{ij}\ssp,$ we
take surjectivity $\rng\rho^1_i=\vecs E_i\ssp,$ then again (2) holds.
                                        
\newProCla 3 Theorem.

If $\tilde f$ has a projective $\CPi k\ssp$--\,extension$\,,$ then $\tilde f\in\CPi k\,$.

\Prooff With the notations of Definition 2\ssp, let $\cal F$ be a projective $
\CPi k\ssp$--\,extension of $\tilde f$ via $\cal R\,$. We shall show $\tilde f
\in\CPi k\,$. By induction on $l\ssp$, the reader easily verifies $
\rho^{\,2}_{\sp i\sp j}\circ(\sp\delta^{\,l\sn}f_j\sp)=(\sp\delta^{\,l\sn}
{f_{}}_i\sp)\circ\LHB{.3}{\roman X}\ssp(\ssp l^{\ssp+\!}\times\{\,
\rho^{\,1}_{\sp i\sp j}\ssp\})\ssp$. By induction on $l\ssp$, we next show $
\rho^{\,2}_{\sp i}\circ(\sp\delta^{\,l\sn}f\sp)$ $=(\sp\delta^{\,l\sn}
{f_{}}_i\sp)\circ\LHB{.3}{\roman X}\ssp(\ssp l^{\ssp+\!}\times\{\,
\rho^{\,1}_{\sp i}\ssp\})$ for all $i\in I_1\ssp$. Since each $d^{\,l\sn}
{f_{}}_i\sp(\sp\rho^{\,1}_{\sp i}\sp x\sp)$ is multilinear, this implies that
so is $d^{\,l\sn}f\sp(x)$ too. Here $l\in k+1\ssp$.

The case $l=0$ being our assumption, let the assertion hold for a fixed $l<
k\sp$. To prove it for $l+1\ssp$, we fix $x\in\dom f$ and $\biit w=\biit u\sp
\conc\seq{\sp h\sp}=\seq{\,u\ar 1\sp\ldots\,u\ssp\ai l\sp,h\,}\in(\sp\vecs
E\sp)^{\,l\ssp+\sp 1}$, and consider $\varDelta\ssp(t)=t^{-1}(\sp
d^{\,l\sn}f\sp(\sp x+t\,h\sp)-d^{\,l\sn}f\sp(x))\ssp\biit u\ssp$. We must find
some $y\in\vecs F$ with $
\taurd F\text{\,-\,}\lim_{\,t\to 0}\sp\varDelta\ssp(t)=y$ and $\rho^{\,2}_i\sp
y=y_i=d^{\,l\ssp+\sp 1}{f_{}}_i\sp(\sp\rho^{\,1}_{\sp i}\sp x\sp)\ssp(\sp
\rho^{\,1}_{\sp i}\circ\biit w\ssp)$ for all $i\ssp$.

We have $\rho^{\,2}_{ij}y_j=y_i$ for all $(i,j)\in\Delta\,$, whence there is a
unique $y$ with $\rho^2_iy=y_i\ssp$. To prove $\lim_{\,t\to 0}\sp\varDelta\ssp
(t)=y\ssp$, arbitrarily picking $W\in\ymp F\sp$, we have some $i\in I\ar 1$
with $W_i\in\ymp F_i$ and $(\sp\rho^{\,2}_i)\inve[\,W_i\sp]\inc W_{}$. Then
there is some $\delta>0$ with $\Delta^i(t)-y_i=t^{-1}(\sp d^{\,l\sn}f_i
(\sp\rho^1_ix+t\,\rho^1_ih\sp)(\sp\rho^1_i\circ\biit u\sp)-d^{\,l\sn}f_i
(\sp\rho^1_ix\sp)(\sp\rho^1_i\circ\biit u\sp))-y_i\in W_i$ for $0<|\ssp t\ssp|
<\delta\ssp$. By the inductive assumption, we have $\rho^{\,2}_i\sp(\varDelta
\ssp(t))=\Delta^i(t)\ssp$, whence for $0<|\ssp t\ssp|<\delta\ssp$, we obtain $
\varDelta\ssp(t)-y\in(\sp\rho^{\,2}_i)\inve[\,W_i\sp]\inc W_{}$, and the
induction is completed.

To show that the multilinear map $d^{\,l\sn}f\sp(x)$ is continuous, let $W\in
\ymp F\sp$. There is some $i\in I_1$ with $W_i\in\ymp F_i$ and $(\sp
\rho^{\,2}_i)\inve[\,W_i\sp]\inc W_{}$. Since $d^{\,l\sn}f_i(\sp\rho^1_ix\sp)$
is continuous, we have some \hfil $U_i\in\ymp E_i$ \hfil with \hfil $
d^{\,l\sn}f_i(\sp\rho^1_ix\sp)\ssp[\,U_i^{\,l}\ssp]\inc W_i\ssp$. \hfil Then \hfil
$U=(\sp\rho^{\,1}_i)\inve\sp[\,U_i\sp]\in\mathbreak\ymp E\ssp$, and to prove $
d^{\,l\sn}f\sp(x)\ssp[\,U^{\,l}\ssp]\inc W_{}$, arbitrarily picking $\biit u
\in U^{\,l}$, it suffices to have\Newline $\rho^{\,2}_i(\sp d^{\,l\sn}f\sp(x)
\ssp(\biit u))\in W_i\ssp$. \hfil This indeed holds, \hfil since \hfil $
\rho^{\,1}_i\circ\biit u\in U_i^{\,l}$, \hfil whence we obtain\linebreak $
\rho^{\,2}_i\sp(\sp d^{\,l\sn}f\sp(x)\ssp(\biit u))=d^{\,l\sn}f_i(\sp
\rho^{\,1}_i\sp x\sp)(\sp\rho^{\,1}_i\circ\biit u\sp)\in W_i\ssp$.

Now, to prove $\tilde f\in\CPi k\,$, we verify $(\pi)$ for it. Fixing $l<k+1$
and $x\in\dom f$ and $W\in\ymp F\sp$, we find some $W_i\in\ymp F_i$ with $
(\sp\rho^{\,2}_i)\inve[\,W_i\sp]\inc W_{}$. Then applying $(\pi)$ to the map $
(E_i,F_i,f_i)\ssp$, we take a corresponding $U_i\in\ymp E_i\ssp$, and put $U=
(\sp\rho^{\,1}_i)\inve[\,U_i\sp]\ssp$. Next we fix $\eps>0\,$, and still using
$(\pi)$ take a corresponding $V_i\in\ymp E_i\ssp$, and put $V=(\sp
\rho^{\,1}_i)\inve[\,V_i\sp]\ssp$. \hfil For $v\in V$ fixed, \hfil to prove \hfil
$(\sp d^{\,l\sn}f\sp(\sp x+v\sp)-d^{\,l\sn}f\sp(x))\ssp[\,U^{\,l}\ssp]\inc
\eps\,W_{}$,\linebreak arbitrarily picking $\biit u\in U^{\,l}$, it suffices
to show $\rho^{\,2}_i\sp((\sp d^{\,l\sn}f\sp(\sp x+v\sp)-d^{\,l\sn}f\sp(x))\ssp
\biit u\sp)\in\eps\,W_i\ssp$. This holds, since $\rho^{\,1}_iv\in V_i$ and $
\rho^{\,1}_i\circ\biit u\in U_i^{\,l}$, whence $\rho^{\,2}_i\sp((\sp
d^{\,l\sn}f\sp(\sp x+v\sp)-d^{\,l\sn}f\sp(x))\ssp\biit u\sp)$ $=(\sp
d^{\,l\sn}f_i(\sp\rho^{\,1}_ix+\rho^{\,1}_iv\sp)-d^{\,l\sn}f_i(\sp\rho^{\,1}_i
x\sp))\ssp(\sp\rho^{\,1}_i\circ\biit u\sp)\in\eps\,W_i\ssp$.          \newQED

The Keller class $C_c^{\ssp k}$ for $k\in\No\cup\{\infty\}$ has as members
exactly the maps $\tilde f=$ $(E\ssp,F\sp,f\sp)\in\CPi{\ssp 0}$ with $\dom\sn
f=\dom\sn(\sp d_{EF}^{\,i}f\sp)$ and $\delta^{\ssp i\sn}\tilde f\in\CPi{\ssp 0}
$ for all $i<k+1\ssp$. From [\,Ke\ssp; p.\ 99, Prop.\ 2.7.1\,] we obtain
\NSN{\bf
4 Proposition.} {\it
The equality $C_c^{\,\infty}=\CinftyPi$ is valid.}
\NS
We write $\tilde g\ssp\Circ\tilde f=(E\ssp,G\ssp,g\circ f\sp)$ for vector maps
$\tilde f=(E\ssp,F\sp,f\sp)$ and $\tilde g=(F\sp,G\ssp,\sp g\sp)\ssp$.\Newline
With the class $\cal L$ of continuous linear maps of locally convex spaces,
note $\cal L\inc\CinftyPi$ for the applications of the following

\newProCla 5 Theorem.

Let $\ssp\cal C\inc\CPi{\sp 0}$ be such that every $\tilde f\in\cal C$ is
almost Gateaux differentiable\Newline and has some $\sp\tilde g\in\cal C$
and some $\sp\tilde h\in C^{\ssp\infty}_c$ with $\sp\delta\sp\tilde f=\tilde g\ssp
\Circ\tilde h\,$. Then $\ssp\cal C\inc\CinftyPi\ssp$.

\Prooff If we have $(*)$ \ $\{\,\tilde f\in\cal C:\delta\sp\tilde f\in
C_c^{\ssp k\,}\}\inc C_c^{\,k\ssp+\sp 1}$ for all $k\in\No\ssp,$ \hfil using
the\linebreak chain rule for $C_c^{\ssp k}\ssp,$ a trivial induction gives $
\cal C\inc C_c^{\ssp k}$ for all $k\in\No\ssp,$ hence $\cal C\inc\CinftyPi$ by
Proposition 4\ssp. Thus it suffices to prove $(*)$ by induction on $k\in\No\ssp$.

To prove $(*)$ for $k=0\,,$ let $\tilde f=(E\ssp,F\sp,f\sp)\in\cal C$ with $
\delta\sp\tilde f\in C_c^{\sp\,0}=\CPi{\ssp 0}\,$. Since we have $\cal C\inc
\CPi{\ssp 0}\,,$ and $\delta\sp f\sp(\sp x\ssp,t\,u\sp)=t\,\delta\sp f\sp(\sp
x\ssp,u\sp)$ for all $x\in\dom f$ and $u\in\vecs E$ and $t\in\Ke\sp,$ to
obtain $\dom(\sp d^{\,1\sn}f\sp)=\dom f\sp,$ it suffices for fixed $x\in\dom f
$ and $u\ssp,v$ $\in\vecs E$ to prove $w=\delta\sp f\sp(\sp x\ssp,u+v\sp)-
\delta\sp f\sp(\sp x\ssp,u\sp)-\delta\sp f\sp(\sp x\ssp,v\sp)=\bnull F\ssp$.

To proceed via reductio ad absurdum, assuming $w\not=\bnull F\ssp,$ by Hahn\,--\,
Banach, we have some $\ell\in\cal L\ssp(F\sp,\biit K\,)$ with $\roman{Re\,}(\ssp
\ell\,w\sp)=1\ssp$. Then putting $g=\ell\circ f\sp,$ by continuity of $\delta\sp
f\sp,$ we have some $\delta>0$ with $\roman{Re\,}(\sp\delta\ssp g\ssp(\sp x+t\,
u+s\,v\ssp,v\sp)-\delta\ssp g\ssp(\sp x\ssp,v\sp))<\frac12$ for $0<$ $t\ssp,s<
\delta\ssp$. By the definition of $\delta\sp f\sp,$ we further have some $t\in
{]}\,0\,,\delta\,{[}$ with $\alpha=\roman{Re\,}($ $t^{-1\sp}(\sp g\ssp(\sp x+t
\,u\sp)-g\ssp(x))-\delta\ssp g\ssp(\sp x\ssp,u\sp)-t^{-1\sp}(g\ssp(\sp x+t\,
(\sp u+v\sp))-g\ssp(x))+\delta\ssp g\ssp(\sp x\ssp,u+v\sp))$ $<\frac12\ssp$.
By the mean value theorem of classical differential calculus, we have some $s
\in{]}\,0\,,t\,{[}$ with $t^{-1\sp}\roman{Re\,}(\sp g\ssp(\sp x+t\,(\sp u+v\sp
))-g\ssp(\sp x+t\,u\sp))=\roman{Re\,}(\sp\delta\ssp g\ssp(\sp x+t\,u+s\,v\ssp,
v\sp))\ssp$. Then we get \ $1=\roman{Re\,}(\ssp\ell\,w\sp)=\roman{Re\,}(\sp
\delta\ssp g\ssp(\sp x\ssp,u+v\sp)-\delta\ssp g\ssp(\sp x\ssp,u\sp)-\delta\ssp
g\ssp(\sp x\ssp,v\sp))$\newline$\mhyppy{5}
=\alpha+\roman{Re\,}(\sp t^{-1\sp}(\sp g\ssp(\sp x+t\,(\sp u+v\sp))-g\ssp(\sp
x+t\,u\sp))-\delta\ssp g\ssp(\sp x\ssp,v\sp))$\newline$\mhyppy{5}
=\alpha+\roman{Re\,}(\sp\delta\ssp g\ssp(\sp x+t\,u+s\,v\ssp,v\sp)-\delta\ssp
g\ssp(\sp x\ssp,v\sp))<\frac12+\frac12=1\ssp,$ a contradiction.

Now assuming $(*)$ for a fixed $k\in\No\ssp,$ to prove it for $k+1\ssp,$ let
$\tilde f\in\cal C$ with $\delta\sp\tilde f$ $\in C_c^{\,k\ssp+\sp 1}\ssp$. By
$C_c^{\,k\ssp+\sp 1}\inc C_c^{\ssp k}$ and the inductive assumption, we then
have $\tilde f\in C_c^{\,k\ssp+\sp 1}\ssp$. Hence to prove $\tilde f\in
C_c^{\,k\ssp+\sp 2}\ssp,$ it suffices to show $\dom\sn(\sp d^{\,k\ssp+\sp 2}f
\sp)=\dom\sn f$ and $\delta^{\,k\ssp+\sp 2}\tilde f\in\CPi{\ssp 0}\,$. Writing $
\delta^{\,l\sp}(\sp\delta\sp\tilde f\ssp)=((E\sp\yr 2)^{\,l\ssp+\sp 1\!},F\sp,
\delta^{\,l\sn}f\sn\ar 1)\ssp,$ by induction on $l\in k^{\,++},$ we see\vskip.5mm\centerline{$
\delta^{\,l\sn}f\sn\ar 1\sp\seq{\,(\sp x\ar 0\sp,x\ar 1)\ssp,
 (\sp x\ar 2\sp,\bnull E)\ssp,\ldots\,(\sp x\ssp\ai l\ar{\sn+1}\ssp,\bnull E)\,}
=\delta^{\,l\ssp+\sp 1\sn}f\ssp\seq{\,x\ar 0\sp,x\ar 1\sp,x\ar 2\sp,\ldots\,
 x\ssp\ai l\ar{\sn+1}\ssp}$}\vskip.5mm
\noin
for all $\seq{\,x\ai i\sp}\in\vecs E^{\,l\sp+\sp 2}$ with $x\ar 0\in\dom\sn f
\sp$. Taking $l=k+1\ssp,$ we immediately obtain $\delta^{\,k\ssp+\sp 2}
\tilde f\in\CPi{\ssp 0}\,,$ and also see that each function $x\sp\ai i\mapsto
\delta^{\,l\ssp+\sp 1\sn}f\ssp\seq{\,x\ar 0\sp,x\ar 1\sp,x\ar 2\sp,\ldots\,
x\ssp\ai l\ar{\sn+1}\ssp}$ is linear $\sigrd E\to\sigrd F$ for $i=2\ssp,\ldots\,
l+1\ssp$. Linearity with respect to $x\ar 1$ is easily seen to hold by
induction on $l\ssp,$ and then we are done.                              \QED

\newProCla 6 Theorem.

The membership $(E\sqcap F\sn\ar 1\ssp,F\sn\ar 2\ssp,f\sp)\in\CinftyPi$ holds
in {\rm\,Example 2.6\ssp.}

\Prooff Keeping $Q\ssp,\varPi\ssp,\varPi\sn\ar 2$ fixed, we let $P$ be the
class of all pairs $\pi=(\varPi\ar{\sn 1}\sp,O\sp)$ such that all this data is
as we have assumed in Example 2.6\ssp. Then denoting by $\tilde f{_{\!}}_\pi$
the corresponding map, we define $\cal C=\{\,\tilde f{_{\!}}_\pi:\pi\in P\,\}
\ssp$. Applying Theorem 5\ssp, it now suffices to show that $
\tilde f{_{\!}}_\pi$ is continuous with open domain and almost Gateaux
differentiable, and that we can write $\delta\sp\tilde f{_{\!}}_\pi=
\tilde f{_{\!}}_{\pi_1\!}\Circ\sp\tilde\ell$ for some $\pi\ar 1\in P$ and $
\ssp\tilde\ell\in\cal L\ssp$.

Continuity, openness, and directional differentability follow by a simple
deduction from Lemma 2.7\ssp, as does $\delta\sp f\sp(\sp x\ssp,y\,;u\ssp,
v\sp)=u\circ[\,\roman{id\,},y\,]+\partial\ar 2\ssp x\circ[\,\roman{id\,},y\,]
\,.\,v\ssp$. De- fining $p\,(\sp\eta\,;\sp\xi\ar 1\sp,\xi\ar 2)=(\sp\eta\ssp,
\xi\ar 1)$ and $d\ar 2\ssp x\ssp(\sp\eta\,;\sp\xi\ar 1\sp,\xi\ar 2)=
\partial\ar 2\ssp x\ssp(\sp\eta\ssp,\xi\ar 1)\,\xi\ar 2\ssp,$ and writing $
\varPi\ar 3=$\Newline$\varPi\sqcap(\varPi\ar{\sn 1}\sn\sqcap\varPi\ar{\sn 1})\ssp,
$ and with $O\yr 1=\{\ssp(\sp\eta\,;\xi\ar 1\sp,\xi\ar 2):(\sp\eta\ssp,
\xi\ar 1)\in O$ and $\xi\ar 2\in\vecs\varPi\ar{\sn 1}\ssp\}$ further putting $
\pi\ar 1=(\varPi\ar{\sn 1}\sn\sqcap\varPi\ar{\sn 1}\ssp,O\yr 1)\ssp,$ we have a
continuous linear map\vskip.5mm\centerline{$
\tilde\ell=((E\sqcap F\sn\ar 1)\sqcap(E\sqcap F\sn\ar 1)\ssp,\Cinfty(\sp
O\yr 1_{\varPi_3}\sp,\varPi\ar 2)\sqcap C^{\ssp i\sp}(\sp Q_\varPi\sp,
\varPi\ar{\sn 1}\sn\sqcap\varPi\ar{\sn 1})\ssp,\sp\ell\,)$}

\noin defined by $(\sp x\ssp,y\,;u\ssp,v\sp)\mapsto(\sp u\circ p+d\ar 2\ssp
x\ssp,[\,y\ssp,v\,]\sp)\ssp$. A straightforward deduction shows\Newline that
we now have $\delta\sp\tilde f{_{\!}}_\pi=\tilde f{_{\!}}_{\pi_1\!}\Circ\sp
\tilde\ell\,$.                                                        \newQED

\noin{\bf
7} \Remarkss From Theorem 6\ssp, one almost directly obtains several particular
  results as {\it corollaries\ssp} by simple observations. For example, we
have the following.

(a) \ With $i\in\No$ and $\varPi=\varPi\ar{\sn 1}=\varPi\sn\ar 2=\biit R\,,$ let
$Q=I\inc\Re$ be a compact interval, and put $O=I\sn\times\sn\Re\ssp$. Fixing a
smooth $\varphi:\Re\to\Re\ssp,$ consider $g:G=\Cinfty(\Re\sp)\sqcap F\to$\Newline
$F=C^{\ssp i\sp}(I\sp)$ given by $(\sp x\ssp,y\sp)\mapsto x-\varphi\circ y\ssp
,$ where we understand $\dom g=\vecs G$ and $\dom(\sp u-v\sp)=(\dom u\sp)\cap
(\dom v\sp)\ssp$. The map $p:G\to E\sqcap F$ given by $(\sp x\ssp,y\sp)\mapsto
(\sp\bar x\ssp,y\sp)\ssp,$ where $\bar x\ssp(\sp\eta\ssp,\xi\sp)=x\ssp(\eta)-
\varphi\ssp(\xi)\ssp,$ is a sum of a continuous linear and constant map, hence
smooth. Since we have $g=f\circ p\,,$ from Theorem 6 in conjunction with the
chain rule [\,Hi1\ssp; Prop.\ 0.11, p.\ 240\,]\ssp, we obtain $(\sp G\sp,F\sp,
\sp g\sp)\in\CinftyPi\ssp$.

(b) \ Now that we have Theorem 6 available, to retrieve the variation formula\par\noin
(v)\hfill$\delta\sp f\sp(\sp x\ssp,y\,;u\ssp,v\sp)=u\circ[\,\roman{id\,},y\,]
  +\partial\ar 2\ssp x\circ[\,\roman{id\,},y\,]\,.\,v\ssp,$\hfill\phantom{(v)}

\noin we do not any longer have to make a recourse to Lemma 2.7 since we can
argue as follows. Fixing $x\ssp,y\ssp,u\ssp,v\ssp,$ for $t\not=0$ such that $
(\sp x+t\,u\ssp,y+t\,v\sp)\in\dom f\sp,$ define $\varDelta\ssp(t)=t^{-1}(\sp
f\ssp(\sp x+t\,u\ssp,y+t\,v\sp)-f\ssp(\sp x\ssp,y\sp))\ssp$. \hfil From
Theorem 6\ssp, we know that the li- mit $\lim_{\,t\to 0}\varDelta\ssp(t)=
\delta\sp f\sp(\sp x\ssp,y\,;u\ssp,v\sp)$ exists as a vector in the space $F_2
\ssp$. Since a con- tinuous linear map ev$_\eta:F_2\to\varPi_2$ is defined by
$y\mapsto y\ssp(\eta)\ssp,$ we can compute the li- mit "pointwise\ssp" as an
exercise in finite dimensions\par$\mhyppy{5.2}
 \delta\sp f\sp(\sp x\ssp,y\,;u\ssp,v\sp)\ssp(\eta)
=\roman{ev}_\eta\sp(\sp\delta\sp f\sp(\sp x\ssp,y\,;u\ssp,v\sp))
=\roman{ev}_\eta\sp(\sp\lim_{\,t\to 0}\varDelta\ssp(t))$\par$\mhyppy{31}
=\lim_{\,t\to 0}\roman{ev}_\eta\sp(\sp\varDelta\ssp(t))
=\lim_{\,t\to 0}\varDelta\ssp(t)\ssp(\eta)$\par$\mhyppy{31}
=\lim_{\,t\to 0}t^{-1}((\sp x+t\,u\sp)\ssp(\sp\eta\ssp,(\sp y+
 t\,v\sp)\ssp(\eta))-x\ssp(\sp\eta\ssp,y\ssp(\eta)))$\par$\mhyppy{31}
=u\ssp(\sp\eta\ssp,y\ssp(\eta))+
 \partial\ar 2\ssp x\ssp(\sp\eta\ssp,y\ssp(\eta))\ssp(\sp v\ssp(\eta))$\par$\mhyppy{31}
=(\sp u\circ[\,\roman{id\,},y\,]+\partial\ar 2\ssp x\circ
 [\,\roman{id\,},y\,]\,.\,v\sp)\ssp(\eta)\ssp$.

\noin Since this holds for every $\eta\in Q\ssp,$ we get (v)\ssp. In
particular for partial derivatives, we get $\partial\ar 2\sp f\sp(\sp x\ssp,
y\sp)\ssp v=\delta\sp f\sp(\sp x\ssp,y\,;\bnull E\ssp,v\sp)=\partial\ar 2\ssp
x\circ[\,\roman{id\,},y\,]\,.\,v\ssp$. Similarly, using induction on $k\in\N
\sp,$ one derives the higher order variation formula\vskip.5mm\centerline{$
\delta^{\,k\sn}f\sp(\biit z)=\partial_2^{\,k}x\circ[\,\roman{id\,},y\,]\,.\,
\bar v+\sum_{j=1}^{\,k}\partial_2^{\,k\sp-1}u_j \circ[\,\roman{id\,},y\,]
\,.\,\bar v_j\ssp,$}

\noin where $\biit z=\seq{\,(\sp x\ssp,y\sp)\ssp,(\sp u\ar 1\sp,v\ar 1)\ssp,
\ldots\,(\sp u\ai k\sp,v\ai k)\,}\ssp,$ and we have written $\bar v=
[\,v\ar 1\sp,\ldots\,v\ai k\ssp]$ and\Newline $\bar v_j=[\,v_1\sp,\ldots\,
v_{j-1},v_{j+1},\ldots\,v_k\ssp]\ssp$.

(c) \ In the proofs of Theorems 5.2 and 5.4 below, we need some information of
    the map $\tilde f\sn\ar 1=(E\sqcap C^{\,i\ssp+\sp 1\sp}(I\sp)\ssp,
F\ssn\ar 2\ssp,f\sn\ar 1)\ssp,$ where with $\varPi=\varPi\sn\ar 2=\biit R$ and
$\varPi\ar{\sn 1}=\biit R\sp\sqcap\biit R\,,$ we have $Q=I=[\,0\,,1\,]\ssp,$ and
$E$ and $F\ssn\ar 2$ are as in Example 2.6\ssp, and\par\centerline{$
f\sn\ar 1=\{\ssp(\sp\varphi\ssp,y\ssp,z\sp):\varphi\in\vecs E$ and $y\in\vecs
C^{\,i\ssp+\sp 1\sp}(I\sp)$ and $z=\varphi\circ[\,\roman{id\,};y\ssp,y\ssp'\,]
\in\vecs F\ssn\ar 2\ssp\}\ssp$.}

We then have $\tilde f\sn\ar 1\in\CinftyPi\ssp,$ and the formula\vskip.5mm\centerline{$
\partial\ar 2\ssp f\sn\ar 1\sp(\sp\varphi\ssp,y\sp)\ssp v
=\partial\ar 2\sp\varphi\circ[\,\roman{id\,};y\ssp,y\ssp'\,]\cdot v
+\partial\ar 3\sp\varphi\circ[\,\roman{id\,};y\ssp,y\ssp'\,]\cdot v\ssp'$}

\noin holds for $(\sp\varphi\ssp,y\sp)\in\dom f\sn\ar 1$ and $v\in\vecs
C^{\,i\sp+\sp 1\sp}(I\sp)\ssp,$ \hfil when labelling anew the partial
derivatives so that if $\partial_{_{\sp\roman{II}}}\sp\varphi$ denotes the
original $\partial\ar 2\sp\varphi$ with values in $\cal L\ssp(\varPi\ar{\sn 1}\sp,
\biit R\ssp)\ssp,$ the functions $\partial\ar 2\sp\varphi$ and $
\partial\ar 3\sp\varphi$ are defined by $\xi=(\sp\eta\,;\xi\ar 1\sp,\xi\ar 2)
\mapsto\partial_{_{\sp\roman{II}}}\sp\varphi\ssp(\xi)\ssp(\sp 1\ssp,0\ssp)$
and $\xi\mapsto\partial_{_{\sp\roman{II}}}\sp\varphi\ssp(\xi)\ssp(\sp 0\,,1\sp
)\ssp,$ respectively. For the proof, we note that $f\sn\ar 1=f\circ(\sp
\roman{id}\risti2\ell\ssp)\ssp,$ where $f$ is as in Example 2.6\ssp, and the
continuous linear $\ell:C^{\,i\ssp+\sp 1\sp}(I\sp)\to F\sn\ar 1$ is given by
$y\mapsto[\,y\ssp,y\ssp'\,]\ssp$. Then we can calculate

$\mhyppy{10.2}\partial\ar 2\ssp f\sn\ar 1\sp(\sp\varphi\ssp,y\sp)\ssp v
=\delta\sp f\sn\ar 1\sp(\sp\varphi\ssp,y\,;\bnull E\ssp,v\sp)
=\delta\sp f\ssp(\sp\varphi\ssp,\ell\,y\,;\bnull E\ssp,\ell\,v\sp)$

$\mhyppy{29}
=\partial_{_{\sp\roman{II}}}\sp\varphi\circ[\,\roman{id\,},\sp\ell\,y\,]\,.\,(\ssp\ell\,v\sp)
=\partial_{_{\sp\roman{II}}}\sp\varphi\circ[\,\roman{id\,};y\ssp,y\ssp'\,]\,.\,[\,v\ssp,v\ssp'\,]
$

$\mhyppy{29}
=\partial\ar 2\sp\varphi\circ[\,\roman{id\,};y\ssp,y\ssp'\,]\cdot v
+\partial\ar 3\sp\varphi\circ[\,\roman{id\,};y\ssp,y\ssp'\,]\cdot v\ssp'$.

\newProCla 8 Theorem.

The inclusion $\cal H_{_T}\inc\CinftyPi(\biit C\ssp)$ is in force.

\Prooff To apply Theorem 5\ssp, it suffices to verify that every $\tilde f=
(E\ssp,F\sp,f\sp)\in\cal H_{_T}$ has $\delta\sp\tilde f\in\cal H_{_T}\ssp$.
This we show using the Cauchy\,--\,formula, for which we introduce the
following particular weak integral concept. For a (continuous) function $
\gamma:\{\,\zeta\in\Ce:|\ssp\zeta\sp|=1\,\}\to\vecs F\sp,$ we write $
\oint\gamma\ssp(\zeta)\,d\ssp\zeta=v$ if{}f $v\in\vecs F$ and in the sense of
Riemann\,--\,integration, we have $\int_{\,0}^{\,1}\ell\,(\sp\gamma\ssp(\sp
e^{\,2\sp\pi\sp i\sp s\sp}))\,e^{\,2\sp\pi\sp i\,s}\,ds=\ell\,v$ for all $\ell
\in\cal L\ssp(F\sp,\biit C\,)\ssp$.

Now fix $\tilde f\sp$. If $x\in\dom f$ and $u\in\vecs E$ are such that $D\ar 2
=\{\,\zeta\in\Ce:|\ssp\zeta\sp|<2\,\}\inc\{\,t:x+t\,u\in\dom f\,\}\ssp,$ for
any $\ell\in\cal L\ssp(F\sp,\biit C\,)\ssp,$ we have $t\mapsto\ell\,(\sp f\sp(
\sp x+t\,u\sp))$ holomor- phic $D\ar 2\to\Ce$ in the ordinary sense. Hence the
ordinary Cauchy formula gives\par
\noin
(c)\hfill$f\sp(\sp x+t\,u\sp)=\oint\ssp(\sp\zeta-t\sp)^{-1\sn}
   f\sp(\sp x+\zeta\ssp u\sp)\,d\ssp\zeta\,$ for $\,|\ssp t\ssp|<1\,,$\hfill\phantom{(c)}

\noin whence again composing with $\ell\,,$ and then differentiating with
respect to $t\ssp,$ where we recall that $f$ is continuous\ssp, we obtain\par
\noin
(d) $\mhyppy{20.6}\delta\sp f\sp(\sp x\ssp,u\sp)
  =\oint\sp\zeta^{-2\sn}f\sp(\sp x+\zeta\ssp u\sp)\,d\ssp\zeta\,$.

A Hahn\,--\,Banach argument shows that if $\rng\gamma\inc V\in\ymp F$ with $V$
closed and absolutely convex, then also $\oint\gamma\ssp(\zeta)\,d\ssp\zeta\in
V_{}$. Using this, continuity of $f$ and (d) show that every $x\in\dom f$ has
some $U\in\ymp E$ with $\delta\sp f\,|\,((\sp x+U\sp)\sn\times\sn U\sp)$
continuous.\Newline Then using $\delta\sp f\sp(\sp x\ssp,u\sp)=t^{-1}\delta\sp
f\sp(\sp x\ssp,t\,u\sp)$ with $t\not=0\,,$ we obtain $\delta\tilde f\in\CPi{\ssp 0}\ssp$.

To prove that $\delta\sp\tilde f$ is almost Gateaux\,--\,differentiable,
fixing $x\ssp,u\ssp,u\ar 1\ssp,u\ar 2$ and $\ell\,,$ it suffices to show that
we have an ordinary holomorphic function $\varphi$ given by\vskip.5mm\centerline{$
t\mapsto\ell\,(\sp\delta\sp f\sp(\sp x+t\,u\ar 1\sp,u+t\,u\ar 2))
=\oint\sp\zeta^{-2}\sp\ell\,(\sp f\sp(\sp x+t\,u\ar 1\sn+
 \zeta\ssp(\sp u+t\,u\ar 2)))\,d\ssp\zeta$}\par$\mhyppy{22}
=\oint\sn\oint\sp\zeta^{-2}\sp(\sp\zeta\sp\ar 1\sn-t\sp)^{-1}
 \sp\ell\,(\sp f\sp(\sp x+\zeta\sp\ar 1\ssp u\ar 1\sn+\zeta\ssp(\sp u+\zeta\sp
 \ar 1\ssp u\ar 2)))\,d\ssp\zeta\sp\ar 1\ssp d\ssp\zeta\,,$

\noin where we obtained the last expression by using (c) with $u$ replaced by
$u\ar 1\sn+\zeta\ssp u\ar 2\ssp$. Since we assume $f$ continuous,
differentiability of $\varphi$ is a standard fact of finite dimensional
analysis, as in the deduction of (d)\ssp. We are done.                \newQED

\noin{\bf
9} \Examplee \hfil With \hfil $\roman D\sp(r)=\{\,\eta\in\Ce:|\ssp\eta\ssp|<r
  \,\}$ \hfil and \hfil $\bar{\roman D}\sp(r)=\{\,\eta\in\Ce:|\ssp\eta\ssp|\le
r\,\}\ssp$, \hfil let\linebreak $\Omega=\roman D\sp(1)\sn\times\sn\Ce\ssp$. \hfil
Hence \hfil taking \hfil $\biit K=\biit C\ssp$, \hfil and \hfil fixing $0<r<1
\ssp$, \hfil then \hfil $E=H\sp(\Omega)$ \hfil is\linebreak Fr\'echet, and $F=
H\ai b\sp(\sp\bar{\roman D}\sp(r))$ is a Banach space. Noting that the
topologies of both spaces $E\ssp,F$ require only uniform boundedness (on
compacta) of the functions themselves, and no derivatives are involved, using
the principles of the proof of Lemma 2.7\ssp, the reader may easily show that
$(E\sqcap F\sp,F\sp,f\sp)\in\cal H_{_T}\inc\CinftyPi(\biit C\ssp)$ holds for the
function $f$ defined by $\vecs(E\sqcap F\sp)\owns(\sp x\ssp,y\sp)\mapsto x
\circ[\,\roman{id\,},y\,]\ssp$.


\subhead 4

                   Implicit and inverse function theorems

If we have a continuous function $g$ between subsets of some topological
spaces with $(\sp x\ssp,y\sp)\in g\inc W$ and $W$ open in the product topology,
then we can choose open sets $U$ and $V$ with $(\sp x\ssp,y\sp)\in U\times V
\inc W$ and $g\ssp[\,U\,]\inc V$. Hence also taking into ac-\Newline count
[\,Hi1\ssp; p.\ 244, Lemma 1.3\,(1)\,]\ssp, we see that an equivalent
formulation of the {\it implicit function\ssp} theorem in [\,Hi1\ssp; p.\ 235\,]
is the following

\newProCla 1 Theorem.

Let $\tilde f=(E\sqcap F\sp,F\sp,f\sp)\in\CPi k$ with $k\in\N\cup\{\infty\}$
and $F$ a Banach space. If $(\sp z\ssp,b\ssp)\in f$ with $\partial\ar 2\sp f\sp
(z)$ bijective $\vecs F\to\vecs F\sp,$ there is some $W$ satisfying $z\in W\in
\taurd(E\sqcap F\sp)$ and $(E\ssp,F\sp,f\sp\inve\sp[\sp\{\sp b\sp\}\sp]\cap W\sp
)\in\CPi k\ssp,$ and such that also $\partial\ar 2\sp f\sp(w)$ is bijective $
\vecs F\to\vecs F$ for all $w\in W\sp$.\rm

\noin{\bf
2} {\font\=cmssi10\D\sp e\sp f\sp i\sp n\sp i\sp t\sp i\sp o\sp n\ssp}.
  For a vector map $\tilde f=(E\sqcap F\sp,F\sp,f\sp)\ssp,$ we say that $
\cal F$ is a Bp$\ar 2\ssp${\it--\,extension\ssp} of $\tilde f$ by $\cal R$
if{}f we here have functions $\cal R=\{\,(\sp i\ssp,\sp j\ssp,\sp
{\rho_{\sp}}_{i\sp j}\sp):(\sp i\ssp,\sp j\sp)\in\Delta\ssp\}$ and $\cal F=
\seq{\,({F_{}}_i\ssp,{f_{}}_i\ssp,{\rho_{\sp}}_i\sp):i\in I\sn\ar 1\ssp}\ssp,$
where $I\sn\ar 1=\dom\Delta\,,$ such that $\seq{\,(F\sp,{F_{}}_i\ssp,
{\rho_{\sp}}_i\sp):i\in I\sn\ar 1\ssp}$ is an LCS projective limit of $\{\ssp
(\sp i\ssp,\sp j\,;F_j\ssp,{F_{}}_i\ssp,\sp{\rho_{\sp}}_{i\sp j}\sp):(\sp
i\ssp,\sp j\sp)\in\Delta\ssp\}\ssp,$ and for all $i\in I\sn\ar 1$ we have $
{F_{}}_i$ a Banachable space together with ${\rho_{\sp}}_i\circ f\inc{f_{}}_i
\circ(\sp\roman{id}\risti2{\rho_{\sp}}_i\sp)\ssp$.

Assuming we can obtain existence of an implicit function, its
differentiability may be guaranteed using the following

\newProCla 3 Theorem.

Let $\seq{\,({F_{}}_i\ssp,{f_{}}_i\ssp,{\rho_{\sp}}_i\sp):i\in I\sn\ar 1\ssp}$
be a $\ssp\roman{Bp}\ar 2\ssp$--\,extension of $(E\sqcap F\sp,F\sp,f\sp)$ by $
\{\ssp(\sp i\ssp,\sp j\ssp,\sp{\rho_{\sp}}_{i\sp j}\sp):(\sp i\ssp,\sp j\sp)
\in\Delta\ssp\}\ssp,$ and let $b\in\vecs F\sp$. If $\sp\{\ssp(E\sqcap{F_{}}_i
\ssp,{F_{}}_i\ssp,{f_{}}_i\sp):i\in I\sn\ar 1\ssp\}\inc\CPi k$ with\Newline $k
\in\N\cup\{\infty\}\ssp,$ and if {\ssp\rm(1)} and {\ssp\rm(2)} below hold,
then $(E\ssp,F\sp,\sp g\sp)\in\CPi k\,$.\vskip.5mm\noin{\rm%
(1) \ }$g=f\sp\inve[\ssp\{\ssp b\ssp\}\ssp]$ and $\,\all{i\in I\sn\ar 1\sn}\,$\newline$\mhyppy{39.3}
   {g_{}}_i={f_{\sn}}_i\KN{.5}\inve[\ssp\{\,{\rho_{}}_i\ssp b\ssp\}\ssp]$ is a
   function and $\ssp\dom{g_{}}_i\inc\dom g\,,$\vskip0mm\noin{\rm%
(2) \ }$\all{i\ssp,x\ssp,y}\,i\in I\sn\ar 1$ and $(\sp x\ssp,y\ssp,b\ssp)\in f
   \imply\partial\ar 2\sp{f_{\sn}}_i\sp(\sp x\ssp,{\rho_{}}_i\ssp y\sp)$ is
   bijective $\vecs{F_{\sn}}_i\to\vecs{F_{\sn}}_i\,$.

\Prooff Using Theorem 3.3\ssp, it suffices to construct some $\cal G$ and $
\cal R$ such that $\cal G$ is a projective $\CPi k\ssp$--\,extension of $
(E,F,g)$ via $\cal R\,$. Let us show that we may take $\cal G=\seq{\,(E\ssp,
{F_{\sn}}_i\ssp,\sp{g_{}}_i\,;\idv E\ssp,\sp{\rho_{\sp}}_i\sp):i\in I\sn\ar 1
\ssp}$ and $\cal R=\{\,(\sp i\ssp,\sp j\,;\sp\idv E\ssp,\sp
{\rho_{\sp}}_{i\sp j}\sp):(\sp i\ssp,\sp j\sp)\in\Delta\,\}\ssp$. To check the
conditions of Theorem 3.3 and Definition 3.2\ssp, we must show that $g$ is a
function with $g_i=\rho_i\circ g\ssp,$ and that $(E,F_i,g_i)\in\CPi k$ holds
for all $i\in I\sn\ar 1\ssp$.

As each ${g_{}}_i$ is a function by (1)\ssp, for all $i\in I_1$ and all $x,y,z
\ssp,$ we have the impli- cations $(x,y),(x,z)\in g\imply(x,y,b),(x,z,b)\in f
\imply(x,\rho_iy,\rho_ib),(x,\rho_iz,\rho_ib)\in f_i\imply$ $(x,\rho_iy),
(x,\rho_iz)\in{g_{}}_i\imply\rho_iy=\rho_iz\,$. This immediately gives $\rho_i
\circ g\inc{g_{}}_i\ssp,$ and since by the projective limit property, we have
$[\ \all i\,\rho_iy=\rho_iz\ ]\imply y=z\ssp,$ also $g$ is a function. By $
\dom{g_{}}_i\inc\dom g\ssp,$ we further obtain ${g_{}}_i=\rho_i\circ g\ssp$.
                                                
To prove \hfill $(E\ssp,F{_{\ssn}}_i\ssp,g{_{\sn}}_i)\in\CPi k\,,$ \hfill
arbitrarily \hfill picking \hfill $w=(\sp x\ssp,z\sp)\in g{_{\sn}}_i\ssp,$ \hfill
it \hfill suffices \hfill to\linebreak have \hfill some \hfill $h$ \hfill with
\hfill $w\in h\inc g{_{\sn}}_i$ \hfill and \hfill $(E\ssp,F{_{\ssn}}_i\ssp,
h\sp)\in\CPi k\,$. \hfill Using \hfill $g{_{\sn}}_i=\rho{_{}}_i\circ g\ssp,$ \hfill
we \hfill see\linebreak that \hfill there \hfill is \hfill some \hfill $y$ \hfill
with \hfill $(\sp x\ssp,y\ssp,b\sp)\in f$ \hfill and \hfill $(\sp y\ssp,z\sp)
\in\rho{_{}}_i\ssp,$ \hfill whence \hfill Theorem 1 in\linebreak conjunction
with assumption (2) gives us the required $h={g_{}}_i\cap W\sp$.      \newQED

To exemplify the different conditions of Theorem 3\ssp, we modify a classical
example, and consider
\NS
\Examplee We let $\Delta=\{\ssp(\sp i\ssp,\sp j\sp):i\ssp,\sp j\in\Zepp$ and $
i\le j\,\}$ and $\Eps=\seq{\,e^{\,s}:s\in\Re\,}\ssp$. Writing $I_i=[-i\ssp,i\,
]\ssp,$ we now put $E=F=\Cinfty(\Re\sp)$ and $F_i=C^{\ssp i\sp}(\sp I_i)\ssp,$
and define $f=\{\ssp(\sp x\ssp,y\ssp,$ $x-\Eps\circ y\sp):x\ssp,y\in\vecs E\,
\}$ and $f_i=\{\ssp(\sp x\ssp,y\ssp,x-\Eps\circ y\sp):x\in\vecs E$ and $y\in$
$\vecs F_i\ssp\}\ssp$. Then we take $b=\bnull F=\Re\sn\times\sn\{\sp 0\sp\}$
and $k=\infty\ssp,$ and $\rho_{ij}=\{\ssp(\sp y\ssp,y\,|\,I_i\sp):y\in$ $\vecs
F_j\ssp\}$ and $\rho_i=\{\ssp(\sp y\ssp,y\,|\,I_i\sp):y\in\vecs F\,\}\ssp$.

From Remarks 3.7\ssp, we know $(E\sqcap F_i,F_i,f_i)\in\CinftyPi\ssp,$ and $
\partial\ar 2\sp{f_{\sn}}_i\sp(\sp\Eps\circ y\ssp,y\,|\,{I_{}}_i\sp)\,v=\Eps
\circ y\cdot v$ when $y\in\vecs E$ and $v\in\vecs F_i\ssp$. We also have $\dom
{g_{}}_i=\{\,x\in\vecs E:x\ssp[\,{I_{}}_i\ssp]\inc\Rep\ssp\}\not\subseteq\{\,x
\in\vecs E:\rng x\inc\Rep\ssp\}=\dom g\ssp$. Hence other conditions of Theorem
3 are satisfied, except that in (1) the requirement $\dom{g_{}}_i\subseteq\dom
g$ fails for all $i\ssp$.

\NSN{\bf
4} {\font\=cmssi10\D\sp e\sp f\sp i\sp n\sp i\sp t\sp i\sp o\sp n\ssp}.
  We say that $\cal F$ is a {\it B\ssp$\partial\ar 2\ssp$--\,extension\ssp} of
$\tilde f$ in $\CPi k$ by $\cal R$ around $z\yr 0$ if{}f it is a Bp$\ar 2\ssp
$--\,extension and with the notations of Definition 2 we have some $i\ar 0\in
I\ar 1$ and $W\ssn\ar 0\in\taurd(E\sqcap F_{i_0})$ and $z\yr 0=(\sp x\yr 0\!,
y\yr 0\!,b\ssp)\in f$ and $(\sp x\yr 0\!,\sp{\rho_{\sp}}_{i_{\sp 0\sp}}y\yr 0
\sp)\in W\ssn\ar 0\inc\dom{f_{\sn}}_{i_{\sp 0}}$ and\Newline $\{\ssp(E\sqcap
{F_{}}_i\ssp,{F_{}}_i\ssp,{f_{}}_i\sp):i\in I\sn\ar 1\ssp\}\inc\CPi k$ with $k
\in\N\cup\{\infty\}\ssp,$ and for $(\sp i\ar 0\ssp,\sp j\sp)\in\Delta$ we have
${\rho_{\sp}}_{i_0\sp j}$ injective and $\rho_j\circ f=f_j\circ(\sp\roman{id}
\risti 2\rho_j)$ and ${\rho_{\sp}}_{i_0\sp j}\circ f_j={f_{\sn}}_{i_{\sp 0\!}}
\circ(\sp\roman{id}\risti 2{\rho_{\sp}}_{i_0\sp j\sp})\ssp,$ and\vskip.5mm\noin
(1) \ $\partial\ar 2\sp f_{i_0}(x^0,\rho_{i_0}y^0)$ is bijective $\vecs F_{i_0}
    \to\vecs F_{i_0}\,,$\par\noin
(2) \ $\all{x\ssp,y\ssp,\sp\ell}\,(\sp x\ssp,y\ssp,\sp\ell\,)\in\partial\ar 2\sp
   f_j$ and $z=(\sp x\ssp,\sp{\rho_{\sp}}_{i_0\sp j\,}y\sp)\in
   W\ssn\ar 0$ and\newline\hyppy{7mm}$\partial\ar 2\sp f_{i_0}(z)$ is
   bijective $\vecs F_{i_0}\to\vecs F_{i_0}\,\imply\,\sp\ell\ssp$ is bijective
   $\vecs F_j\to\vecs F_j\,,$\par\noin
(3) \ $\all{x\yr 1\!,y\yr 1\!,\biit y}\,(\sp x\yr 1\!,y\yr 1)\in\vecs(E\sqcap
   F_{i_0})$ and $\biit y=\seq{\,{y_{}}_n\ssp}\in(\vecs F_j)^{\,I\!\!N}$ and\newline\hyppy{7mm}
   [ $\all{n\in\N}\,(\sp n^{-1}x\yr 0\sn+(\sp 1-n^{-1})\,x\yr 1\!,\sp
   {\rho_{\sp}}_{i_0\sp j\,}{y_{}}_n\ssp,{\rho_{\sp}}_{i_0\sp}b\ssp)\ssp,
   (\sp x\yr 1\!,y\yr 1\!,{\rho_{\sp}}_{i_0\sp}b\ssp)\in f_{i_0}$ ]\,,\newline\hyppy{7mm}
   and ${\rho_{\sp}}_{i_0\sp j}\circ\biit y\to y\yr 1$ in the topology $\taurd
   F_{i_0}\,\imply\,y\yr 1\in\rng{\rho_{\sp}}_{i_0\sp j}\,$.

\noin We say that $\tilde f$ is {\it B\ssp$\partial\ar 2\ssp$--\,extensible\ssp} 
in $\CPi k$ around $z\yr 0$ if{}f it is has one such extension.
\NS
The main part of Theorem 1 is generalized in the following {\it implicit function}

\newProCla 5 Theorem.

If $\tilde f=(E\sqcap F\sp,F\sp,f\sp)$ is B\ssp$\partial\ar 2\ssp$--\,
extensible in $\CPi k$ around $(\sp z\ssp,b\sp)\ssp,$ there is $W$ with $z\in
W\in\taurd(E\sqcap F\sp)$ and $(E\ssp,F\sp,f\sp\inve\sp[\sp\{\sp b\sp\}\sp]
\cap W\sp)\in\CPi k\ssp$.

\Prooff With the notations of Definition 4\ssp, let $\cal F$ be a {\it B\ssp}$
\partial\ar 2\ssp$--\,extension of $\tilde f$ in $\CPi k$ by $\cal R$ around $
z\yr 0=(\sp z\ssp,b\sp)=(\sp x\yr 0\!,y\yr 0\!,b\ssp)\ssp$. Applying Theorem 1
to $(E\sqcap F_{i_0}\sp,F_{i_0}\sp,{f_{\sn}}_{i_{\sp 0}})$ at $z\ar 2=(\sp
x\yr 0\!,\sp{\rho_{\sp}}_{i_{\sp 0\sp}}y\yr 0\sp)\ssp,$ we obtain some $
W\ssn\ar 2\inc W\ssn\ar 0$ with $z\ar 2\in W\ssn\ar 2\in\taurd(E\sqcap F_{i_0}
)$ and $(E\ssp,F_{i_0}\sp,\sp g\ar 0)\in\CPi k\ssp,$ when $g\ar 0=
f_{i_0}^{-\iota}[\ssp\{\,\rho_{i_0}b\ssp\}\ssp]\cap W\ssn\ar 2\ssp$. In
addition, we have $\partial_2f_{i_0}(w)$ bijective $\vecs F_{i_0}\to\vecs
F_{i_0}$ for $w\in W\ssn\ar 2\ssp$. Then we take any $U\in\taurd E$ with $
x\yr 0\in U\inc\dom g\ar 0$ and $U-x\yr 0$ circled, and further put $W=(\sp
\roman{id}\risti 2{\rho_{\sp}}_{i_0})\inve\sp[\,W\ssn\ar 2\ssp]\,|\,U\sp$. We
now have $z\in W\in\taurd(E\sqcap F\sp)\ssp,$ and writing $g=f\sp\inve\sp[\sp
\{\sp b\sp\}\sp]\cap W\sp,$ it remains to prove that also $\tilde g=(E\ssp,
F\sp,\sp g\sp)\in\CPi k$ holds.

Writing \hfil $I_2=\Delta\,[\sp\{\ssp i\ar 0\sp\}\sp]$ \hfil and \hfil $
\Delta\sn\yr 2=\Delta\cap I_2^{\times2}\ssp,$ \hfil and \hfil further \hfil
putting \hfil $W_j=\mathbreak(\sp\roman{id}\risti 2\rho_{i_0\sp j})\inve\sp[\,
W\ssn\ar 2\ssp]\,|\,U_{}$, \hfil we apply \hfil Theorem 3 \hfil with \hfil $
\Delta$ \hfil replaced by \hfil $\Delta\yr 2$ \hfil and \hfil $h_i=\mathbreak
f_i\ssp|\,W_i$ in place of $f_i$ and $h=f\,|\,W$ in place of $f\sp$. To obtain
$\tilde g\in\CPi k\ssp,$ we must verify conditions (1) and (2)\ssp, and show $
\rho_i\circ h\inc h_i\circ(\sp\roman{id}\risti2\rho_i)$ for all $i\in
I\ar 2\ssp$. For the last one, letting $(x,y,v)\in h\ssp,$ we must show $
(x,y,\rho_iv)\in h_i\circ(\sp\roman{id}\risti2\rho_i)\ssp$. By $(x,y)\in W_{}
$, we have $(x,\rho_{i_0i}(\rho_iy))=(x,\rho_{i_0}y)\in W\ssn\ar 2$ and $x\in
U_{}$, hence $(x,\rho_iy)\in W_i\ssp$. Now the conclusion follows from $
\rho_i\circ f\inc f_i\circ(\sp\roman{id}\risti2\rho_i)\ssp$.

Condition $g=h\inve[\sp\{\sp b\sp\}\sp]$ in (1) means $f\sp\inve\sp[\sp\{\sp b
\sp\}\sp]\cap W=(\sp f\,|\,W\sp)\inve\sp[\sp\{\sp b\sp\}\sp]\ssp,$ which is
trivial. For $i\in I\ar 2$ fixed, to show that $g_i=h_i^{-\iota}[\ssp\{\,
\rho_ib\ssp\}\ssp]$ is a function, letting $(\sp x\ssp,y\ai{\char'027})$ $\in g_i$ for $\nu=1
\ssp,2\ssp,$ we have $(x,y_\nu,\rho_ib)\in f_i$ with $(x,y_\nu)\in W_i\ssp,$
hence $(x,\rho_{i_0i}y_\nu)\in$ $W\ssn\ar 2\ssp,$ and by $\rho_{i_0i}\circ f_i
\inc f_{i_0}\circ(\sp\roman{id}\risti 2\rho_{i_0i})$ also $(x,\rho_{i_0i}y_\nu
,\rho_{i_0}b)\in f_{i_0}$ holds. Consequently\Newline $(x,\rho_{i_0i}y_\nu)\in
g\ar 0\ssp,$ hence $y\ar 1=y\ar 2\ssp,$ as $g\ar 0$ is a function and $
\rho_{i_0i}$ is injective.

To prove \hbox{$\dom g{_{}}_i\inc\dom g$} in (1)\ssp, arbitrarily fixing \hbox{$
j\in I\ar 2\ssp$}, we first show $\dom g_j$ $\in\taurd E\ssp$. Thus letting \hbox{$
(\sp x\ssp,y\sp)\in g_j\ssp$}, we have \hbox{$(\sp x\ssp,y\ssp,\rho_j\ssp b\sp)
\in h_j\inc f_j$} with \hbox{$(\sp x\ssp,\rho_{i_0j}\sp y\sp)\in W\ssn\ar 2
\ssp$}. By (2) of Definition 4\ssp, we then have \hbox{$
\partial\ar 2\sp h_j\sp(\sp x\ssp,y\sp)=\partial\ar 2\ssp f_j\sp(\sp x\ssp,
y\sp)$} bijective \hbox{$\vecs F_j\to\vecs F_j\ssp$}, whence Theorem 1 gives $
W\!\ar 1$ with \hbox{$(\sp x\ssp,y\sp)\in W\!\ar 1$} and \hbox{$(E\ssp,F_j\sp,
g_j\cap W\!\ar 1)\in\CPi k\ssp$}. In particular, we have \hbox{$x\in U\sn\ar 1
=\dom(g_j\cap W\!\ar 1)\in\taurd E$} and \hbox{$U\sn\ar 1\inc\dom g_j\ssp$}.
We have obtained $\dom g{_{\sn}}_j\in\taurd E$ since $x$ is arbitrary.

To proceed in the proof of $\dom g_i\inc\dom g\ssp,$ we next show $\dom g_j=U
\sp$. Trivially having $\dom g_j\inc U\sp,$ we fix $x\in U$ and show $x\in\dom
g_j\ssp$. For this, writing $S=\mathbreak\{\,s:0\le s\le 1$ and $\all r\,0\le
r\le s\imply(\sp 1-r\sp)\,x\yr 0\sn+r\,x\in\dom g_j\,\}$ and $t=\sup\sp S\ssp,
$ it suffices to show $1=t\in S\ssp$. Note that by $(x^0,y^0,b)\in f$ and $
\rho_j\circ f\inc f_j\circ(\sp\roman{id}\risti2\rho_j)\ssp,$ we have $0\in S
\ssp,$ hence $S\not=\emptyset$ and $t\ge 0\,$. If we have $t\in S\ssp,$ then
by $\dom g_j\in\taurd E\ssp,$ we indirectly see $t=1\ssp,$ whence it remains
to prove $t\in S\ssp$.

To \hfil get \hfil $t\in S\ssp,$ \hfil we \hfil use \hfil (3) of Definition 4 \hfil
with \hfil $x^1=(\sp 1-t\sp)\,x\yr 0\sn+t\,x$ \hfil and \hfil $y^1=\mathbreak
g\ar 0(x^1)$ and $y_n=g_j(x_n)\ssp,$ where $x_n=n^{-1}x^0+(1-n^{-1})\,x^1\in
\dom g_j\ssp$. By $\rho_{i_0j}\circ f_j$ $\inc f_{i_0}\circ(\sp
\roman{id}\risti2\rho_{i_0j})\ssp,$ we then have $(x_n,\rho_{i_0j}y_n,
\rho_{i_0}b)\in f_{i_0}\ssp,$ and observing $(x_n,\rho_{i_0j}y_n)$\Newline$\in
W_2\ssp,$ hence also $\rho_{i_0j}y_n=g\ar 0(x_n)\ssp$. Since $(E,F_{i_0},
g\ar 0)\in\CPi k\ssp,$ we have $g\ar 0$ continuous $\taurd E\to\taurd F_{i_0}
\ssp,$ which gives $\rho_{i_0j}\circ\biit y\to y^1$. Now (3) gives existence
of some $y$ with $(y,y^1)\in\rho_{i_0j}\ssp$. Then using $\rho_{i_0j}\circ f_j
\iinc f_{i_0}\circ(\sp\roman{id}\risti2\rho_{i_0j})$ and injectivity of $
\rho_{i_0j}\ssp,$ we see\Newline $(x^1,y)\in g_j\ssp,$ hence $x^1\in\dom g_j$
and $t\in S\ssp$.

Now to prove $\dom g_i\inc\dom g\ssp$, fixing $x\in\dom g_i=U_{}$, we write $
y_j=g_j(x)$ for $j\in I_2\ssp$. Using $\rho_{i_0j}\circ f_j\inc f_{i_0}\circ(\sp
\roman{id}\risti2\rho_{i_0j})\ssp$, we first obtain $y_{i_0}=\rho_{i_0j}y_j$
for all $j\in$ $I_2\ssp$.\Newline Then utilizing $\rho_{i_0l}=\rho_{i_0j}\circ
\rho_{jl}$ and injectivity of $\rho_{i_0j}\ssp$, we get $y_j=\rho_{jl}y_l$ for
$(j,l)$ $\in\Delta\yr 2\ssp$. Hence there is $y$ with $\{\ssp y\ssp\}=\bigcap\ssp
\{\,\rho_j^{-\iota}\sp[\ssp\{\ssp y_j\sp\}\ssp]:j\in I_2\ssp\}\ssp$, and to
prove $(x,y)\in$ $ g\ssp$, we must show $(x,y,b)\in f$ and $(x,\rho_{i_0}y)\in
W_2\ssp$. Using $\rho_j\circ f\iinc f_j\circ(\sp\roman{id}\risti2\rho_j)\ssp$,
we get $(x,y,b)\in f\sp$, and further note $(x,\rho_{i_0}y)=(x,y_{i_0})\in
g_{i_0}\inc W_2\ssp$.

For (2) assuming $(x,y,b)\in h\ssp,$ we must show $\partial_2
h_i(x,\rho_iy)$ bijective $\vecs F_i\to\vecs F_i\ssp$. Using (2) of Definition
4\ssp, it suffices to show $w=(x,\rho_{i_0i}(\rho_iy))\in W\ssn\ar 2\ssp$. This
indeed holds, since $(x,y)\in W\imply w=(x,\rho_{i_0}y)\in W\ssn\ar 2\ssp$.
Now all conditions of Theorem 3 are verified,
whence we obtain $\tilde g\in\CPi k\ssp$.                             \newQED

Let us say that a vector map $\tilde f=(E\ssp,E\ssp,f\sp)$ is {\it BI\,--\,
extensible\ssp} in $\CPi k$ about $z$ if{}f $(E\sqcap E\ssp,E\ssp,f\sn\ar 1)$
is {\it B\ssp}$\partial\ar 2\ssp$--\,extensible in $\CPi k$ around $(\sp
z\ssp,\bnull E)\ssp,$ when $f\sn\ar 1=\{\,(\sp x\ssp,v\ssp,u-x\ssp):{}$\Newline
$(\sp v\ssp,u\sp)\in f$ and $x\in\vecs E\,\}\ssp$. We say that $\tilde f$ is a
local $\CPi k\ssp$--\,diffeomorphism $U\to V$ if{}f $g=f\,|\,U$ is bijective $
U\to V$ with $(E\ssp,E\ssp,\sp g\sp)\ssp,(E\ssp,E\ssp,\sp g\inve\sp)\in
\CPi k\ssp$. From Theorem 5\ssp, we now obtain an {\it inverse function
theorem} as

\newProCla 6 Corollary.
                                        
If $\tilde f=(E\ssp,E\ssp,f\sp)$ is BI\,--\,extensible in $\CPi k$ about $z\ssp
,$ \hfill then some $U,V$ satisfy $z\in V\times U$ with $\tilde f$ a local $
\CPi k\ssp$--\,diffeomorphism $U\to V_{}$.

\Prooff Let $\tilde f$ be {\it BI\,--}\,extensible in $\CPi k$ about $z=
(x_0,y_0)\ssp$. Noticing that we have $f_1^{-\iota}[\sp\{\sp\bnull E\}\sp]=
\{\ssp(\sp x\ssp,y\sp):(\sp y\ssp,x\sp)\in f\,\}=f\sp\inve,$ from Theorem 5\ssp,
we get an open $W$ with $z\in W$ and $(E,E,g_2)\in\CPi k\ssp,$ where we put $
g_2=f\sp\inve\sn\cap W\sp$. By continuity of $g_2\ssp,$ we can choose open $
U_1,V_1$ with $z\in V_1\times U_1\inc W$ and $g\ar 2\KN1\image V\ssn\ar 1\inc
U\sn\ar 1\ssp$.

Letting $\seq{\,(E_i\ssp,f_{1i}\ssp,\rho_i):i\in I\sn\ar 1\ssp}$ be a {\it
B\ssp}$\partial\ar 2\ssp$--\,extension of $(E\sqcap E\ssp,E\ssp,f\sn\ar 1)$ in
$\CPi k$ by $\{\,(\sp i\ssp,\sp j\ssp,\sp{\rho_{\sp}}_{i\sp j}\sp):(\sp i\ssp,
\sp j\sp)\in\Delta\ssp\}$ around $(\sp z\ssp,\bnull E)\ssp,$ we define $f_i=
f_{1i}(\sp\bnull E\ssp,\,\cdot\,)\ssp$. We now have $\rho_i\circ f=f_i\circ
\rho_i$ and $\rho_{i_0j}\circ f_j=f_{i_0}\circ\rho_{i_0j}$ for $(\sp
i\ar 0\ssp,\sp j\sp)\in\Delta\ssp$. If we had $(*)\ \ \rho_{ij}\circ f_j=$ $
f_i\circ\rho_{ij}$ for $(i,j)\in\Delta\sn\yr 0=\Delta\ssp|\ssp I\ar 0\ssp$,
where $I\ar 0=\Delta\!\image\sn\{\ssp i\ar 0\}\ssp$, then $\seq{\,
({E_{}}_i\ssp,{E_{}}_i\ssp,{f_{\sn}}_i\,;\sp\rho^{}_{\sp i}\ssp,
\rho^{}_{\sp i}\sp):{}$ $i\in I\ar 0\ssp}$ would be a projective $\CPi k$
extension of $\tilde f$ via $\{\,(\sp i\ssp,\sp j\,;\sp\rho^{}_{\sp i\sp j}
\ssp,\rho^{}_{\sp i\sp j}\sp):(\sp i\ssp,\sp j\sp)\in\Delta\sn\yr 0\ssp\}\ssp,
$ and Theorem 3.3 would give $\tilde f\in\CPi k\ssp$. To prove this, it thus
suffices to obtain ($*$)\ssp. Recalling that $\rho_{i_0i}$ is injective, we
have\par$\mhyppy{10}
 \rho_{i_0i}\circ(\sp\rho_{ij}\circ f_j)
=(\sp\rho_{i_0i}\circ\rho_{ij})\circ f_j
=\rho_{i_0j}\circ f_j$\newline$\mhyppy{36}
=f_{i_0}\circ\rho_{i_0j}=f_{i_0}\circ(\sp\rho_{i_0i}\circ\rho_{ij})
=(\sp f_{i_0}\circ\rho_{i_0i})\circ\rho_{ij}$\newline$\mhyppy{36}
=(\sp\rho_{i_0i}\circ f_i)\circ\rho_{ij}
=\rho_{i_0i}\circ(\sp f_i\circ\rho_{ij})\imply(*)\ssp$.

Now having $\tilde f\in\CPi k\ssp,$ using continuity of $f\sp,$ we pick an
open $U$ with $y_0\in U\inc$ $(\dom f\sp)\cap U_1$ and $V=f\sp\image\sn U\inc
V_1\ssp$. Putting $g_1=g_2\ssp|\,g_2\inve[\,U\ssp]$ and $f\,|\,U=g\ssp,$ then
$z\in V\times U$ and $(E\ssp,E\ssp,\sp g\sp)\ssp,(E\ssp,E\ssp,\sp g\ar 1)\in
\CPi k\ssp,$ and it remains to show $g_1=g\inve$.

For $g_1\inc g\inve,$ we first observe that $\rng g_1\inc U\sp,$ since $g_2$
is a function. Hence for $(x,y)\in g_1\inc g_2\inc f\sp\inve,$ we have $y\in U
\sp,$ whence $(y,x)\in f\,|\,U=g\ssp,$ and $(x,y)\in g\inve$. Conversely, for
$g\inve\inc g\ar 1\ssp,$ we deduce

$(x,y)\in g\inve\imply(y,x)\in f\,|\,U\imply(x,y)\in f\inve$ and $y\in U\inc
 U_1$ and $x\in f\sp\image\sn U\inc V_1$\par$\mhyppy{13}
\imply y\in U$ and $(x,y)\in f\inve$ and $(x,y)\in V_1\times U_1\inc W$\par$\mhyppy{8}
\imply y\in U$ and $(x,y)\in f\sp\inve\sn\cap W=g_2\imply(\sp x\ssp,y\sp)\in
 g\ar 2\ssp|\,g\ar 2\KN1\inve\sp[\,U\ssp]=g\ar 1\ssp$.                   \QED
                                

\subhead5

                           Smooth solution maps

Our goal in this section is to give various examples of applications of
Theorems 4.1, 4.3, and 4.5 to prove results about smooth dependence of the
solution of a nonlinear (partial) differential equation on the data.
"Data\ssp" here refers not only to initial or boundary conditions, but also to
the "nonlinearity\ssp" of the equation. We only want to illustrate the
principles of application, and we do not aim at hard results. To be able to
give fairly complete proofs within moderate space, we have striven for
examples that, say, do not require too "heavy\ssp" apriori estimes.

\NSN{\bf
1} \Examplee We consider the equation $\varphi\ssp(\sp s\,;y\ssp(s)\ssp,
y\ssp'(s))=0$ for $s\in I=[\,0\,,1\,]\ssp,$ requiring the initial condition $
y\ssp(0)=\eta\in\Re\,$. Recalling our basic definitions, we write this
equation as $\varphi\circ[\,\roman{id\,};y\ssp,y\ssp'\,]=I\sn\times\sn\{\sp 0
\sp\}\ssp$. Here we assume $y\in\vecs C^{\ssp 1}(I\sp)$ and $\varphi:O\to\Re$
smooth for a fixed $O\ssp,$ which is open in the natural topology of $I\sn
\times\sn(\Re\sn\times\sn\Re\sp)\ssp,$ induced from $\taurd(\biit R\sp\sqcap(
\biit R\sp\sqcap\biit R\,))\ssp$.

We put $E=E\ar 0\sn\sqcap E\ar 1=\biit R\sqcap\Cinfty(O)\ssp,$ and denote by $
F_i$ the topological linear subspace of $G=C^{i+1}(I\sp)$ with $\vecs F_i=\{\,
y\in\vecs G:y\ssp(0)=0
                      \,\}\ssp$. The {\it solution relation} of our equation
we call the set $\varSigma$ of all $w=(\sp x\ssp,y\sp)=(\sp\eta\ssp,
\varphi\ssp,y\sp)\in\vecs(E\sqcap C^{\ssp 1}(I\sp))$ with\Newline $\varphi
\circ[\,\roman{id\,};y\ssp,y\ssp'\,]=I\sn\times\sn\{\sp 0\sp\}$ and $y\ssp(0)=
\eta\ssp$. We let $\varSigma\ar 0$ be the set of all {\it regular\ssp} $w\in
\varSigma\ssp,$ this meaning that if we write $p_\iota=\partial_\iota\varphi
\circ[\,\roman{id\,};y\ssp,y\ssp'\,]$ for $\iota=2\ssp,3\ssp,$ then $0\not\in
\rng p_3$ and every $v\in\vecs F_0$ has a unique $u\in\vecs F_0$ with $
p\sp\ar 3\sp u\ssp'+p\sp\ar 2\sp u=v\ssp'$.

An elementary exercise in ordinary differential equations shows that for $w$
to\linebreak be \hfill regular \hfill it \hfill suffices \hfill just \hfill to \hfill
require \hfill $0\not\in\rng p\sp\ar 3\ssp,$ \hfill since \hfill the \hfill $
C^{\ssp 1\,}$--\,solution \hfill $u$ \hfill to\linebreak $p\sp\ar 3\sp u\ssp'+
p\sp\ar 2\sp u=v\ssp'$ with $u\ssp(0)=0$ is $u=e^{-\ssmb A}\!\int_{\ssp 0}\sp
(\sp e^{\ssp\ssmb A\sp}p\sp\ar 3\KN1^{-1}\sp v\ssp')\ssp$, where $\smb A=
\int_{\ssp 0}\sp(\ssp p\sp\ar 3\KN1^{-1}\sp p\sp\ar 2)\ssp$.

As an application of Theorems 4.1 and 4.3, we now show that for every regular
$w\in\varSigma$ there are $U\in\taurd E$ and $V\in\taurd C^{\ssp 1}(I\sp)$
with $w\in U\sn\times\sn V\sp,$ and such that the relation $g=\varSigma\cap(\sp
U\ssn\times\sn V\sp)$ actually is a function $U\to\vecs\Cinfty(I\sp)$ with $(
E\ssp,\Cinfty(I\sp)\ssp,g\sp)$\Newline a smooth map. Taking into account the
three lines at the beginning of Section 4\ssp, we see that this is
equivalently encoded in the following succinct formulation of

\newProCla 2 Theorem.

For every $w\in\varSigma\ar 0$ there is some $W$ with $w\in W\in\taurd(E\sqcap
C^{\ssp 1}(I\sp))$ and $(E\ssp,\Cinfty(I\sp)\ssp,\varSigma\cap W\sp)\in\CinftyPi\ssp$.

\Prooff Utilizing the notation $\bar\eta=I\sn\times\sn\{\eta\}\ssp,$ for $i\in\No
\ssp,$ we define the maps $\tilde f_i$ $=(E\sqcap F_i,F_i,f_i)$ by $f_i(w)=
\int_{\,0}\varphi\circ[\,\roman{id\,};y+\bar\eta\ssp,y\ssp'\,]\ssp,$ thus
requiring that $f_i(w)(t)=\int_{\,0}^{\,t}\varphi\ssp(\sp s\,;\eta+y\ssp(s)
\ssp,y\ssp'(s))\,ds$ must hold for all $t\in I\sp$. Here we take $w=(\eta,
\varphi,y)\in\vecs(E\sqcap F_i)$ such that $\dom(\varphi\circ[\,\roman{id\,};
y+\bar\eta\ssp,y\ssp'\,]\sp)=I\sp,$ i.e., such that for every $s\in I\sp,$ we
have $(\sp s\,;\eta+y\ssp(s)\ssp,y\ssp'(s))\in O\ssp$. 

We now have the equivalence $(\sp\eta\ssp,\varphi\ssp,\bar\eta+y\sp)\in
\varSigma\equivv(\sp\eta\ssp,\varphi\ssp,y\sp)\in f_0^{\sp-\iota\sp}[\sp\{\sp
b\sp\}\sp]$ for all $\eta\ssp,\varphi\ssp,y\ssp,$ when $b=I\sn\times\sn\{\sp 0
\sp\}=\bnull F\sn\LHB{.3}{_{_i}}\ssp$. Since $\int_{\,0}$ defines a continuous
linear map\par\centerline{$
C^{\ssp i}(I\sp)\to{F_{\sn}}_i\ssp,$ and $(\eta,y)\mapsto[\,y+\bar\eta\ssp,
y\ssp'\,]$ respectively $\biit R\sqcap F_i\to C^{\ssp i}(\sp I\sp,\biit R\sp
\sqcap\biit R\sp)\ssp,$}\vskip-.5mm
\noin
the chain rule [\,Hi1\ssp; Prop.\ 0.11, p.\ 240\,] in conjunction with Theorem
3.6 gives us the result $\tilde f_i\in\CinftyPi$ for all $i\ssp$.

Fixing $(\sp\eta\ssp,\varphi\ssp,y\sp)\in\varSigma\ar 0\ssp,$ we prove $y\in
\vecs\Cinfty(I\sp)\ssp$. Since $0\not\in\rng p_3$ holds by regu- larity, using
Theorem 4.1\ssp, or even the classical finite dimensional implicit function
theorem in its utterly simple form, we see that every $t\in I\setminus\sn
\{\ssp 0\ssp,1\ssp\}$ has a smooth ${\psi_{}}_t$ defined on an open
neighborhood of $(\sp t\ssp,y\ssp\sp(t))$ such that (i) \ $y\ssp'(s)=
{\psi_{}}_t(\sp s\ssp,y\ssp(s))$ holds for $s\in I$ close to $t\ssp$.
Induction then shows that $y$ is smooth on ${]}\,0\,,1\ssp{[}\,$. To prove $y
\in\vecs\Cinfty(I\sp)\ssp,$ it suffices to obtain (i) also for $t=0\,,1\ssp$.

Using continuity of the partial derivatives of $\varphi\ssp,$ for each fixed $
k\in\N\sp,$ we can extend $\varphi$ to a $C^{\ssp k}$ function $\bar\varphi
\ssp,$ for example $\bar\varphi\ssp(\sp s\ssp,\xi\sp)=\sum_{\,i=0\,}^{\,k}(\sp
i\ssp!\sp)^{-1}s^{\ssp i}\sp\partial_1^{\,i}\varphi\ssp(\sp 0\,,\xi\sp)$ when
$s<0\,,$ with $(\sp 0\,;\eta\ssp,y\ssp'(0))\ssp,(\sp 1\,;y\ssp(1)\ssp,
y\ssp'(1))\in\dom\bar\varphi\in\taurd(\biit R\sqcap(\biit R\sp\sqcap\biit R\,)
)\ssp$. Then an application of the implicit function theorem gives (i) with $
{\psi_{}}_t$ now a $C^{\ssp k}$ function. As above, we deduce $y\in$ $\vecs
C^{\ssp k}(I\sp)\ssp,$ and noting that here $k\in\N$ is arbitrary, then $y\in
\vecs\Cinfty(I\sp)$ immediately follows.

Now preparing ourselves to the application of Theorem 4.3\ssp, we let $F$ be
the topological linear subspace of $G=\Cinfty(I\sp)$ with $\vecs F=\{\,y\in
\vecs G:y\ssp(0)=0\,\}\ssp,$ and further put $\Delta=\{\ssp(i,j):i,j\in\No$
and $i\le j\,\}\ssp,$ thus having $I_1=\No\ssp$. Also let $\rho_i=\idv F$ and
$\rho_{ij}=\idv F_j$ for $(i,j)\in\Delta\ssp$. Now fixing $w=(\sp\eta\ssp,
\varphi\ssp,y\sp)\in\varSigma\ar 0\ssp,$ we have $w\ar 0=(\sp\eta\ssp,\varphi
\ssp,y-\bar\eta\sp)\in\varSigma\ar 1\ssp,$ when $\varSigma\ar 1$ is the set of
all $w_1=(\eta_1,\varphi_1,y_1)$ such that $\bar w_1=(\eta_1,\varphi_1,
\bar\eta_1+y_1)\in\varSigma\ar 0\ssp$.

Using the result $\rng\varSigma\ar 0\inc\vecs\Cinfty(I\sp)$ from above, we
obtain $\varSigma\ar 1\inc\dom f_i$ for all $i\in\No\ssp$. Recalling Remarks
3.7\sp(c)\ssp, we see that if $w_1\in\varSigma_1\ssp,$ then $\partial_2f_i
(w_1):\vecs F_i\to$ $\vecs F_i$ \hfill is \hfill the \hfill function \hfill $u
\mapsto\int_{\,0}(\sp p\ar 3\sp u\ssp'+p\ar 2\sp u\sp)\ssp,$ \hfill which \hfill
for \hfill $i=0$ \hfill is \hfill bijective \hfill by \hfill the\linebreak
assumed regularity. Since here $p{_{}}_\iota={\partial_{}}_\iota\varphi\ar 1\sn
\circ[\,\roman{id\,};\bar\eta\ar 1\sn+y\ar 1\sp,y\ar 1\KN1'\,]$ are smooth, a
simple induction using $u'=p_3^{-1}(v'-p_2u)$ shows that bijectivity holds for
all $i\in\No\ssp$.

Applying Theorem 4.1 to $\tilde f_0$ at $(\sp w\ar 0\sp,b\ssp)\ssp,$ we get $
W_2$ with $w\ar 0\in W_2\in\taurd(E\sqcap F_0)$ and $(E,F_0,f_0^{-\iota}[\sp
\{\sp b\sp\}\sp]\cap W\ssn\ar 2)\in\CinftyPi$ and $\partial_2f_0(w_1)$
bijective $\vecs F_0\to\vecs F_0$ for $w_1\in$ $W_2\ssp$. We want $W_1$ with $
w_0\in W_1\in\taurd(E\sqcap F_0)\ssp,$ and such that $\bar w_1$ would be
regular\Newline for $w_1\in f_0^{\sp-\iota\sp}[\sp\{\sp b\sp\}\sp]\cap W_1\ssp
$. From the preceding, we see that the unique solvability holds for $w_1\in
W_2\ssp$. Using continuity of $(\sp s\,;\eta\ar 1\ssp,\varphi\ar 1\ssp,y\ar 1)
\mapsto\partial\ar 3\sp\varphi\ar 1\sp(\sp s\,;\eta\ar 1\sn+y\ar 1\sp(s)\ssp,
y\ar 1\KN1'(s))$\Newline in conjunction with compactness of $I\sp,$ we see
that suitably shrinking $W_2$ to $W_1\ssp,$ we have $\bar w_1$ regular for $
w_1\in f_0^{\sp-\iota\sp}[\sp\{\sp b\sp\}\sp]\cap W_1\ssp$. Thus $
f_0^{\sp-\iota\sp}[\sp\{\sp b\sp\}\sp]\cap W_1\inc\varSigma\ar 1\ssp$.

Now with \hbox{$f\sNorr{\snn\sixroman1i}=f\sNorr{\snn i}\ssp|\,W\ssn\ar 1$}
and \hbox{$f=\bigcap\ssp\{\,f\sNorr{\snn\sixroman1i}:i\in\No\ssp\}\ssp$}, we
put \hbox{$g\sNor i=f_{\sixroman1i}^{\sp-\iota\sp}[\sp\{\sp b\sp\}\sp]$} and $
g$\linebreak \hbox{$=f\sp\inve\sp[\sp\{\sp b\sp\}\sp]\ssp$}. Then taking \hbox{$
\tilde f=(E\sqcap F_{\sp},F_{\sp},f\sp)\ssp$}, a straighforward verification
shows that $\seq{\,(F\sNorr{\snn i}\ssp,f\sNorr{\snn\sixroman1i}\ssp,
\rho\sNor{i\sp}):i\in\No\ssp}$ is a Bp$\ar 2\ssp$--\,extension of $\tilde f$
by \hbox{$\{\ssp(\sp i\ssp,\sp j\sp,\rho\sNor{ij\sp}):(\sp i\ssp,\sp j\sp)\in
\Delta\,\}\ssp$}. As we have $W\!\ar 1\in\taurd(E\sqcap F\sNorr 0)\ssp$, also
$(E\sqcap F\sNorr{\snn i}\ssp,F\sNorr{\snn i}\ssp,f\sNorr{\snn\sixroman1i\sp})
\in\CinftyPi$ holds for $i\in\No\ssp$.

\def\oseoy{\overset{_{_{\kern1mm\roman o}}}\to y}
A continuous linear map $\ell:E\sqcap C^{\ssp 1}(I\sp)\to E\sqcap F_0$ is
defined by $(\sp\eta\ar 1\sp,\varphi\ar 1\sp,y\ar 1\sp)\mapsto(\sp
y\ar 1(0)\ssp,\varphi\ar 1\sp,\oseoy\sn\ar 1)\ssp,$ \hfil where we write \hfil $
\oseoy\sn\ar 1=y\ar 1\sn-I\sn\times\sn\{\sp y\ar 1(0)\sp\}\ssp$. \hfil Hence \hfil
$W=\,\ell\ssp\inve\sp[\,W\ssn\ar 1\ssp]\in\mathbreak\taurd(E\sqcap C^{\ssp 1}(
I\sp))\ssp$. Since $\eta=y\ssp(0)$ and $(\sp\eta\ssp,\varphi\ssp,y-\bar\eta\sp
)=w\ar 0\in W\ssn\ar 1\ssp,$ also $w\in W$ holds. We are done once we show $(
E\ssp,F\sp,g\sp)\in\CinftyPi$ and (x) \ that $\varSigma\cap W$ is the function
$\dom g\owns x\ar 1=(\sp\eta\ar 1\sp,\varphi\ar 1)\mapsto\bar\eta\ar 1\sn+
g\ssp(\sp x\ar 1)\ssp$.

To prove $(E\ssp,F\sp,g\sp)\in\CinftyPi\ssp,$ it suffices that we verify
conditions (1) and (2) of Theorem 4.3\ssp. Since by $f_0^{-\iota}[\sp\{\sp b
\sp\}\sp]\cap W_1\inc f_0^{-\iota}[\sp\{\sp b\sp\}\sp]\cap W_2\ssp,$ we have $
f_0^{-\iota}[\sp\{\sp b\sp\}\sp]\cap W_1$ a function, for (1) it suffices to
show $g_i\inc g_0\inc f_0^{-\iota}[\sp\{\sp b\sp\}\sp]\cap W_1$ and $g_0\inc g
$ for all $i\in\No\ssp$. We trivially have $g\inc g_i\inc g_0\ssp,$ and $g_0
\inc g$ is a simple consequence of $g_0\inc f_0^{\sp-\iota\sp}[\sp\{\sp b\sp\}
\sp]\cap W_1\inc\varSigma_1$ and $\rng\varSigma\ar 1\inc\vecs F\sp$. For $g_0
\inc f_0^{-\iota}[\sp\{\sp b\sp\}\sp]\cap W_1$ we only observe $g_0=
f_{10}^{-\iota}[\sp\{\sp b\sp\}\sp]=(f_0|W_1)^{-\iota}[\sp\{\sp b\sp\}\sp]=
f_0^{-\iota}[\sp\{\sp b\sp\}\sp]\cap W_1\ssp$. Condition (2) follows from what
we have already shown when noticing $\rho_i\inc\roman{id\,},$ we utilize the
implication $(w_1,b)\in f\imply w_1\in f_0^{\sp-\iota\sp}[\sp\{\sp b\sp\}\sp]
\cap W_1\inc\varSigma_1\ssp$. Theorem 4.3 gives $(E\ssp,F\sp,g\sp)\in\CinftyPi\ssp$.

To prove (x)\ssp, for $x\ar 1=(\sp\eta\ar 1\sp,\varphi\ar 1)\ssp,$ it suffices
to observe that we have\newline$\mhyppy{20}
(\sp x\ar 1\sp,\bar\eta\ar 1\sn+y\ar 1)\in\varSigma\cap W\equivv
(\sp x\ar 1\sp,y\ar 1)\in f_0^{\sp-\iota\sp}[\sp\{\sp b\sp\}\sp]\cap
W\ssn\ar 1=g\ssp$.                                                    \newQED

If we insisted on expressing the content of Theorem 2 in a language very akin
to the standard formalism, we would be lead for example to the following
a bit clumsy formulation: {\it Let $I=[\,0\,,1\,]$ and let $O$ be open in $I
\times\Re^{\,2}$. Suppose the $C^1$ function $y_1:I\to\Re$ and the smooth $
\varphi_1:O\to\Re$ satisfy\par\noin{\rm
(1)} \ $y_1(0)=\eta_1$ and $(s,y_1(s),y_1'(s))\in O$ and $\varphi_1(s,y_1(s),
                                                y_1'(s))=0$ for all $s\in I,$\par\noin{\rm
(2)} \ $\partial_3\varphi_1(s,y_1(s),y_1'(s))\not=0$ for all $s\in I$.

Then there are $\eps>0$ and an open neighborhood $V$ of $\varphi_1$ in $
\Cinfty(O)\ssp,$ and a neighborhood $V_1$ of $y_1$ in $C^1(I)$ such that for
each $x=(\eta,\varphi)$ with $|\eta-\eta_1|<\eps$ and $\varphi\in V,$ there is
a unique $y\in V_1$ with $y(0)=\eta$ and $(s,y(s),y'(s))\in O$ and $\varphi
(s,y(s),y'(s))=0$ for all $s\in I$. This assignment $x\mapsto y$ defines a map\vskip.5mm\centerline{$
\Re\times\Cinfty(O)\iinc\{\,(\eta,\varphi):|\eta-\eta_1|<\eps$ and $\varphi\in 
V\,\}\to\Cinfty(I)\,,$}\par
\noin
which is smooth in the Keller $\CinftyPi$ sense.} (Taking into account that
the spaces in question are all Fr\'echet, smoothness then holds also in
various other senses, cf.\ [\,Hi1\ssp; p.\ 241, Remarks 0.12\,]\ssp.)

We refrain from expressing any of our further results in this kind of longish
form, thus leaving the possible formulations to the reader if one such wishes.
Next modifying Example 1 a little, we consider a boundary value problem.

\NSN{\bf
3} \Examplee We study ($*$) \ $y\ssp''=\varphi\circ[\,\roman{id\,};y\ssp,
  y\ssp'\,]$ with $y\ssp(0)=\eta\ar 0$ and $y\ssp(1)=\eta\ar 1$ for $C^{\ssp 2}$
functions $y:I=[\,0\,,1\,]\to\Re\,$. Here we assume $\eta=(\sp\eta\ar 0\ssp,
\eta\ar 1)\in\Re\times\Re\,,$ and $\varphi:O\to\Re$ smooth as in Example
1\ssp. The solution relation $\varSigma$ of our present equation is the set of
all $w=(x,y)=(\eta,\varphi,y)=(\eta_0,\eta_1,\varphi,y)$ with ($*$) and $
\eta_0,\eta_1,\varphi,y$ as we have just specified.

Now putting $E=\biit R\sp\sqcap\biit R\sp\sqcap\Cinfty(O)\ssp,$ and letting $
{F_{\sn}}_i$ be the topological linear subspace of $G=C^{\,i\sp+\sp 2\sp}(I\sp
)$ with $\vecs F_i=\{\,y\in\vecs G:y\ssp(0)=y\ssp(1)=0\,\}\ssp,$ we here call
$w\in\varSigma$ {\it regular\ssp} if{}f the equation $u''-p_3u'-p_2u=v''$ has
a unique solution $u\in\vecs F_0$ for any $v\in\vecs F_0\ssp$. Letting $
\varSigma\ar 0$ be the set of all regular $w\ssp,$ we prove

\newProCla 4 Theorem.

For every $w\in\varSigma\ar 0$ there is some $W$ with $w\in W\in\taurd(E\sqcap
C^{\ssp 2}(I\sp))$ and $(E\ssp,\Cinfty(I\sp)\ssp,\varSigma\cap W\sp)\in\CinftyPi\ssp$.
                                
\Prooff We first observe that $y\mapsto y''=z$ is a linear homeomorphism $F_0
\to C\ssp(I\sp)$ with inverse $\ell:z\mapsto y\ssp$, where $y(s)=z\yr 2(s)-
z\yr 2(1)\ssp s$ with $z^2(s)=\int_0^s(\int_0^tz)\,dt\ssp$. Also $\ell\,|\,(\sp
\vecs C^{\ssp i\sp}(I\sp)):C^{\ssp i\sp}(I\sp)\to F_i$ is a linear
homeomorphism for all $i\in\No\ssp$. For $\eta=\mathbreak(\sp\eta\ar 0\ssp,
\eta\ar 1)\ssp$, \hfill now \hfill defining \hfill $\bar\eta=\seq{\,\eta\ar 0\sn
+(\sp\eta\ar 1\sn-\eta\ar 0)\ssp s:s\in I\,}\ssp$, \hfill we \hfill let \hfill $
\tilde f_{0i}=\mathbreak(E\sqcap F_i,F_i,f_{0i})\ssp$, where $
f\sn\LHB{.2}{\ar 0}\LHB{.33}{_i}\sp(w)=y-\ell\,(\sp\varphi\circ
[\,\roman{id\,};y+\bar\eta\ssp,y\ssp'+\bar\eta\ssp'\,]\sp)$ for $\,w=(x,y)=
(\eta,\varphi,y)=(\eta_0,\eta_1,\varphi,y)\in\vecs(E\sqcap F_i)$ such that $
\rng[\,\roman{id\,};y+\bar\eta\ssp,y\ssp'+\bar\eta\ssp'\,]\inc O\ssp$.

Taking $b=I\sn\times\sn\{\sp 0\sp\}\ssp,$ the equivalence \ $(\sp\eta\ssp,
\varphi\ssp,\bar\eta+y\sp)\in\varSigma\equivv(\sp\eta\ssp,\varphi\ssp,y\sp)\in
f_{0i}^{\sp-\iota\sp}[\sp\{\sp b\sp\}\sp]$\newline now holds for all $
\eta\ssp,\varphi\ssp,y\ssp,$ and all $i\in\No\ssp,$ since an induction based
on ($*$) shows that we have $\rng\varSigma\inc\Cinfty(I\sp)\ssp$. Theorem 3.6
in conjunction with the chain rule\,\footnote{\,and the fact: $(E\ssp,F\sp,
  f\sp)\ssp,(E\ssp,G\sp,g\sp)\in C^{\ssp k}_\pi\imply(E\ssp,F\sqcap G\sp,[\,
  f\sp,g\,]\sp) \in C^{\ssp k}_\pi\ssp,$\newline\hyppy{2mm}see [\,Hi1\ssp;
  Prop.\ 0.10, p.\ 240\,]} 
and $\cal L\inc\CinftyPi$ gives us $\{\,\tilde f_{0i}:i\in\No\ssp\}\inc
\CinftyPi\ssp$. \hfill From Remarks 3.7\sp(c)\ssp, we obtain\linebreak $
\partial\ar 2\sp f\sn\ar 0\LHB{.3}{_i}\sp(w)\ssp u=u-\ell\,(\sp p\ar 3\sp
u\ssp'\sn+p\ar 2\sp u\sp)\ssp$. \hfill We see that $\partial_2f_{0i}(w):\vecs
F_i\to\vecs F_i$ is bijective\linebreak for $w\in\varSigma\ar 1=\{\ssp(\sp\eta
\ssp,\varphi\ssp,y\sp):(\sp\eta\ssp,\varphi\ssp,\bar\eta+y\sp)\in
\varSigma\ar 0\ssp\}\ssp$.

Now fixing $w=(\sp\eta\ar 0\ssp,\varphi\ar 0\ssp,y\ar 0)\in\varSigma\ar 0\ssp
$, \hfill we apply Theorem 4.1 to $\tilde f\sNorr{\sixroman00}$ at $(\sp
w\ar 0\ssp,b\ssp)$ for $w\ar 0=(\sp\eta\ar 0\ssp,\varphi\ar 0\ssp,y\ar 0-
\bar\eta\ar 0)\ssp$. We then obtain some $W\ssn\ar 0$ with $w\ar 0\in
W\ssn\ar 0\in\taurd(E\sqcap F\sNorr 0)$ and\Newline $(E\ssp,F\sNorr 0\ssp,
f_{\sixroman00}^{\sp-\iota\sp}[\sp\{\sp b\sp\}\sp]\cap W\ssn\ar 0)\in\CinftyPi
$ and $\partial\ar 2\ssp f\sNorr{\sixroman00\sp}(\sp w\ar 1)$ bijective $\vecs
F\sNorr 0\to\vecs F\sNorr 0$ for $w\ar 1\in$ $W\ssn\ar 0\ssp$.\Newline Then we
have $g=f_{\sixroman00}^{\sp-\iota\sp}[\sp\{\sp b\sp\}\sp]\cap W\ssn\ar 0\inc
\varSigma\ar 1\ssp$, \hfill and hence also $
\partial\ar 2\ssp f\sNorr{\sixroman0i}\ssp(\sp w\ar 1):\vecs F\sNorr{\snn i}
\to\vecs F\sNorr{\snn i}$\linebreak is bijective for $w\ar 1\in g\ssp$. A
continuous linear map $\ssp\ell\ar 0\sn:E\sqcap C^{\ssp 2}(I\sp)\to E\sqcap
F\sNorr 0$ is given by $(\sp\eta\ssp,\varphi\ssp,y\sp)\mapsto(\sp
\eta{_{}}_y\ssp,\varphi\ssp,y-\bar\eta{_{}}_y\sp)\ssp$, where $\eta{_{}}_y=(\sp
y\value 0\,,y\value 1\sp)\ssp$. Taking $W=
\ell\ar 0\KN1\inve\sp[\,W\ssn\ar 0\ssp]\ssp$, we now have $w\in W\in
\taurd(E\sqcap C^{\ssp 2}(I\sp))\ssp$, \hfill and writing $\,\tilde g\ar 1=
(E\ssp,\Cinfty(I\sp)\ssp,\varSigma\cap W\sp)\ssp$, \hfill we are done once we
show $\tilde g\ar 1\in\CinftyPi\ssp$.

To prove $\tilde g\ar 1\in\CinftyPi\ssp,$ we first construct the setting for
application of Theorem 4.3\ssp.\Newline Let $F$ be the topological linear
subspace of $G=C^{\,\infty}(I\sp)$ with $\vecs F=\{\,y\in\vecs G:y\value 0=
y\value 1=0\,\}\ssp$. We take $\Delta$ and $I\ar 1$ as in the proof of Theorem
2\ssp, and also define $\rho_i$ and $\rho_{ij}$ as there. We put $f=\bigcap
\ssp\{\,f_i:i\in\No\ssp\}\ssp,$ where $f_i=f_{0i}\ssp|\,W_0\ssp$. As before,
one verifies the conditions of Theorem 4.3\ssp, which then gives us $
(E\ssp,F\sp,g\sp)\in\CinftyPi\ssp$. Now $\tilde g\ar 1\in\CinftyPi$ follows
once the reader has verified that the relation $\varSigma\cap W$ is ex- actly
the function $\dom g\owns x=(\sp\eta\ssp,\varphi\sp)\mapsto\bar\eta+g\ssp(x)
\ssp$.                                                                \newQED

As an example, if we have $w=(\eta,\varphi,y)\in\varSigma$ with $\rng p\ar 2=
\{\sp r\}$ and $\rng p\ar 3=\{0\}\ssp,$ an elementary exercise gives $w\not\in
\varSigma\ar 0$ if{}f $r=-n^{\ssp 2\sp}\pi^{\ssp 2}$ for some $n\in\Zepp$. In
fact, when $r=-n^2\pi^2,$ instead of bijectivity of $\ell_w=\partial_2
f_{0i}(w):\vecs F_i\to\vecs F_i\ssp,$ we have $\ell_w^{\sp-\iota\sp}[\sp\{\sp
\bnull F\sn\LHB{.3}{_{_i}}\}\sp]=\{\,\smb C\ssp s_n:\smb C\in\Re\,\}\ssp,$
where $s_n=\seq{\,\sin\ssp(n\pi t):t\in I\,}\ssp,$ and $\rng\ell_w=$ $\{\,v\in
\vecs F_i:\int_{\,0}^{\,1}v\sn\cdot\sn s_n=0\,\}\ssp$. Thus $\rng\ell_w$ and $
\ell_w^{\sp-\iota\sp}[\sp\{\sp\bnull F\sn\LHB{.3}{_{_i}}\}\sp]$ are
(topologically) comp- lemented subspaces of $F_i\ssp$.

More interesting (and more laborious, of course) results analogous to Theorem
4 might be obtained by following the above lines and, for example, developing
further the study in [\,W\ssp; pp.\ 368\,--\,369\,]\ssp.

In the previous examples we had a "fully nonlinear\ssp" (or implicit) first
order initial value problem and an explicit second order boundary value
problem for an ordinary scalar differential equation. We proved that in both
cases the corresponding {\it solution relation for smooth data\ssp} actually {\it
is a smooth map locally\ssp} around each "regular\ssp" pair (\sp data\ssp,
solution\sp)\ssp. We used Theorems 4.1 and 4.3\ssp, and verification of their
conditions was a relatively simple matter.

Corresponding results for partial differential equations require much more
work since we need to apply various "apriori estimates\ssp". This is already
illustrated in our next example where we study nonlinear initial value problems
corresponding to a particularly simple linear partial differential operator,
namely $y\mapsto{y_{}}_t+{y_{}}_x\ssp$.

\NSN{\bf
5} \Examplee We study the equation ${y_{}}_t+{y_{}}_\theta=\varphi\ssp(\sp
  t\ssp,\theta\ssp,y\ssp(\sp t\ssp,\theta\sp))$ on the closed half-cylinder $
\Repp\sn\times\roman S^1$ with initial condition $y\ssp(\sp 0\,,\theta\sp)=
y_0(\theta)\ssp$. To have a precise setting which can be written in a simple
manner, we interpret the equation for functions $y:\Repp\sn\times\Re\to\Re$ by
imposing the periodicity requirements $y\ssp(\sp t\ssp,\eta+1\sp)=y\ssp(\sp
t\ssp,\eta\ssp)$ and $\varphi\ssp(\sp t\ssp,\eta+1\ssp,\xi\sp)=\varphi\ssp(\sp
t\ssp,\eta\ssp,\xi\sp)\ssp$. The equation then becomes

\noin$(*)\mhyppy{25}\partial\ar 1\sp y+\partial\ar 2\sp y=\varphi\circ
[\,\roman{id\,},y\,]\,$ with $\,y\ssp(\sp 0\,,\cdot)=y\ar 0\ssp$.

We can further rewrite this $y=\bar y\ar 0\sn+\cal I\ssp(\sp 0\,,\varphi\circ
[\,\roman{id\,},y\,]\sp)\ssp,$ where $\bar y_0(t,\eta)=y_0(\eta-t)$\Newline
and $\cal I(a,z)(t,\eta)=\int_{\,a}^{\,t}(z\circ{}?_\tau)(t,\eta)\,d\tau$ with
${}?_\tau(t,\eta)=(\tau,\eta-t+\tau)\ssp$. Indeed,\par$\mhyppy{7.5}
 (y-\bar y_0)(t,\eta)=y(t,\eta)-y_0(\eta-t)=y(t,\eta)-y(0,\eta-t)$\par$\mhyppy{27}
=\big/_{\!\tau=0\,}^{\ t}y\ssp(\sp\tau\sp,\eta-t+\tau\sp)=\int_{\,0}^{\,t}\seq{\,
  y\ssp(\sp\tau\sp,\eta-t+\tau\sp):\tau\in\Repp\ssp}\ssp'$\par$\mhyppy{27}
=\int_{\,0}^{\,t}(\sp\partial\ar 1y\ssp(\sp\tau\sp,\eta-t+\tau\sp)+
 \partial\ar 2\sp y\ssp(\sp\tau\sp,\eta-t+\tau\sp))\,d\ssp\tau$\par$\mhyppy{27}
=\int_{\,0}^{\,t}\varphi\ssp(\sp\tau\sp,\eta-t+\tau\sp),
                       y\ssp(\sp\tau\sp,\eta-t+\tau\sp))\,d\ssp\tau$\par$\mhyppy{27}
=\int_0^t(\varphi\circ[\,\roman{id\,},y\,]\circ{}?_\tau)(t,\eta)\,d\tau
=\cal I(0,\varphi\circ[\,\roman{id\,},y\,]\sp)(t,\eta)\ssp,$

\noin from which one easily deduces the equivalence $(*)\equivv y=\bar y\ar 0\sn
+\cal I\ssp(\sp 0\,,\varphi\circ[\,\roman{id\,},y\,]\sp)\ssp,$ assuming for
example that $\varphi$ is smooth and $y$ is of class $C^{\ssp 1},$ as we shall do.

If we have a locally convex space $S\ssp(\ldots\Re\sn\times\ldots)$ of
functions $y:\ldots\Re\sn\times\ldots\to$ $\Re\,,$ we denote by $S\ssp(\ldots
\roman S\yr 1\ssn\times\ldots)$ its topological linear subspace with
underlying set formed by the $y$ with $y\ssp(\ldots\,\eta\ssp,\ldots)=y\ssp(
\ldots\,\eta+1\ssp,\ldots)$ for all $\eta\in\Re\,$.

{\def\leLCS-{{\le}{}_{_{\roman{LCS}}}\text{\sp-\sp}}
For $i\in\No$ and an interval $I\inc\Re\ssp,$ we write $\,C^{\,0\sp,\ssp i\sp}
(\sp I\sn\times\sn\Re\sp)=$\vskip.5mm\centerline{$
=\leLCS-\inf\ssp\{\,G:C\ssp(\sp I\sn\times\sn\Re\sp)\le G$ and $\all{l\in
i\sp\yplus\!}\,\roman p^{\,l}\KN{1.3}\LHB{.5}{\ai G}\in\cal L\ssp(\sp G\sp,
C\ssp(\sp I\sn\times\sn\Re\sp))\ssp\}\ssp,$}

\noin where $\roman p^{\,l}_G=\{\ssp(\sp y\ssp,\partial_2^{\,l}y\sp):y\in\vecs
G\,\}$ with $\partial_2^{\,l}y$ defined by the recursion $\partial_2^{\,0}y=y$
and $\partial_2^{\,k\sp+1}y=\partial_2(\partial_2^{\,k}y\sp)\ssp$. \hfil Generally,
the space $C^{\,0\sp,\ssp i\sp}(\sp I\sn\times\sn\Re\sp)$ is Fr\'echet. \hfil For $I
$ compact, we have $C^{\,0\sp,\ssp i\sp}(\sp I\sn\times\sn\roman S\yr 1)$ a
Banach space, and as a canonical norm for it we use\par\centerline{$
y\mapsto\|\ssp y\ssp\|\sp\LHB{.2}{\ai C}\LHB{.3}{_{^{0,i}}}=\sup\ssp\{\,
|\ssp r\sp|:\exi{l\in i^+}\,r\in\rng(\partial_2^{\,l}y\sp)\ssp\}\ssp$.}

We put $E_0=\Cinfty(\roman S^1)\ssp,$ and let $E_1$ be the topological linear
subspace of $E_{_A}=\Cinfty(\Repp\sn\times\roman S^1\times\Re)$ with $\vecs
E_1$ formed by the $\varphi$ such that every $T\in\Repp$ has some $M<\infty$
with $\varphi(t,\eta,\xi)=0$ for $t\le T$ and $M\le|\xi|\ssp$. Note that $
\vecs E_1$ is not a $\taurd E_{_A}\ssp$--\,closed set, whence $E_1$ is not a
Fr\'echet space. We put $E=E_0\sqcap E_1\ssp,$ and $F=\Cinfty(\Repp\sn\times
\roman S^1)\ssp$. We consider the solution relation $\Sigma$ of $(*)$ having
as members exactly the $(x,y)=(y_0,\varphi,y)$ such that $x=(y_0,\varphi)\in
\vecs E$ and $y\in\vecs C^1(\Repp\sn\times\roman S^1)$ and $y=\bar y_0+
\cal I(0,\varphi\circ[\,\roman{id\,},y\,]\sp)\ssp$. We now have

\newProCla 6 Theorem.

The membership $(E\ssp,F\sp,\Sigma\sp)\in\CinftyPi$ with $\dom\Sigma=\vecs E$ holds.

\Prooff Let $b=\bnull F$ and $g=\Sigma\,$. Taking $\Delta=\{\,(\sp i\ssp,\sp
j\sp):i\ssp,\sp j\in\N$ and $i\le j\,\}\ssp,$ we apply Theorem 4.3 with $F_i=
C^{\,i}(\sp[\,0\,,i\,]\times\roman S^1)$ and the functions $f_i$ and $f$
defined by the prescription $(x,y)=(y_0,\varphi,y)\mapsto y-\bar y_0-\cal I(0,
\varphi\circ[\,\roman{id\,},y\,]\sp)\ssp$. The functions $\rho_{ij}$ and $
\rho_i$ we define by $y\mapsto y\,|\,(\sp[\,0\,,i\,]\sn\times\sn\Re\ssp)\ssp$.
Since $v\mapsto\cal I\ssp(\ssp 0\,,v\sp)$ gives a continuous linear map $F_i
\to F_i$ for $i\in I_1=\N\sp,$ by Theorem 3.6 and Proposition 3.1\ssp, we then
have $(E\sqcap F_i,F_i,f_i)\in\CinftyPi\ssp$. Lemmas c and e below together
give $\dom\Sigma=\vecs E\ssp$. It remains to prove $\tilde g=(E\ssp,F\sp,
\Sigma\sp)\in\CinftyPi\ssp,$ for which we must verify conditions (1) and (2)
of Theorem 4.3\ssp.

For $g=f\sp\inve\sp[\sp\{\sp b\sp\}\sp]$ in (1)\ssp, we trivially have "$
\supseteq\ssp$", and the converse follows from Lemma e below. That ${g_{}}_i$
is a function follows from Lemma a below, and we have $\dom{g_{}}_i\inc\vecs E
=\dom g$ trivially. For condition (2) from Remarks 3.7\ssp(b) with $a=
\partial\ar 3\sp\varphi\circ[\,\roman{id\,},y\,]\ssp,$ we obtain $
\partial\ar 2\sp{f_{\sn}}_i\sp(\sp x\ssp,{\rho_{}}_i\sp y\sp)\,u=u-\cal I\ssp
(\ssp 0\,,a\cdot u\sp)\ssp,$ whence Lemmas d\ssp, e\ssp, and a show that (2)
holds. Now Theorem 4.3 gives $\tilde g\in\CinftyPi\ssp$.              \newQED

For the following lemmas, when $I=\Repp$ or $I=[\,t\ar 0\ssp,t\ar 0+
\smb T\,\sp]$ with $t\ar 0\ssp,\smb T\in\Repp\ssp,$ writing $S=\vecs C\ssp(\sp
I\sn\times\sn\roman S\yr 1\ssn\times\Re\sp)\ssp$, we let $S\ar 3\sp(I\sp)=\{\,
\varphi:\varphi\ssp,\partial\ar 3\sp\varphi\in S\,\}\ssp$. We further put $
S\ar 2\sp(I\sp)=\{\,\varphi:\varphi\ssp,\partial\ar 2\sp\varphi\,,
\partial\ar 3\sp\varphi\,,\partial\ar 3\sp\partial\ar 2\sp\varphi\,,
\partial\ar 3\sp\partial\ar 3\sp\varphi\in S\,\}\ssp,$ and let $S\ar 1(I\sp)$
be the subset of $S\ar 2\sp(I\sp)$ formed by $\varphi$ such that every $t\ar 1
\in I$ has some $M<\infty$ with $\varphi(t,\eta,\xi)=0$ for $t\le t\ar 1$ and
$M\le|\xi|\ssp$. Then $\vecs E_1=S_1(\Repp)\cap\vecs\Cinfty
(\Repp\!\times\sn\Re\sn\times\sn\Re\sp)\ssp$.

\newpage 

\newProCla a Lemma.

{\rm(uniqueness)} Let $\varphi\in S\ar 3\sp(I\sp)$ and $a\,,u\ssp,v\in\vecs
C^{\,0\sp}(\sp I\sn\times\sn\roman S\yr 1)\ssp,$ where $I=\Repp$\Newline or $I
=[\,0\,,\smb T\,]$ with $\smb T\in\Rep$. Then $\,[\ u-v=\cal I\ssp(\sp 0\,,
\varphi\circ[\,\roman{id\,},u\,]\sp)-\cal I\ssp(\sp 0\,,\varphi\circ
[\,\roman{id\,},v\,]\sp)\imply$\Newline $u=v\ ]$ and $\,[\ u=\cal I\ssp(\sp
0\,,a\cdot u\sp)\imply\rng u=\{\sp 0\sp\}\sp\,]$ are valid.

\Prooff Let \hfil $y=u-v=\cal I\ssp(\sp 0\,,\varphi\circ[\,\roman{id\,},u\,]-
\varphi\circ[\,\roman{id\,},v\,]\sp)$ \hfil for the first implication.\linebreak
Defining $a(t,\eta)=\int_0^1(\partial_3\varphi\circ[\,\roman{id\,},s\,u+(1-s)
\,v]\sp)(t,\eta)\,ds\ssp,$ we then have

$y(t,\eta)=\int_0^t(\varphi\circ[\,\roman{id\,},u]\circ
{}?_\tau(t,\eta)-\varphi\circ[\,\roman{id\,},v]\circ{}?_\tau(t,\eta))\,d\tau$

$\mhyppy{9}=\int_0^t\int_0^1((\partial_3\varphi\circ[\,\roman{id\,},\,
su+(1-s)v]\cdot y)\circ{}?_\tau)(t,\eta)\,ds\,d\tau
=\cal I(0,a\cdot y)(t,\eta)\ssp$.

Fixing $L\in I\sp,$ and writing $\sigma(t)=\sup\ssp\{\,|\,y(t,\eta)\ssp|:\eta
\in\Re\,\}\ssp,$ it suffices to show that $\sigma(t)=0$ for all $t\in J=
[\,0\,,L\,]\ssp$. Put $M=\sup\ssp\{\,|\sp r\sp|:r\in a\,[\,J\times\Re\,]\,\}
\ssp,$ and note that $\sigma$ is continuous. By $y=\cal I(0,a\cdot y)\ssp,$ we
have $\sigma(t)\le M\int_{\,0}^{\,t}\sigma\ssp,$ hence $\frac d{dt}(e^{-Mt}
\int_{\,0}^{\,t}\sigma)\le 0\,,$ and consequently $\int_{\,0}^{\,t}\sigma\le 0
\,,$ whence $\sigma(t)=0$ follows. This deduction at the same time verified
the second implication.                                                  \QED

\newProCla b Lemma.

{\rm(apriori estimates)} With $\smb L\ssp,\smb R\ssp,t\ar 0\in\Repp$ let $I=
[\,t\ar 0\ssp,t\ar 0\sn+\smb L\,]\ssp,$ and let $\varphi\in S_2(I)$ and $F_i=
C^{0,i}(I\times\roman S^1)$ with $i=0\,,1\ssp$. Writing\newline$\mhyppy{5}
\smb M\ar 0=\sup\ssp\{\,|\ssp r\sp|:r\in\varphi\,[\,I\sn\times\sn\Re\sn\times\sn\{\sp
    0\sp\}\,]\cup\partial\ar 3\sp\varphi\,[\,I\sn\times\sn\Re\sn\times\sn\Re
    \sp\,]\,\}\,,$\newline$\mhyppy{5}
\smb M\ar 1=\sup\ssp\{\,|\ssp r\sp|:r\in\varphi\,[\,I\sn\times\sn\Re\sn\times\sn\{\sp
    0\sp\}\,]\cup(\partial\ar 2\sp\varphi\cup\partial\ar 3\sp\varphi\sp)\,[\,
    I\sn\times\sn\Re\sn\times\sn\Re\sp\,]\,\}\,,$ \,and\newline$\mhyppy{6.4}
\smb M  =\sup\ssp\{\,|\ssp r\sp|:r\in(\sp\partial\ar 3\sp\varphi\cup\partial\ar 2\sp
    \partial\ar 3\sp\varphi\cup\partial\ar 3\sp\partial\ar 3\sp\varphi\sp)\,[\,
    I\sn\times\sn\Re\sn\times\sn\Re\sp\,]\,\}\,,$ \,then hold\vskip.5mm\noin{\rm
(1)} \ if $t\ar 0=0$ and $v\ssp,y\in\vecs F_i$ and\newline$\mhyppy{7}
   y=v+\cal I\ssp(\sp 0\,,\varphi\circ[\,\roman{id\,},y\,]\sp),$ then $
   \|\ssp y\ssp\|_{F_i}\le(\sp 1+\|\ssp v\ssp\|_{F_i})\,e^{\,\ssmb{L\sp M}_i}-1\ssp,$\vskip0mm\noin{\rm
(2)} \ if $u\ssp,v\in\bar B_{F_1}(\smb R)\ssp,$ then\newline$\mhyppy{7}
   \|\,\cal I\ssp(\sp t\ar 0\ssp,\varphi\circ[\,\roman{id\,},u\,]\sp)-
   \cal I\ssp(\sp t\ar 0\ssp,\varphi\circ[\,\roman{id\,},v\,]\sp)\ssp\|_{F_1}
   \le\smb{L\,M}\,(\sp 2+\smb R)\,\|\,u-v\,\|_{F_1}\,$.

\Prooff We first note that for (1) we may assume $\smb M\ai i\not=\infty\ssp,$
since in the case $\smb L\not=0$ the assertion otherwise is trivial, and a
moment's thought shows that it also holds if $\smb L=0\,$. Likewise, for (2)
we may assume $\smb M\not=\infty\ssp$. Now to prove (1)\ssp, we observe that $
z=\partial\ar 2\sp y$ in the case $i=1$ satisfies the equation $z=c+
\cal I\ssp(\sp 0\,,a\sn\cdot\sn z+b\ssp)\ssp,$\Newline where $c=
\partial\ar 2\sp v$ and $b=\partial\ar 2\sp\varphi\circ[\,\roman{id\,},y\,]$
and $a=\partial\ar 3\sp\varphi\circ[\,\roman{id\,},y\,]\ssp$. Writing $
\smb C\ai i=\|\ssp v\ssp\|_{F_i}$ and $\sigma\ar 1(t)=\sup\ssp\{\,
\sigma\ar 0\sp(t)\ssp,\sigma\ssp(t)\ssp\}\ssp,$ where $\sigma\ar 0\sp(t)=\sup\ssp
\{\,|\,y\ssp(\sp t\ssp,\eta\sp)\ssp|:\eta\in\Re\,\}$ and $\sigma\ssp(t)=\sup\ssp
\{\,|\,z\ssp(\sp t\ssp,\eta\sp)\ssp|:\eta\in\Re\,\}\ssp,$ we have $
\|\ssp y\ssp\|_{F_i}=\sup\ssp\{\,\sigma\ai i(t):t\in I\,\}\ssp,$ whence it
suffices to prove $\sigma\ai i(t)\le(\sp 1+\smb C\ai i)\,e^{\,\ssmb{L\sp M}_i}-1\ssp$.

From $z=c+\cal I\ssp(\sp 0\,,a\sn\cdot\sn z+b\ssp)\ssp,$ we get $\sigma\ssp(t)
\le\smb C\ar 1+\smb M\ar 1\int_{\,0}^{\,t}(\sp\bold l+\sigma\sp)\ssp$, where
$\bold l=I\sn\times\sn\{1\}\ssp$.\Newline Using $\varphi\ssp(\sp t\ssp,
\eta\ssp,y\ssp(\sp t\ssp,\eta\sp))=\varphi\ssp(\sp t\ssp,\eta\ssp,0\sp)+
\int_{\,0}^{\,1}\partial\ar 3\sp\varphi\ssp(\sp t\ssp,\eta\ssp,s\,y\ssp(\sp
t\ssp,\eta\sp))\,y\ssp(\sp t\ssp,\eta\sp)\,ds\ssp,$ we get $\sigma\ar 0\sp(t)$\Newline
$\le\smb C\ar 0+\smb M\ar 0\int_{\,0}^{\,t}(\sp\bold l+\sigma\ar 0)$ \hfill
from \hfill $y=v+\cal I\ssp(\sp 0\,,\varphi\circ[\,\roman{id\,},y\,]\sp)\ssp$. \hfill
\hfill Altogether, \hfill these \hfill give\linebreak $\sigma\ai i(t)\le
\smb C\ai i+\smb M\ai i\int_{\,0}^{\,t}(\sp\bold l+\sigma\ai i)\ssp,$ whence $
\seq{\,e^{-\ssmb M_it}\int_{\,0}^{\,t}(\sp\bold l+\sigma\ai i):t\in
I\,}\sp'(t)\le(\sp 1+\smb C\ai i)\,e^{-\ssmb M_it}$,\Newline and consequently
$\int_{\,0}^{\,t}(\sp\bold l+\sigma\ai i)\le e^{\ssp\ssmb M_it}
\int_{\,0}^{\,t}(\sp 1+\smb C\ai i)\,e^{-\ssmb M_is}ds$ for $t\in I\sp$. Hence
finally\par$\mhyppy{13}
\sigma\ai i(t)\le\smb C\ai i+\smb M\ai ie^{\ssp\ssmb M_it}\int_{\,0}^{\,t}(\sp
  1+\smb C\ai i)\,e^{-\ssmb M_is\,}d\ssp s$\newline$\mhyppy{24}
= \smb C\ai i+e^{\ssp\ssmb M_it}(\sp 1+\smb C\ai i)\,(\sp 1-e^{-\ssmb M_it})
= (\sp 1+\smb C\ai i)\,e^{\ssp\ssmb M_it}-1\ssp$.
                                                
To prove (2)\ssp, putting $y=\cal I(t_0,\varphi\circ[\,\roman{id\,},u\,]\sp)-
\cal I(t_0,\varphi\circ[\,\roman{id\,},v\,]\sp)$ with $u,v$ fixed, note that
we can write

$y\ssp(t,\eta)=\int_{t_0}^t\int_0^1(\partial_3\varphi\circ[\,\roman{id\,},s\,u
+(1-s)\,v\,]\cdot(u-v))\ssp(\tau,\eta-t+\tau)\,ds\,d\tau\,,$

\newpage 

\noin and hence \ $\partial_2y\ssp(t,\eta)=\int_{t_0}^t\int_0^1(\partial_2
\partial_3\varphi\circ[\,\roman{id\,},s\,u+(1-s)\,v\,]\cdot(u-v)$

$\mhyppy{10}+\partial_3^2\varphi\circ[\,\roman{id\,},s\,u+(1-s)\,v\,]
\cdot(s\,\partial_2u+(1-s)\,\partial_2v)\cdot(u-v)$

$\mhyppy{10}+\partial_3\varphi\circ[\,\roman{id\,},s\,u+(1-s)\,v\,]\cdot(
\partial_2u-\partial_2v))\ssp(\tau,\eta-t+\tau)\,ds\,d\tau\,$.

\noin From these we get $\sigma\ar 0\sp(t)\le\smb{L\ssp M}\,\|\,u-v\,\|_{F_0}$
and $\sigma\ssp(t)\le\smb{L\ssp M}\,(\|\,u-v\,\|_{F_0}+{}$\newline $\smb R\ssp
\|\,u-v\,\|_{F_0}+\|\,u-v\,\|_{F_1})\le\smb{L\ssp M}\,(\sp 2+\smb R)\,
\|\,u-v\,\|_{F_1}\ssp,$ and (2) then follows.                            \QED

\newProCla c Lemma.

{\rm(existence)} Let $I=\Repp$ or $I=[\,0\,,\smb T\,\sp]$ with $\smb T\in\Rep.
$ Then all $v\in C$ $=\vecs C^{\sp\,0\sp,\sp 1\sp}(\sp I\sn\times
\roman S\yr 1)$ and $\varphi\in S\ar 1\sp(I\sp)$ have $y\in C$ with $y=v+
\cal I\ssp(\sp 0\,,\varphi\circ[\,\roman{id\,},y\,]\sp)\ssp$.

\Prooff Fix $v\ssp,\varphi\,,$ and for $t\in I\sp,$ let $S(t)$ be the set of
all $z\in\vecs C^{\sp\,0\sp,\sp 1\sp}(\sp[\,0\,,t\,]\times\roman S^1)$ with
$z=v+\cal I(0,\varphi\circ[\,\roman{id\,},z\,]\sp)\ssp$. By Lemma a\ssp, then
$S(t)=\emptyset$ or singleton. Writing $T=\{\,t\in I:S(t)\not=\emptyset\,\}$
and $t\ar 1=\sup T,$ it suffices to prove $t\ar 1=\sup I$ and that $T$ is
closed, since then $y=\bigcup\bigcup\ssp\{\,S(t):t\in I\,\}$ is the required
solution. Noticing that $T$ is an interval, to proceed indirectly, we assume $
t_1<\sup I$ or $t_1=\sup I\not\in T,$ and show that a contradiction follows.
We first observe that $t_1\ge 0\,,$ since we have $v\,|\,(\{\sp 0\sp\}\sn
\times\sn\Re\sp)\in S(0)\ssp$. It suffices to show that a contradiction
follows in each of the following cases: (a) \ $t_1=0\,,$ (b) \ $0<t_1<\sup I
\sp,$ (c) \ $t_1=\sup I\in\Re\setminus\sn T$.

Writing $I_1=I\cap[\,0\,,t\ar 1\sn+1\,]$ and $\Phi=\bigcup\ssp\{\,
\partial_3^{\,j}\partial_\nu^{\,i}\varphi:i\ssp,\sp j=0\,,1\ssp;\nu=2\ssp,3\,
\}\ssp,$ we put $\smb A=\|\,v\,|\,(\sp I\sn\ar 1\ssn\times\sn\Re\sp)\,\|\sp
\LHB{.2}{_{C^{\,0\sp,1}}}$ and $\smb M=\sup\ssp\{\ssp 1+|\ssp r\sp|:r\in\Phi\,
[\,I\ar 1\ssn\times\sn\Re\sn\times\sn\Re\sp\,]\,\}\ssp$. Now with $\smb R=4\,
(\sp 1+\smb A)\,(\sp 1+\smb M\,t\ar 1\ssp e^{\,M\sp t_1})\ssp,$ we take $J=[\,
t\ar 0\ssp,t\ar 2\ssp]\ssp,$ where $t_0$ and $t_2$ in the above three cases
are defined according to\par\noin
(a) \ $t\ar 0=0$ and $t\ar 2=\smb L=\inf\ssp\{\ssp\sup I\sp,\frac13\ssp(\smb{M
    \ssp R})^{-1}\sp\}\ssp,$\par\noin
(b) \ $t\ar 0=t\ar 1\sn-\smb T\sn\ar 1$ and $t\ar 2=t\ar 1\sn+\smb T\sn\ar 1
    \ssp,$ where $\smb T\sn\ar 1=\inf\ssp\{\,t\ar 1\sp,\sup I-t\ar 1\sp,\frac
    16\ssp(\smb{M\ssp R})^{-1}\sp\}\ssp,$\par\noin
(c) \ $t\ar 0=t\ar 1\sn-\smb L$ and $t\ar 2=t\ar 1\ssp,$ where we now have $
    \smb L=\inf\ssp\{\ssp t\ar 1\sp,\frac13\ssp(\smb{M\ssp R})^{-1}\sp\}\,$.

Then we pick $y\in S\sp(\sp t\ar 0)\ssp,$ define $z^0(t,\eta)=v(t,\eta)+
\int_{\,0}^{\,t_0}(\varphi\circ[\,\roman{id\,},y\,]\circ{}?_\tau)(t,\eta)\,
d\tau$ for $(t,\eta)\in J\sn\times\sn\Re\,,$ and put $F\sn\RHB{.3}{_{_!}}=
C^{0,1}(J\times\roman S^1)\ssp$. If now $z=z^0+\cal I(t_0,\varphi\circ
[\,\roman{id\,},z\,]\sp)$ for some $z\in\vecs F\sn\RHB{.3}{_{_!}}\ssp,$ a
simple exercise shows that we then obtain $y\cup z\in S\sp(\sp t\ar 2)\ssp,$ a
contradiction. Thus it suffices to construct $z\ssp$.

We consider the function $\varrho:\vecs F\sn\RHB{.3}{_{_!}}\to\vecs
F\sn\RHB{.3}{_{_!}}$ given by $z\mapsto z^0+\cal I(t_0,\varphi\circ
[\,\roman{id\,},z\,]\sp)\ssp$. If we can show that $\varrho$ is a contractor $
\bar B\ai F\LHB{.4}{_{_!}}(\smb R)\to\bar B\ai F\LHB{.4}{_{_!}}(\smb R)\ssp,$
the Banach fixed point theorem gives $z$ with $(z,z)\in\varrho\,,$ and we are
one. Using (1) and (2) of Lemma b\,, the reader may verify that indeed $
\varrho\ssp\,|\,\bar B\ai F\LHB{.4}{_{_!}}(\smb R)$ is a contractor.

As a hint in this verification we mention the following. With $\alpha=
\smb{L\,M\,}(\sp 2+\smb R)\le\frac12\ssp$, it suffices to obtain the
inequalities $\|\,\varrho\,(\bnull F\LHB{.4}{_{_!}})\ssp\|+\alpha\,\smb R\le
\smb R$ and $\|\,\varrho\,(u)-\varrho\,(v)\ssp\|\le\alpha\,\|\ssp u-v\ssp\|\ssp
$. The latter one is just (2) of Lemma b\ssp. With\par$\mhyppy{20}
z\ar 0\sp(\sp t\ssp,\eta\sp)=\int_{\,0}^{\,t_0}\varphi\ssp(\sp\tau\ssp,
\eta-t+\tau\ssp,y\ssp(\sp\tau\ssp,\eta-t+\tau\sp))\,d\sp\tau\,$, \hfil we have\linebreak\centerline{$
 \varrho\,(\bnull F\LHB{.4}{_{_!}})\ssp(\sp t\ssp,\eta\sp)
=v\ssp(\sp t\ssp,\eta\sp)+z\ar 0\sp(\sp t\ssp,\eta\sp)+\int_{\,t_0}^{\,t}
 \varphi\ssp(\sp\tau\ssp,\eta-t+\tau\ssp,0\sp)\,d\sp\tau\,$,}

\noin and to estimate $\|\ssp z\ar 0\sp\|\ssp$, one utilizes (1) of Lemma b\ssp.
Then it is a bit tedious veri- fication to show that also $\,\|\,\varrho\,
(\bnull F\LHB{.4}{_{_!}})\ssp\|+\alpha\,\smb R\le\smb R\,$ holds.        \QED

\newProCla d Lemma.

{\rm(existence)} If $i\in\N\sp,$ and $\sp I=\Repp$ or $I=[\,0\,,\smb T\sp\,]$
with $\smb T\in\Rep,$ then any $\,a\ssp,v\in C=\vecs C^{\sp\,0\sp,\ssp i\sp}
(\sp I\sn\times\sn\roman S\yr 1)$ have $u\in C$ with $u=v+\cal I\ssp(\sp 0\,,
a\cdot u\sp)\,$.

\Prooff By Lemma a\ssp, it suffices to show the assertion for $\smb T\in\Rep$.
For a while assuming it for $i=1\ssp,$ we now prove that the full result $\all
i\,\all{a\,,v}\ldots$ follows then by induction on $i\ssp$. Indeed, suppose we
have the result for fixed $i\in\N\sp,$ and let $u=v+\cal I\ssp(\sp 0\,,a\cdot u
\sp)$ with $\,a\ssp,v\in D=\vecs C^{\sp\,0\sp,\ssp i\ssp+1\sp}(\sp I\sn\times
\sn\roman S\yr 1)\ssp$. Then we have $u\in C\ssp,$ and $z=\partial\ar 2\sp u$
satisfies $z=\partial\ar 2\sp v+\cal I\ssp(\sp 0\,,\partial\ar 2\ssp a\cdot u
\sp)+\cal I\ssp(\sp 0\,,a\cdot z\sp)\ssp,$ whence by the inductive assumption
and by Lemma a\ssp, we have $z\in C\ssp,$ consequently $u\in D\ssp$.

It thus suffices to treat the case $i=1\ssp$. With $F_0=C^{\,\sp 0}(\sp
I\sn\times\sn\roman S\yr 1)\ssp,$ let $\smb A=
\|\ssp v\ssp\|\sp\LHB{.2}{\ai F}\LHB{.5}{\sn_{_0}}$ and define $\chi\ssp(s)=(
\int_{-\infty}^{\,\infty}\chi\ar 0){}^{-1}\int_{-\infty}^{\,s}(\sp\chi\ar 0\sp
(-\sigma\sp)-\chi\ar 0\sp(\sigma))\,d\sp\sigma\ssp,$ where we have $\chi\ar 0
\sp(s)=\exp\,((\sp s-1\sp)^{-1\sp}(\sp s-2\sp)^{-1\sp})$ for $1<s<2\ssp,$ and
$\chi\ar 0\sp(s)=0$ otherwise. Then with $\smb M=$ $2\,\sup\ssp(\rng\chi\ssp'
\sp)\sn\cdot\sn\|\,a\ssp\|\sp\LHB{.2}{\ai F}\LHB{.5}{\sn_{_0}}\ssp,$ we put $
\smb B=(\sp 1+\smb A)\,e\yi{\,M\,T}$.

Now taking $\varphi:(\sp t\ssp,\eta\ssp,\xi\sp)\mapsto\chi\ssp(\sp\smb B^{-1\sp}
\xi\sp)\,a\,(\sp t\ssp,\eta\sp)\,\xi$ in Lemma c\,, we have some $u$ with $u=v
+\cal I\ssp(\sp 0\,,\varphi\circ[\,\roman{id\,},u\,]\sp)\ssp,$ and to verify $
u=v+\cal I\ssp(\sp 0\,,a\cdot u\sp)\ssp,$ it suffices to show $\|\ssp u\ssp
\|_{F_0}\le\smb B\ssp$. This indeed is the case since Lemma b\ssp(1) gives\par$\mhyppy{17}
\|\ssp u\ssp\|_{F_0}\le(\sp 1+\smb A)\,e\yi{\,M\,T}-1<(\sp 1+\smb A)\,
e\yi{\,M\,T}=\smb B\ssp$.                                                \QED

\newProCla e Lemma.

{\rm(regularity)} Let $i\in\N\sp,$ and $\sp I=\Repp$ or $I=[\,0\,,\smb T\sp\,]
$ with $\smb T\in\Rep$. Also\Newline let $\varphi\in\vecs\Cinfty(\sp
I\sn\times\sn\roman S\yr 1\ssn\times\sn\Re\sp)$ and $y\ar 0\in\vecs\Cinfty(\sp
\roman S\yr 1)$ and $a\,,v\in C=\vecs C^{\ssp i\sp}(\sp I\sn\times\sn
\roman S\yr 1)\ssp$. The\Newline implications $\,[\ y=\bar y\ar 0+\cal I\ssp(\sp
0\,,\varphi\circ[\,\roman{id\,},y\,]\sp)\in\vecs C^{\,\sp 0\sp,\ssp 1\sp}(\sp
I\sn\times\sn\roman S\yr 1)\imply y\in\vecs\Cinfty(\sp I\sn\times\sn
\roman S\yr 1)\ ]$ and $\,[\ u\in\vecs C^{\,\sp 0\sp,\ssp i\sp}(\sp
I\sn\times\sn\roman S\yr 1)$ and $v=u-\cal I\ssp(\sp 0\,,a\cdot u\sp)\imply u
\in C\ ]\sp$ then are valid.

\Prooff Let $y=\bar y\ar 0+\cal I\ssp(\sp 0\,,\varphi\circ[\,\roman{id\,},y\,]
\sp)\in\vecs C^{\,\sp 0\sp,\ssp 1\sp}(\sp I\sn\times\sn\roman S\yr 1)\ssp$.
Putting $z=\partial_2y\ssp,$ we then have $z=\bar y_0'+\cal I(0,\partial_2
\varphi\circ[\,\roman{id\,},y\,]+\partial_3\varphi\circ[\,\roman{id\,},y\,]
\cdot z\sp)\ssp$. Using this, we show inductively for $i\in\N$ that we have
$y\in D=\vecs C^{0,i}(\sp I\sn\times\roman S^1)\ssp$. The case $i=1$ being our
assumption, assume the assertion for $i\ssp$. By Lemmas d and a\ssp, we have $
z\in D\ssp,$ whence $y\in\vecs C^{0,i+1}(\sp I\sn\times\roman S^1),$ and the
induction is completed.

To prove $y\in\vecs\Cinfty(\sp I\sn\times\sn\roman S\yr 1)\ssp,$ we (again)
recall H.\ A.\ Schwarz' theorem which (roughly) says that if $
\partial_{{\iota_{\sn}}_1}z$ and $\partial_{{\iota_{\sn}}_1}
\partial_{{\iota_{\sn}}_2}z$ are defined and $\partial_{{\iota_{\sn}}_1}
\partial_{{\iota_{\sn}}_2}z$ is continuous, then $\partial_{{\iota_{\sn}}_2}
\partial_{{\iota_{\sn}}_1}z=\partial_{{\iota_{\sn}}_1}
\partial_{{\iota_{\sn}}_2}z\ssp$. Writing $e_1=(1,0)$ and $e_2=(0,1)$ and $
c\,(\bit\eeta)=\Card(\sp\bit\eeta\sp\inve\sp[\sp\{\sp e\ar 1\}\sp]\sp)\ssp,$
it now suffices that inductively on $k\in\No\ssp,$ we prove

\noin(r$_k)\mhyppy{13}\all{\smb N\sp,\bit\eeta}\,\smb N\in\No$ and $\bit\eeta
\in\{\ssp e\ar 1\sp,e\ar 2\sp\}\,\yi N$ and $c\,(\bit\eeta)\le k\imply$\newline$\mhyppy{31}
\dom\sn(\sp d\ssp(\bit\eeta)\ssp y\sp)=I\sn\times\sn\Re$ and $
d\ssp(\bit\eeta)\ssp y$ is continuous\,.

Above, we already proved (r$_0$)\ssp. So suppose we have (r$_k$)\ssp, and
fix $\bit\eeta\in\{e_1,e_2\}\yi N$ with $c(\bit\eeta)=k+1\ssp$. Considering
any $\bit\eeta_1$ with $c(\bit\eeta_1)=k\ssp,$ we have\vskip.5mm\centerline{$
d(\bit\eeta_1)\partial_1y=d(\bit\eeta_1)(-\partial_2y+\varphi\circ
[\,\roman{id\,},y\,]\sp)=-d\ssp(\sn\seq{\ssp e\ar 2\sp}\conc\bit\eeta\ar 1)
\ssp y+d(\bit\eeta_1)(\varphi\circ[\,\roman{id\,},y\,]\sp)\ssp,$}

\noin where $d(\bit\eeta_1)(\varphi\circ[\,\roman{id\,},y\,]\sp)$ can be
written (via appropriate recursive definition, and inductive proof) as a
finite sum of terms $\partial_1^{l_1}\partial_2^{l_2}\partial_3^{l_3}\varphi
\circ[\,\roman{id\,},y\,]\cdot\prod_{i=1}^{l_3}d(\bit\xxi_i)y\ssp,$ where $
c(\bit\xxi_i)\le k\sp$. By (r$_k$) then $\dom(d(\bit\eeta_1)\partial_1y)=
I\sn\times\sn\Re$ holds and $d(\bit\eeta_1)\partial_1y$ is continuous. Now
applying Schwarz' theorem in a finite induction, we see that we can write $
d(\bit\eeta)y=d(\bit\eeta_1)(\partial_1y)$ for some $\bit\eeta_1$ with $
c(\bit\eeta_1)=k\sp$. The induction is completed.

The latter implication follows by induction similar to the above, once we
first note $\,v=u-\cal I\ssp(\sp 0\,,a\cdot u\sp)\imply\partial\ar 1u=
\partial\ar 1v+\partial\ar 2v-\partial\ar 2u+a\cdot u\,$.             \newQED

To give an application of Theorem 4.5\ssp, we modify Example 5 thus considering
\NSN{\bf
7} \Examplee We take the equation ${y_{}}_t+{y_{}}_\theta=\varphi\ssp(\sp
  t\ssp,\theta\ssp,y\ssp(\sp t\ssp,\theta\sp))$ now on the compact cylinder $
I\times\roman S^1$ with $I=[\,0\,,\smb T\,]$ for $\smb T\in\Rep$ fixed,
initial condition as before. We put \hfil $E=E\ar 0\sqcap E\ar 1$ \hfil with \hfil
$E\ar 0=\Cinfty(\sp\roman S\yr 1\sn)$ \hfil and \hfil $E\ar 1=\Cinfty(\sp
I\times\roman S\yr 1\sn\times\Re\sp)\ssp,$ \hfil and \hfil $F=$\linebreak $
\Cinfty(\sp I\times\roman S\yr 1)\ssp,$ which all are Fr\'echet spaces.
Letting the solution relation $\Sigma$ be the set of $(\sp y\ar 0\ssp,
\varphi\ssp,y\sp)\in\vecs(E\sqcap C^{\ssp 1}(\sp I\times\roman S\yr 1))$ with
$y=\bar y_0+\cal I\ssp(\sp 0\,,\varphi\circ[\,\roman{id\,},y\,]\sp)\ssp,$ we show

\newProCla 8 Theorem.

The membership $(E\ssp,F\sp,\Sigma\sp)\in\CinftyPi$ is valid.

\Prooff The set $\Sigma$ is a function by Lemma a\ssp, and Lemma e gives $\rng
\Sigma\inc\vecs F\sp$. Hence arbitrarily given $z\in\Sigma\ssp,$ it suffices
to obtain $W$ with $z\in W\in\taurd(E\sqcap F\sp)$ and $(E\ssp,F\sp,\Sigma\cap
W\sp)\in\CinftyPi\ssp$. As in the proof of Theorem 6\ssp, with $b=\bnull F\ssp
,$ we have $\Sigma=f\sp\inve\sp[\sp\{\sp b\sp\}\sp]\ssp,$ when taking $F_i=
C^{\,i}(\sp I\sn\times\sn\roman S^1)\ssp,$ we define $f$ and $f_i$ as before.
By Theorem 4.5\ssp, it thus suffices to show that with the notations of
Definition 4.4\ssp, we obtain a {\it B\ssp}$\partial\ar 2\ssp$--\,extension of
$(E\sqcap F\sp,F\sp,f\sp)$ in $\CinftyPi$ around $(\sp z\ssp,b\ssp)\ssp,$ when
we define $\rho_i$ and $\rho_{ij}$ by $y\mapsto y\ssp$. Choosing $i\ar 0=1$
and $W\ssn\ar 0=\vecs(E\sqcap{F_{\sn}}_1)\ssp,$ conditions (1) and (2) of
Definition 4.4 are verified as in the proof of Theorem 6\ssp. Also condition
(3) is satisfied, since by Lemma e we have $y\yr 1\in\vecs F\inc\rng
{\rho_{\sp}}_{1j}$ for all $j\in I_1=\N$ whenever $(\sp x\yr 1\!,y\yr 1\!,
\bnull F)=(\sp x\yr 1\!,y\yr 1\!,\sp{\rho_{\sp}}_1\sp b\ssp)\in{f_{\sn}}_1\ssp
$. We are done.                                                       \newQED

For example, with $\smb T<2$ taking $y\ar 0\sp(\eta)=\frac12$ and $\varphi\ssp
(\sp t\ssp,\eta\ssp,\xi\sp)=\xi^{\,2}$ and $y\ssp(\sp t\ssp,\eta\sp)=
(\sp 2-t\sp)^{-1},$ we have $(\sp y\ar 0\ssp,\varphi\ssp,y\sp)\in\Sigma\,,$
whence by Theorem 8 an open neighborhood $U$ of $(\sp y\ar 0\ssp,\varphi\sp)$
in $E$ exists with $U\inc\dom\Sigma\,$. Put loosely, in particular the
equation in question has a unique solution $z$ for each data $x=(\sp
y\ar 1\sp,\varphi\ar 1)\in U_{},$ and $z$ depends smoothly on $x\ssp,$ the map
$U\owns x\mapsto z$ is smooth.

Of course, the result of Theorem 8 with equally simple proof could also have
been obtained by using Theorem 4.3 instead of 4.5\ssp.                    \NS

\Remarkk Since the spaces $E\ssp,F$ are now Fr\'echet, we have $\tilde g=
(E\ssp,F\sp,\Sigma\sp)\in\CinftyPi$ exactly in case $\Sigma:E\iinc\dom\Sigma
\to F$ is a $\Lip\infty$ map in the sense [\,KM\ssp; p.\ 118, Def.\Newline
12.1\,]\ssp. In principle, we could then prove $\tilde g\in\CinftyPi$ using
[\,FK\ssp; Thm.\ 4.8.5, p.\ 153\,]\ssp, but this would not be very practical.
Indeed, we would be lead to the following.

Arbitrarily fix smooth functions $\hat c_0:\Re\times\Re\to\Re$ and $\hat c_1:
\Re\times I\times\Re\times\Re\to\Re$ with $\hat c_0(\sigma,\eta+1)=
\hat c_0(\sigma,\eta)$ and $\hat c_1(\sigma,t,\eta+1,\xi)=\hat c_1(\sigma,t,
\eta,\xi)\ssp,$ and define $c_0(\sigma)(t,\eta)=\hat c_0(\sigma,\eta-t)$ and $
c_1(\sigma)(t,\eta,\xi)=\hat c_1(\sigma,t,\eta,\xi)\ssp$. Letting $\check y$
be the set of all\Newline pairs $(\sigma,y)\in\Re\times(\sp\vecs F\sp)$ with $y=
c_0(\sigma)+\cal I(0,c_1(\sigma)\circ[\,\roman{id\,},y\,]\sp)\ssp,$ the set $S
=\dom\check y$ should be open, and we should have a $\Lip0$ map $\check y:
\biit R\iinc S\to F$ in the sense referred to above. Further, taking any
compact $K\inc S$ and any bounded set $B$ in $F\sp,$ then also \hfil $\{\,u\in
\vecs F:\exi{\sigma\in K}\,u-\cal I\ssp(\sp 0\,,
\partial\ar 3\ssp c\sp\ar 1(\sigma)\circ
[\,\roman{id\,},\check y\ssp(\sigma)\,]\sn\cdot\sn u\sp)\in B\,\}$ \hfil
should\linebreak be a bounded set in the space $F\sp$.\NS

Theorem 4.5 requires the functions ${\rho_{\sp}}_{i\sp j}$ to be injective. In
the previous example this was trivial by $\vecs C^{\,i\ssp+\sp 1}(\Omega)\inc
\vecs C^{\ssp i}(\Omega)\ssp,$ and verification of condition (3) in De-
finition 4.4 was based on the fact that the solution of a differential
equation is more differentiable than apriori. To present an example where (3)
fails, we take complex scalars, and consider the following

\NSN{\bf
9} \Examplee We take the equation $y\ssp'=\varphi\circ[\,\roman{id\,},y\,]$
  with $y\ssp(0)=y\ar 0\ssp,$ where with the notations of Example 3.9 we have
$\varphi:\Omega=\roman D\sp(1)\sn\times\sn\Ce\to\Ce$ and $y:\roman D\ssp(1)\to
\Ce$ holomorphic. We are interested to know whether $(\sp y\ar 0\ssp,\varphi
\sp)\mapsto y$ defines a holomorphic map $E\to F$ when we put $E=\biit C\sp
\sqcap H\sp(\Omega)$ and $F=H\sp(\sp\roman D\sp(1))\ssp$.

To construct the setting for application of Theorem 4.5\ssp, writing $F_r=
H\ai b\sp(\bar{\roman D}\sp(r))\ssp$, we take $\cal R=\{\,(\sp r\sp,s\ssp,
{\rho_{\sp}}_{r\sp s}\sp):0<r\le s<1\,\}$ and $\cal F=\seq{\,({F_{}}_r\ssp,
{f_{}}_r\ssp,{\rho_{\sp}}_r\sp):0<r<1\,}\ssp,$ where ${\rho_{\sp}}_{r\sp s\,}y
={\rho_{\sp}}_r\ssp y=y\,|\,\bar{\roman D}\sp(r)$ and $f_r(\sp y\ar 0\ssp,
\varphi\,,y\sp)\ssp(\eta)=f\sp(\sp y\ar 0\ssp,\varphi\,,y\sp)\ssp(\eta)=
y\ssp(\eta)-y\ar 0\sn-\int_{\,0}^{\,\eta}\varphi\circ[\,\roman{id\,},y\,]\ssp,
$ here put concisely but a bit imprecisely since $\dom{\rho_{\sp}}_{r\sp s}
\not=\dom{\rho_{\sp}}_r$ and $\dom f_r\not=\dom\sn f\sp$. The functions $
{\rho_{\sp}}_{r\sp s}$ are injective by the familiar uniqueness property of
holomorphic functions.

Let us study whether $\cal F$ is a {\it B\ssp}$\partial\ar 2\ssp$--\,extension 
of $(E\sqcap F\sp,F\sp,f\sp)$ in $\CinftyPi$ by $\cal R$ around\Newline $(\sp
\bnull E\ssp,\bnull F\ssp,\bnull F)\ssp$. From Example 3.9\ssp, we obtain $(E
\sqcap F_r\ssp,F_r\ssp,f_r\sp)\in\CinftyPi\ssp$. Conditions\newline (1) and
(2) of Definition 4.4 are satisfied since they reduce to proving that with $a=
\partial\ar 2\sp\varphi\circ[\,\roman{id\,},y\,]$ the differential equation $
(\sp v\ssp'-u\ssp'+a\cdot u\sp)\ssp(\eta)=0$ for $|\ssp\eta\ssp|<r$ with\Newline
initial condition $u\ssp(0)=v\ssp(0)$ has a unique solution $u\in\vecs F_r$
for all $v\in\vecs F_r\ssp$. The solution is $u=e^{\,A\sp}(\sp u\sp\ar 0+
\int_{\,0}(\sp e^{-A\sp}v\ssp'\sp))\ssp,$ where $A=\int_{\,0}a$ and $u\ar 0=
\bar{\roman D}\sp(r)\sn\times\sn\{\ssp u\ssp(0)\sp\}\ssp$.

We now study whether condition (3) holds. It reduces to the following {\ssp\it
problem\,}. Let $0<r<s<1\ssp,$ fix $(\sp\eta\ar 0\ssp,\varphi\sp)\in\vecs E
\ssp,$ and let ${y_{}}_n\in\vecs F_s$ and $y\yr 1\in\vecs F_r$ be such that
for $|\ssp\eta\ssp|<r$ we have the differential equations\vskip.5mm\centerline{$
{y_{}}_n\KN1'\sp(\eta)=(\sp 1-n^{-1})\,\varphi\ssp(\sp\eta\ssp,{y_{}}_n(\eta))
\,$ with $\,{y_{}}_n\sp(0)=(\sp 1-n^{-1})\,\eta\sp\ar 0$}

$\mhyppy{14}y\yr 1{}'(\eta)=\varphi\ssp(\sp\eta\ssp,y\yr 1(\eta))\mhyppy{16.5}
$ with $\,y\yr 1(0)=\eta\sp\ar 0$

\noin and also ${y_{}}_n\ssp|\,\bar{\roman D}\sp(r)\to y\yr 1$ uniformly on $
\bar{\roman D}\sp(r)\ssp$. Then we should be {\it able to extend\ssp} $y\yr 1$
to a vector of the space $H\ai b\sp(\sp\bar{\roman D\sp}(s))\ssp$. Here is a
counterexample: Take $\eta\sp\ar 0=1$ and de- fine $\varphi\ssp(\sp\eta\ssp,
\xi\sp)=s^{-1}\sp\xi^{\,2}$. Writing ${a_{}}_n=1-n^{-1},$ we have $
{y_{}}_n\sp(\eta)={a_{}}_n\sp s\,(\sp s-a_{\sp n}^{\,2\,}\eta\sp)^{-1}$ for $
|\sp\eta\sp|\le s\ssp,$ and $y\yr 1(\eta)=s\,(\sp s-\eta\sp)^{-1}$ for $
|\sp\eta\sp|\le r\ssp$. There is no $\bar y{}\yr 1$ with $
(\sp\bar y{}\yr 1\!,y\yr 1)\in{\rho_{\sp}}_{r\sp s}\ssp$.

In fact, there is no possibility to achieve $(E\ssp,F\sp,\Sigma\sp)\in\CinftyPi
$ for the solution rela- tion $\Sigma=f\sp\inve\sp[\sp\{\ssp\bnull F\}\sp]$ of
our present equation, since $\bnull E\in\dom\Sigma$ is not an interior point,
i.e., we have $\bnull E\not\in\roman{Int_{}}_{\tau_{rd}E}\ssp(\dom\Sigma\sp)
\ssp$. Indeed, since $\lim_{\,n\to\infty}n^{\,2\,}r^{\,n}=0$ for $0<r<1\ssp$,
given any $U\in\ymp E\ssp,$ choosing $\eps>0$ sufficiently small, and taking $
n\in\Ze$ sufficiently large with $n\,\eps>1\ssp,$ and defining $\varphi\ssp(\sp
\eta\ssp,\xi\sp)=n\,(\sp n+1\sp)\,\eta^{\,n}\sp\xi^{\,2},$ we have $x=(\sp\eps
\ssp,\varphi\sp)\in U,$ but there is no $y$ with $(\sp x\ssp,y\sp)\in\Sigma\,.
$ To see this, we only need to observe that the solution of $y\ssp'(\eta)=n\,(\sp
n+1\sp)\,\eta^{\,n}\sp(\sp y\ssp(\eta))^{\,2}$ with $y\ssp(0)=\eps$ is given
by $y\ssp(\eta)=\eps\,(\sp 1-n\,\eps\,\eta^{\,n\ssp+\sp 1\sp})^{-1},$ which is
not defined for all $|\sp\eta\sp|<1\ssp$.


\subhead6

                         Generalized well-posedness

Let us say that \hfil $\cal M\times\bit{034}$ \hfil is a {\it topologized
class of mappings\ssp} \hfil in case \hfil $\bit{034}$ \hfil is a function\linebreak
$\O\to\Univ$ with $\bit{034}\sp\value\ssn X$ a topology for each $X\in\O\ssp,$
and $\cal M$ is a class of triples $(X\sp,Y,f\sp)\ssp,$ where $X\sp,Y\in\O$
and $f\inc(\sp\bigcup\,(\sp\bit{034}\sp\value\ssn X\sp))\sn\times\sn(\sp
\bigcup\,(\sp\bit{034}\sp\value Y\sp))$ is a function. Now given such a class
$\cal M\times\bit{034}\ssp,$ let us consider a family $\bit{006}=\seq{\,
\widetilde\varSigma\ai A\sn:A\in\cal A\,}\ssp,$ where $\cal A$ is a filter
base on the set $\varOmega=\bigcup\sp\cal A\,$, and $\widetilde\varSigma\ai A=
(X\ar 0\ssp,Y\ssn\ai A\sp,\varSigma\ai A)\in\O^{\times2\!}\times\sn\Univ$ for $
A\in\cal A\,$. Also assume we are given the set $S\inc\bigcup\,(\sp
\bit{034}\sp\value\ssn X\ar 0)$ and the pair $P=(\varSigma\ar 0\ssp,Z\ar 0)\ssp
$, where $\varSigma\ar 0\inc\varSigma\ai{\char'012}$ and $Z\ar 0\in\O$ and \hfil
$\bigcup\,(\sp\bit{034}\sp\value Y\sn\!\ai{\char'012})\inc\bigcup\,(\sp
\bit{034}\sp\value\ssn Z\ar 0)$ \hfil and \hfil $(\sp\bit{034}\sp\value\ssn
Z\ar 0)\lei(\sp\bigcup\,(\sp\bit{034}\sp\value Y\sn\!\ai{\char'012}))\inc
\bit{034}\sp\value Y\sn\!\ai{\char'012}$ \hfil in \hfil case\linebreak $
\varOmega\in\cal A$ holds. We then consider the following {\it well-posed\,}ness
properties:\vskip1mm\noin
(1) \ $\all{x\in S}\,\exi{A\in\cal A\,,U\in\bit{034}\sp\value\ssn X\ar 0}\,x
    \in U\cap(\dom\varSigma\sn\ai A)$ and $(X\ar 0\ssp,Y\sn\!\ai A\sp,
    \varSigma\sn\ai A\ssp|\,U\sp)\in\cal M\ssp$,\par\noin
(2) \ $\all{A\in\cal A\setminus\sn\{\varOmega\}\ssp,x\in S}\,\exi{U\in
    \bit{034}\sp\value\ssn X\ar 0}$\newline$\mhyppy{43.8}x\in U\cap(\dom
    \varSigma\sn\ai A)$ and $(X\ar 0\ssp,Y\sn\!\ai A\ssp,
    \varSigma\sn\ai A\sp|\,U\sp)\in\cal M\ssp$,\par

\newpage 

    \noin
(3) \ $\varOmega\in\cal A$ and $\all{x\in S}\,\exi{U\in
    \bit{034}\sp\value\ssn X\ar 0}$\newline$\mhyppy{42.9}x\in U\cap(\dom
    \varSigma\ai{\char'012})$ and $(X\ar 0\ssp,Y\sn\!\ai{\char'012}\ssp,
    \varSigma\ai{\char'012}\ssp|\,U\sp)\in\cal M\ssp$,\par\noin
(4) \ $\varOmega\in\cal A$ and $\all{w\in\varSigma\ar 0}\,\exi W$\newline$\mhyppy{27.5}
    w\in W\in(\sp\bit{034}\sp\value\ssn X\ar 0))\rist\ssn(\sp
    \bit{034}\sp\value\sn Z\ar 0)$ and $(X\ar 0\ssp,Y\sn\!\ai{\char'012}\ssp,
    \varSigma\ai{\char'012}\sn\cap W\sp)\in\cal M\ssp$.\vskip1mm

We now propose the following terminology: Let us say that $\bit{006}$ is {\it
locally\ssp} well-posed {\it for\ssp} $S$ in $\cal M$ {\it by\ssp} $\bit{034}
\,$ if{}f (1) holds. Analogously, in case of (2) we speak of {\it almost
global\ssp} well-posedness, and we say that $\bit{006}$ is {\it globally\ssp}
well-posed for $S$ in $\cal M$ by $\bit{034}\,$ if{}f\linebreak (3) holds. In
the case of (4) we say that $\bit{006}$ is {\it almost well-posed\ssp} for $P$
in $\cal M$ {\it by\ssp} $\bit{034}\,$.

In case $S=\bigcup\,(\sp\bit{034}\sp\value\ssn X\ar 0)\ssp$, we drop "for $S
$\ssp" in the phrases just introduced associated with (1)\ssp, (2)\ssp,
(3)\ssp, \hfil and dropping "by $\bit{034}$\ " presupposes \hfil $
\bit{034}\ssp=\roman{pr}\ar 2\ssp|\,\roman{dom}\yr 2\cal M\ssp$. \hfil Well-\linebreak
posedness {\it for small data\ssp} shall refer to the case \hfil $S=
\{\sp\bnull X\sn\LHB{.3}{_{_0}}\}\ssp$, \hfil when $\O$ is a class of\linebreak
structured $\bold K\,$--\,vector spaces, i.e., the members of $\O$ are pairs \hfil
$X=(Z\sp,T\sp)$ \hfil where\linebreak $Z$ is a $\bold K\,$--\,vector space and
$T$ is the "structure\ssp" of $X\sp$.

\NSN{\bf
0} \Remarkk One observes that in the case of a singleton $\cal A=\{\varOmega\}
  \ssp$, global and local well-posedness are equivalent, and almost global
then always holds by $\cal A\setminus\sn\{\varOmega\}=$ $\emptyset\,$. Further,
almost global well-posedness of a family $\bit{006}=\seq{\,
\widetilde\varSigma\ai A\sn:A\in\cal A\,}$ holds if{}f the singleton families
$\{\ssp(A\ssp,\widetilde\varSigma\ai A)\ssp\}$ for $A\in\cal A\setminus\sn
\{\varOmega\}$ all are glob\sp/\ssp locally well-posed. Finally, in case $
\varOmega\in\cal A$ holds, global well-posedness of $\bit{006}$ is equivalent
to that of the singleton family $\{\ssp(\varOmega\ssp,\widetilde\varSigma_\varOmega)\ssp\}\ssp$.

\NSN{\bf
1} \Examplee The assertion about "improved local well-posedness in $H^{s_*}
  $\sp" contained in [\,KR\ssp; Thm.\ 1, p.\ 5\,] expressed in our terminology
is the following: Let $\cal M_s$ be the class of all triples $(X\sp,Y,f\sp)\ssp
,$ where $X\sp,Y$ are sets and $f\inc X\sn\times\sn Y$ is a function. Fixing $
s\in\Re$ with $3-\frac{\sqrt3}2<s\ssp,$ and the functions $\seq{\,
\gamma^{\,\iota_1\iota_2}:0\not=\iota\ar 1\sp,\iota\ar 2\in 4\,}$ and $\seq{\,
\varphi^{\,\iota_1\iota_2}:\iota\ar 1\sp,\iota\ar 2\in 4\,}$ satisfying
"suitable conditions\ssp" for a constant $\lambda\in\Rep,$ and taking $\cal A=
\{\,[\,0\,,\smb T\,]:\smb T\in\Rep\ssp\}\ssp,$ we define $\bit{006}=\seq{\,
(X\ar 0\ssp,Y\ssn\ai A\sp,\varSigma\ai A):A\in\cal A\,}$ and the set $S$ and
the family $\bit{034}\,$ as follows.
                                                
We take $X\ar 0=\vecs(H^{\sp s}(I\!\!R^{\,3\sp})\sqcap H^{\ssp s\sp-1\sp}(
I\!\!R^{\,3\sp}))\ssp,$ and for fixed $A=[\,0\,,\smb T\,]\in\cal A$ let $
Y\sn\!\ai A$ be the set of all continuous functions $y:[\,0\,,\smb T\,]\times
I\!\!R^{\,3}\to\Re$ such that the function $\bar y$ on $[\,0\,,\smb T\,]$
defined by $\bar y\value t\value\psi=\int_{\,I\!\!R^{\,3}}y\ssp(\sp t\ssp,
\eta\sp)\sn\cdot\sn\psi\ssp(\eta)\,d\ssp\eta$ for test functions $\psi\ssp,$
i.e., for $\psi\in\vecs\cal D\ssp(I\!\!R^{\,3\sp})\ssp,$ is "continuous into $
H^s$ and of class $C^1$ into $H^{s\sp-1\,}$".

Then we let $\varSigma\ai A$ be the set of all pairs $(\sp x\ssp,y\sp)=
(\sp x\ar 0\ssp,x\ar 1\sp,y\sp)\in X\ar 0\times Y\sn\!\ai A$ such that $
\bar y\ssp(0)=x\ar 0$ and $\bar y\ssp'(0)=x\ar 1\ssp,$ and such that $y$
satisfies\vskip.5mm\noin(?)\hfill " $\partial_0^{\,2}y
-\sum_{i,j=1}^{\,3}\gamma^{\,ij}\circ y\cdot\partial_i\partial_jy
=\sum_{i,j=0}^{\,3}\varphi^{\,ij}\circ y\cdot\partial_iy\cdot\partial_jy$ " \hfill\phantom{(?)}\vskip.5mm

\noin in a suitable sense. Since the authors of [\,KR\,] do not explicitly
define what this "suitable sense\ssp" should mean, we do not bother to make it
clear here either. One possible interpretation, however, is sketched in Remark
below.

We put $\bit{034}\ssp=\seq{\,\cal P\!_s\ssp Q:Q\in\Univ\,}\ssp,$ where $
\cal P\!_s\ssp Q=\{\,U:U\inc Q\,\}\ssp,$ the power set, i.e., the discrete
topology on a set $Q\ssp$. Finally, let $S$ be the set of all $x=(\sp
x\ar 0\ssp,x\ar 1)\in X\ar 0$ having "suitable\ssp" $H^{\sp s}\text{-}
H^{\ssp s\sp-1\,}$--\,norm $\le\lambda\ssp$.

With these "definitions\ssp" now [\,KR\ssp; Thm.\ 1\,] asserts $\bit{006}$ to
be locally well-posed for $S$ in $\cal M_s$ by $\bit{034}\,$. In addition to
this well-posedness result, there is also an assertion (obviously containing
a typographic mess) about some "Strichartz estimate\ssp", but this does not
interest us here.                                                         \NS

\Remarkk \hfil A \hfil precise \hfil meaning \hfil for \hfil (?) \hfil can \hfil
  be \hfil given \hfil along \hfil the \hfil following \hfil lines.\linebreak
Writing \hfil $\Omega={]}\,0\,,\smb T\,{[}\sn\times\sn I\!\!R^{\,3}$, \hfil
for $y\in Y\sn\!\ai A$ \hfil we interpret each term as a vector in the\linebreak
(strong topological) dual space $\cal D\ssp'(\Omega)$ of $\cal D\ssp(\Omega)
\ssp$. Observe that we have the "naturally\ssp" defined continuous linear
injections\vskip.5mm\centerline{$
C\ssp(A\ssp,H^{s-1})\embed C\ssp(A\ssp,H^{s-2})\embed C\ssp(A\ssp,L^2)\embed
L^1\KN{1.5}\LHB{.3}{\ar{loc}}(\Omega)\embed\cal D\ssp'(\Omega)\ssp$,}

\noin and the (continuous) bilinear maps "induced\ssp" by pointwise products

$\mhyppy{20}b_1:C\ssp(\Omega)\sqcap L^1\KN{1.5}\LHB{.3}{\ar{loc}}(\Omega)\to
L^1\KN{1.5}\LHB{.3}{\ar{loc}}(\Omega)\embed\cal D\ssp'(\Omega)\ssp$,\vskip.5mm\centerline{$
b_2:C\ssp(A\ssp,L^2)\sqcap C\ssp(A\ssp,L^2)\to C\ssp(A\ssp,L^1)\embed
L^1\KN{1.5}\LHB{.3}{\ar{loc}}(\Omega)\ssp$.}

With $z=y\,|\,\Omega\ssp$, we directly have $\gamma^{\,ij}\circ z\ssp,
\varphi^{\,ij}\circ z\in\vecs C\ssp(\Omega)\ssp$. \hfil The term $
\partial_0^{\,2}y$ has\linebreak a natural interpretation as a vector of $
\cal D\ssp'(\Omega)\ssp$, and as do $\partial_iy$ for $i=0\,,\ldots\sp 3\ssp$,
and $\partial_i\partial_jy$ for $i\ssp,\sp j=1\ssp,\ldots\sp 3$ in the space $
C\ssp(A\ssp,L^2)\ssp$. \hfil Then we interpret $\,
\gamma^{\,ij}\circ y\cdot\partial_i\partial_jy$\linebreak as $\,
b\ar 1(\sp\gamma^{\,ij}\circ z\ssp,\partial_i\partial_jy\sp)\ssp$, and $\,
\varphi^{\,ij}\circ y\cdot\partial_iy\cdot\partial_jy\,$ as $\,b\ar 1(\sp
\varphi^{\,ij}\circ z\ssp,b\ar 2\sp(\sp\partial_iy\ssp,\partial_jy\sp))\ssp$.

\NSN{\bf
2} \Examplee Let $E=H^\infty(\Re\sp)=\leLCS-\sup\,\{\ssp H^{\sp\sigma\!}:
  \sigma\in\N\,\}\ssp,$ where we consider $H^{\sp\sigma}$ as a (real) locally
convex space of continuous functions $\Re\to\Ce$ via the Sobolev embedding $
H^{\sp\sigma}(\Re\sp)_{\bold c}\embed C.(\Re\sp)_{\bold c}$ for $\frac12<
\sigma\in\Re\,$. \hfil Taking \hfil $\cal A=\{\,[\,0\,,\smb T\,]:\smb T\in\Rep
\ssp\}$ as before, we define $\bit{006}\ar 1=\seq{\,(E\ssp,F\ssn\ai A\ssp,
\varSigma\sn\ai A):A\in\cal A\,}$ as follows.

For fixed $I\in\cal A\,,$ let $S\ai I$ be the set of all pairs $(\sp x\ssp,y\sp
)\in(\sp\vecs E\sp)\sn\times\sn\Ce^{\,I\times I\!\!R}$ such that $\check y=
\seq{\ssp\seq{\,y\ssp(\sp t\ssp,\eta\sp):\eta\in\Re\,}:t\in I\,}\in\vecs
C^{\ssp 1\sp}(\sp I\sp,E\sp)\ssp$. Then we let $F\ssn\ai I\in
\roman{LCS}\ssp(\biit R\ssp)$ be the unique one such that $\vecs F\ssn\ai I\inc
\Ce^{\,I\times I\!\!R}$ with $y\mapsto\check y$ defining a linear
homeomorphism $F\ssn\ai I\to C^{\ssp 1\sp}(\sp I\sp,E\sp)\ssp$. \hfil
Arbitrarily \hfil fixing \hfil a \hfil smooth \hfil $\chi:\Repp\to\Re\,,$ \hfil
we \hfil put \hfil $S=\mathbreak\{\,x\in\vecs E:
                                                \sup\,(\sp\roman{rng\,}(\sp
|\ssp x\ssp|^{\,2}\cdot(\sp\chi\ssp'\sn\circ|\ssp x\ssp|^{\,2\sp})^{\,2\sp}))<
\frac12\,\}\ssp$.

For the nonlinear Schr\"dinger equation

\noin$(*)\mhyppy7
\imag\,\partial_1y+\partial_2^{\,2}y=y\cdot\chi\ssp'\sn\circ|\ssp y\ssp|^{\,2}
\cdot\partial_2^{\,2}(\chi'\circ|\ssp y\ssp|^2)$ \ with \ $y\ssp(\sp 0\,,\cdot)
=\check y\ssp(0)=x\ssp,$

\noin we finally let $\varSigma\ai I$ be set of all $(\sp x\ssp,y\sp)\in
S\ai I$ with $(*)\ssp$.

Letting $\cal M_c$ be the class of maps $(E\ssp,F\sp,f\sp)$ of (real) locally
convex spaces with $\dom f\in\taurd E\ssp,$ by [\,Po1\ssp; Cor.\ 7.6, p.\ 739\,]
now $\bit{006}\ar 1$ is locally well-posed for $S$ in $C^{\ssp 1}_c\,,$ and
almost globally well-posed for small data in $\cal M_c\,$. In [\,Po2\ssp; Thm.\
4.24, p.\ 167, Thm.\ 4.25, p.\ 168\,]\ssp, parallel result are obtained for
local well-posedness in $C^{\,0}_c\ssp,$ and almost global for small data in $
\cal M_c\,$.

We observe that the "general nature\ssp" of the well-posedness result of
Example 1 is considerably weaker than that of this example, since firstly the
class $\cal M_s$ there is only of maps between unstructured sets, whereas here
we have locally convex spaces and even $C_c^{\ssp 1}$ maps. Secondly, there $
\bit{034}\,$ gives the discrete topology on a set, but here it is the locally
convex one.

\NSN{\bf
3} \Examplee With $\smb N\in\N$ fixed, let $\Omega$ be a bounded open set in $
  I\!\!R\,\yi N,$ and writing $k=[\,\frac12\ssp\smb N\,]+2\ssp,$ let then $E=
H_0^{\,k\ssp+\sp 3\sp}(\Omega)\sqcap H_0^{\,k\ssp+\sp 3\sp}(\Omega)$ with the
Sobolev spaces again considered as spaces of continuous functions $\Omega\to
\Re$ via the Sobolev embedding theorem. Now taking $\cal A=\{\,[\,0\,,\smb T\,
]:\smb T\in\Rep\ssp\}\cup\{\ssp\Repp\ssp\}\ssp,$ we shall construct the family
$\bit{006}=\seq{\,(E\ssp,F\ssn\ai A\sp,\varSigma\ai A):A\in\cal A\,}$ as follows.

For each fixed $A\in\cal A\,,$ we let $F\ssn\ai A\in\roman{LCS}\ssp(\biit R\ssp)$
be the unique one having $\vecs F\ssn\ai A\inc C\ssp(A\sn\times\sn\Omega\sp)
\ssp,$ and such that $y\mapsto[\,\bar y\ssp,\bar y\ssp'\,]$ defines a
topological linear isomorphism of $F\ssn\ai A$ onto a topological linear
subspace of $C\ssp(A\ssp,H^{\sp 1}(\Omega)\sqcap H^{\ssp 0}(\Omega))$ with $y
\mapsto\check y\,|\,\sp{]}\,0\,,\smb T\sp\,{[}$\Newline defining a surjection \hfil $
\vecs F\ssn\ai A\to\bigcap\ssp\{\,\vecs H^{\sp i}(\ssp{]}\,0\,,\smb T\sp\,{[}
\,,H^{\,k\sp+\sp 3\sp-\sp i\sp}(\Omega)):i=1\ssp,2\ssp,3\,\}$ \hfil when-\linebreak
ever $0<\smb T\in A\ssp$. The Sobolev spaces in the latter requirement are
considered as\Newline spaces of continuous functions via the respective
Sobolev embeddings, whereas $H^{\sp 1}(\Omega)$ and $H^{\ssp 0}(\Omega)\cong
L^{\sp 2}(\Omega)$ denote Hilbert\sp(iz)\sp able spaces of distributions.

For $i=1\ssp,\ldots\,\smb N$ fixing $\smb A\ai i\in\Rep$ and $2\le r\ai i\in\N
\sp,$ and the functions $\smb F\sp,\sigma\ai i$ with the properties A$_1$ and
A$_3$ in [\,Z\ssp; p.\ 1049\,]\ssp, we let $\varSigma\sn\ai A$ be the set of
all $(\sp x\ssp,y\sp)=(\sp x\ar 0\ssp,x\ar 1\sp,y\sp)\in\vecs(E\sqcap
F\ssn\ai A)$ satisfying $y\ssp(\sp 0\,,\cdot)=x\ar 0$ and $
\partial\ar 0\ssp y\ssp(\sp 0\,,\cdot)=x\ar 1$ and\vskip.5mm\centerline{$
\wave y=\Delta\ssp\partial\ar 0\sp y+\smb F\circ(D\ssp y\sp)
+\sum_{i=1\,}^{\,N}{\partial_{}}_i\sp(\sp\sigma\ai i\sn\circ(\sp
\partial\ai i\sp y\sp)+\smb A\ai i\sp(\sp\partial\ai i\sp
\partial\ar 0\ssp y\sp)^{\,r_i})\,,$}

\noin where $D\ssp y=[\,y\ssp,\partial\ar 0\ssp y\ssp,\partial\ar 1\sp y\ssp,
\ldots\,\partial\ai N\sp y\ssp,\partial\ar 1\sp\partial\ar 0\ssp y\ssp,\ldots\,
\partial\ai N\sp\partial\ar 0\ssp y\,]\ssp$.

Finally letting $S\in\ymp E$ be determined by the description for $W\sn\ar 0
\cap W\ssn\ar 1$ in [\,Z\ssp; p.\ 1049\,]\ssp, now Theorem 2.3 in [\,Z\ssp; p.\
1050\,] asserts that $\bit{006}$ is globally well-posed for $S$ in $C^{\,0}_c
\ssp$. In particular, it asserts that the family $\bit{006}$ is globally
well-posed in $C^{\,0}_c$ for small data.\NS

\Remarkk Possibly one cannot allow general bounded open $\Omega$ in Example 3
above. In [\,Z\ssp; p.\ 1048\,]\ssp, one loosely speaks of a bounded domain
having "sufficiently smooth boundary\ssp". For precise requirents in order to
have $H^{\sp\sigma}(\Omega)\cong W^{\sp\sigma,\sp 2\sp}(\Omega)\ssp,$ see for
example [\,W\ssp; Thm.\ 5.3, p.\ 95, Thm.\ 5.4, pp.\ 95\,--\,96, Thm.\ 2.1,
p.\ 39\,]\ssp.

In addition, even local existence of a solution in [\,Z\,] remains problematic
since not a word is said to prevent $\liminf\,\seq{\,\smb T_n\sp}=0$ when $
[\,0\,,\smb T_n\sp{[}$ denotes the maximal (positive) interval of existence
for the system (5)\ssp, (6) there on p.\ 1051 of ordinary nonlinear second
order differential equations for $\seq{\,T_{jn}:j=1\ssp,\ldots\,n\,}\ssp$.

\NSN{\bf
4} \Examplee Let $\bit{006}\ar 2=\{\ssp(\sp I\ar 0\ssp,\tilde f\sp)\ssp\}\ssp
  $, where $\tilde f=(E\ssp,F\sp,\Sigma\sp)$ with $E\ssp,F\sp,\Sigma$ as in
Example 5.5\ssp, and $I\ar 0$ is any nonempty set. \hfil For example, \hfil we
may fix \hfil $I\ar 0=\Repp\ssp$. \hfil Further,\linebreak let \hfil $
\bit{006}\ar 3=\seq{\,\tilde f\sn\ai I\ssn:I\in\cal A\,}\ssp$, \hfil where \hfil
$\cal A$ \hfil is as in Examples 1 and 2 above, \hfil and we have\linebreak
$\tilde f\sn\ai I=(E\ssp,F\sn\ai I\sp,\varSigma\ai I)$ \hfil defined \hfil as \hfil
follows. \hfil We \hfil take \hfil $E=\Cinfty(\sp\roman S\yr 1\sn)\sqcap
\Cinfty(\Repp\sn\times\roman S\yr 1\sn\times\Re\sp)$\linebreak and $F\ssn\ai I
=\Cinfty(\sp I\times\roman S\yr 1)\ssp,$ and let $\varSigma\sn\ai I$ be the
set of $(\sp y\ar 0\ssp,\varphi\ssp,y\sp)\in\vecs(E\sqcap C^{\ssp 1}(\sp I
\times\roman S\yr 1))$ satisfying $\,\partial\ar 0\ssp y+\partial\ar 1\sp y=
\varphi\circ[\,\roman{id\,},y\,]\,$ with $\,y\ssp(\sp 0\,,\cdot)=y\ar 0\ssp$.

From Theorems 5.6 and 5.8\ssp, we now obtain the result that $\bit{006}\ar 2$
is globally well-posed in $\CinftyPi\ssp$, and $\bit{006}\ar 3$ is almost
globally well-posed for small data in $\CinftyPi\ssp$. We observe that the
"nature\ssp" of these well-posedness results is stronger than that of Example
2 in three respects. First, the class $\cal M$ here is smaller. Second, the
well-posedness here contains also a regularity and a wider uniqueness result
hidden in the definition of $\bit{006}_\iota\ssp$. Third, also the
nonlinearity $\varphi$ is allowed to vary here.

Observe that the above stated well-posedness of $\bit{006}\ar 3$ is not the
exact content of Theorem 5.8 but only a corollary of it. The exact content can
be reformulated as follows. \hfil {\it For any fixed $I\in\cal A\,$, take $E=
\Cinfty(\sp\roman S\yr 1\sn)\sqcap
\Cinfty(\sp I\sn\times\roman S\yr 1\ssn\times\Re\sp)\ssp$, let $F\sn\ai I$ and\linebreak
$\varSigma\ai I$ be as above, and put $S=\dom\varSigma\ai I\ssp$. Then the
singleton family $\{\ssp(\sp I\,;E\ssp,F\sn\ai I\sp,\varSigma\ai I)\ssp\}$ is
globally well-posed for $S$ in $\CinftyPi\ssp$}. \hfill Theorem 5.8 does not
give any information about $S$ but we trivially have $\{\sp\bnull E\}\inc S\ssp$.

\NSN{\bf
5} \Examplee \hfil With \hfil the \hfil notations \hfil of \hfil Example 5.3\ssp, \hfil
  let \hfil $\bit{006}=\{\ssp(\sp I\,;E\ssp,\Cinfty(I\sp)\ssp,\varSigma\sp)\ssp
\}$\linebreak and $P=(\varSigma\ar 0\ssp,C^{\,2\sp}(I\sp))\ssp$. The assertion
of Theorem 5.4 now is that $\bit{006}\sp$ is almost well-\Newline posed for $P
$ in $\CinftyPi\ssp$. Theorem 5.2 can be formulated similarly.

\RunMyHead{S.\ Hiltunen}{}{
}{}


\vskip5mm\centerline{\bf%
                             References}\vskip2mm
{\eightpoint\parskip.5mm\baselineskip3.5mm\leftskip9mm\parindent-9mm
   \def\item#1=#2..{[\,#1\,]\vskip-4mm\noindent#2.\vskip.5mm}
\item AS=Averbukh, V.\ I.\ \& Smolyanov, O.\ G.: The theory of differentiation
  in linear topological spaces, {\it Russian Math.\ Surveys\ssp} {\bf22}
  (1967) no 6, 201--258..
\item Du=Dudley, R.\ M.: {\it Real Analysis and Probability\ssp}, Wadsworth,
  Pacific Grove 1989..
\item FK=Fr\"licher, A.\ \& Kriegl, A.: {\it Linear Spaces and Differentiation
  Theory\ssp}, Wiley, Chichester 1988..
\item Ge=Geba, D.-A.: A local well-posedness result for the quasilinear wave
  equation in $\bold R^{2+1}$, {\it Comm.\ Partial Differential Equations\ssp}
  {\bf29} (2004) no 3-4, 323--360..
\item Hi1=Hiltunen, S.: Implicit functions from locally convex spaces to
  Banach spaces, {\it Studia Math\ssp}.\ {\bf134} (1999) no 3, 235--250..
\item Hi2=Hiltunen, S.: Mapping families, differentiation, and Lie groups of
  diffeomorphisms (article in preparation)\ssp..
\item Hi3=Hiltunen, S.: Banach space situations and the loss of derivatives
  (article in preparation)\ssp..
\item Hi4=Hiltunen, S.: Smooth solution maps for nonlinear wave equations
  (article in preparation)\ssp..
\item Ho=Horv\'ath, J.: {\it Topological Vector Spaces and Distributions\ssp},
  Addison--Wesley, Reading 1966..
\item Jr=Jarchow, H.: {\it Locally Convex Spaces\ssp}, Teubner, Stuttgart 1981..
\item KS=Keel, M.\ \& Smith, H.\ F.\ \& Sogge, C.\ D.: Almost global existence
  for quasilinear wave equations in three space dimensions, {\it J.\ Amer.\
  Math.\ Soc\ssp}.\ {\bf17} (2004) no 1, 109--153..
\item Ke=Keller, H.\ H.: {\it Differential Calculus in Locally Convex
  Spaces\ssp}, LNM {\bf417}, Springer, Berlin 1974..
\item Ky=Kelley, J.\ L.: {\it General Topology\ssp}, GTM {\bf27}, Springer,
  New York 1985..
\item KR=Klainerman, S.\ \& Rodnianski, I.: Improved local well-posedness for
  quasilinear wave equations in dimension three, {\it Duke Math.\ J\ssp}.\
  {\bf117} (2003) no 1, 1--124..
\item KM=Kriegl, A.\ \& Michor, P.\ W.: {\it The Convenient Setting of Global
  Analysis\ssp}, Amer.\ Math.\ Soc.\ Survey {\bf53}, Providence 1997..
\item Pi=Pizanelli, D.: Applications analytiques en dimension infinie, {\it
  Bull.\ Sci.\ Math\ssp}.\ {\bf96} (1972) 181--191..
\item PS=Ponce, G.\ \& Sideris, T.: Local regularity of nonlinear wave
  equations in three space dimensions, {\it Comm.\ Partial Differential
  Equations\ssp} {\bf18} (1993) no 1-2, 169--177..
\item Po1=Poppenberg, M.: Smooth solutions for a class of fully nonlinear
  Schr\"dinger type equations, {\it Nonlinear Analysis\ssp} {\bf45} (2001)
  723--741..
\item Po2=Poppenberg, M.: An inverse function theorem in Sobolev spaces and
  applications to quasi-linear Schr\"dinger equations, {\it J.\ Math.\ Anal.\
  Appl\ssp}.\ {\bf258} (2001) 146--170..
\item W=Wloka, J.: {\it Partial Differential Equations\ssp}, Cambridge Univ.\
  Press, Cambridge 1992..
\item Z=Zhijian, Y.: Initial boundary value problem for a class of non-linear
  strongly damped wave equations, {\it Math.\ Methods Appl.\ Sci\ssp}.\ {\bf26}
  (2003) no 12, 1047--1066..
\par}\vskip3mm

{\eightpoint\baselineskip3.5mm\parskip0mm
  Seppo Hiltunen\vskip.5mm

  Helsinki University of Technology\par
  Institute of Mathematics, U311\par
  P.O.\ Box 1100\par
  FIN-02015 HUT\vskip.5mm
  FINLAND\vskip0mm
  e-mail: shiltune\,\@\,cc.hut.fi}

\enddocument